\documentclass[11pt,a4paper]{article}
\usepackage[utf8]{inputenc}
\usepackage[T1]{fontenc}
\usepackage[english]{babel}
\usepackage{lmodern,microtype}
\usepackage[a4paper,margin=2.55cm,headheight=15pt,headsep=0.65cm]{geometry}
\usepackage{amsmath,amssymb,amsthm,mathtools,mathrsfs}
\usepackage{enumitem,booktabs,array}
\usepackage{tocloft}
\usepackage[hidelinks,unicode]{hyperref}
\usepackage{bookmark,fancyhdr}
\numberwithin{equation}{section}
\newtheorem{theorem}{Theorem}[section]
\newtheorem{lemma}[theorem]{Lemma}
\newtheorem{proposition}[theorem]{Proposition}
\newtheorem{corollary}[theorem]{Corollary}
\theoremstyle{definition}
\newtheorem{definition}[theorem]{Definition}

\theoremstyle{remark}
\newtheorem{remark}[theorem]{Remark}
\newtheorem{example}[theorem]{Example}
\newcommand{\R}{\mathbb R}\newcommand{\C}{\mathbb C}
\newcommand{\N}{\mathbb N}\newcommand{\Q}{\mathbb Q}
\newcommand{\ER}{\mathcal E_{\R}}\newcommand{\M}{\mathcal M}
\newcommand{\Fclass}{\mathscr F_n}\newcommand{\Hclass}{\mathscr H_n}
\newcommand{\norm}[1]{\left\lVert#1\right\rVert}

\newcommand{\dd}{\mathrm d}\newcommand{\out}{\mathsf{s}}
\newcommand{\ask}{\mathsf{c}}
\newcommand{\A}{\mathcal A}
\newcommand{\V}{\mathcal V}
\newcommand{\J}{\mathcal J}
\newcommand{\bad}{\bot}
\DeclareMathOperator{\diag}{diag}\DeclareMathOperator{\dist}{dist}

\hypersetup{
    pdftitle={Proof of the Cordero--Torregrosa Conjecture for Methods with Arbitrary Rules and Its Refinement under Paired Function and Jacobian Evaluations},
    pdfauthor={Rodrigo Alberto Castro Marín}
}

\title{
    \textbf{Proof of the Cordero--Torregrosa Conjecture\\
    for Methods with Arbitrary Rules and Its Refinement\\
    under Paired Function and Jacobian Evaluations}
}

\author{
    Rodrigo Alberto Castro Marín\\[0.5em]
    \small Institute of Mathematics,
    University of Valparaíso, Chile\\[0.3em]
    \small
    \href{mailto:rodrigo.castrom@uv.cl}{\texttt{rodrigo.castrom@uv.cl}}
    \qquad
    \href{https://orcid.org/0000-0003-0010-7775}{\texttt{ORCID: 0000-0003-0010-7775}}
}

\date{November 8, 2026}

\newcommand{\Estar}{\mathcal E_*}
\newcommand{\Pol}{\mathcal P}

\AtBeginDocument{}
\hypersetup{pdfsubject={Cordero--Torregrosa barrier, extension of query budgets, and paired refinement},pdfkeywords={Cordero--Torregrosa conjecture, universal local order, adaptive queries, entire functions, pairing, Jacobian}}

\begin{document}
\maketitle
\begin{abstract}
We study the Cordero--Torregrosa conjecture on the order that can be guaranteed when solving nonlinear systems with prescribed budgets of function and Jacobian evaluations. We consider deterministic, causal, memoryless methods whose rules may be discontinuous, nonmeasurable, or noncomputable. We give an independent proof of a scalar barrier for queries of $f$ and $f'$ by constructing a fixed entire function, a simple root, and a single sequence of initial points that simultaneously exclude every exponent greater than $2^{d-1}$, where $d$ is the maximum number of queries. An exact simulation of vector-valued responses yields the Cordero--Torregrosa inequality $p\le 2^{N_F+N_J-1}$ for function classes containing the analytic test fields. The proof holds for $N_F,N_J\in\N_0$, $N_F+N_J\ge1$, without assuming $N_J\le N_F$; in particular, it covers $N_J>N_F$. Under pairing at the end of the cycle, we obtain the sharper bound $p\le 2^{N_F-s+1}3^{s-1}$, $s=\min\{N_F,N_J\}\ge1$, without requiring temporal precedence. The proofs use indistinguishable perturbations, kernels of constraints, and a construction that preserves responses while producing a single entire witness. These results are upper bounds on the universal local order of a cycle; they do not establish attainability for every budget.
\end{abstract}
\clearpage
{\small\tableofcontents}
\clearpage

\section{Introduction}\label{sec:introduction}
In solving systems of nonlinear equations, we seek methods that reduce the error rapidly without requiring an excessive number of evaluations. The order of convergence describes part of this behavior, but does not by itself determine the information required by each iteration. Indeed, two methods may make the same number of queries and obtain different responses depending on the points at which they evaluate the function and its Jacobian. This leads us to ask what order can be guaranteed when the available information is fixed in advance.

For scalar equations, Kung and Traub \cite{KT1974} studied the relationship between the order of memoryless multipoint methods and their number of evaluations. For nonlinear systems, Arroyo, Cordero, and Torregrosa \cite{ACT2010,ACT2011} proposed a bound that distinguishes function queries from Jacobian queries. For a budget of $N_F$ function values and $N_J$ Jacobians, with $1\le N_J\le N_F$, the proposed barrier is
\begin{equation}\label{eq:intro-general}
 p\le 2^{N_F+N_J-1}.
\end{equation}
Interpreting this inequality requires specifying the problem class, the permitted operations, and the notion of order. We must also distinguish proving an upper bound from constructing, for every budget, a method that attains it. The earlier work on information and maximal order by Wo\'{z}niakowski, Wasilkowski, and Traub \cite{Woz1975,Woz1976,Wasilkowski1981,TraubWoz1980} shows the importance of stating these distinctions explicitly.

The problem becomes more difficult when arbitrary rules for processing responses are allowed. We consider deterministic, causal, memoryless procedures whose rules are fixed independently of the unknown field and its root. Evaluation points may be chosen from previous responses, without requiring them to remain near the root or to be separated from one another. Nor do we assume continuity, rationality, measurability, or computability of the rules. The budget and well-definedness may be required only in a neighborhood of each root, with constants and radii depending on the problem. This model is broader than that of finite rational programs and calls for arguments that do not rely on regularity of the method's decisions.

We study the \emph{universal local order of a cycle}: for each problem with a simple root, we require an estimate of the output error valid for all initial points sufficiently close to that root. The exponent is common to the function class, whereas the constants need not be. Hence, to rule out an order, it is not enough to construct a different function for each initial error that produces an unfavorable output. We need to identify a single problem for which the proposed estimate fails.

The proof of the general barrier begins with an independent scalar result for procedures making at most $d\ge1$ queries of $f$ or $f'$, without requiring both types of information to be obtained at the same points. We prove that there exist a fixed real entire function $f$, a simple root $a$, and a strictly decreasing sequence $s_m>a$, with $s_m\to a$, such that
\[
 \frac{|T_f(s_m)-a|}{(s_m-a)^p}\longrightarrow+\infty
 \qquad\text{for every }p>2^{d-1}.
\]
The function, root, and sequence are chosen before fixing $p$. Thus, the statement is stronger than ruling out each higher exponent separately. Moreover, the witness function can be kept close to the identity on a fixed complex disk while retaining a simple root in the interval under consideration.

To prove this assertion, we construct polynomial perturbations that preserve every response exactly. We use kernels of constraints, their Wronskians, and exceptional sets to control the effect of a perturbation on the root, even when function and derivative queries are unpaired. We then establish an obstruction at each scale, with constants chosen before the query nodes. Finally, we construct a sequence of polynomials converging on compact sets to an entire function that preserves all selected records. Equality of outputs will follow from equality of responses, without assuming continuity of the rules. We develop this argument in full, including the control of separations and passage to a single fixed witness.

Our first principal result is completed by transferring the scalar barrier to systems. We consider diagonal fields of the form
\[
 G_f(z_1,\ldots,z_n)=\bigl(f(z_1),z_2,\ldots,z_n\bigr)^T.
\]
In these fields, a query of the vector-valued function contains only one unknown scalar value, and a query of its Jacobian contains only one unknown first derivative. We show that every vector execution can be simulated while preserving the order, types, and number of queries, even when its nodes do not lie on the first coordinate axis. We thereby obtain \eqref{eq:intro-general} for classes containing our test fields. The second principal result extends the range of budgets: the inequality is proved for all integers $N_F,N_J\ge0$ with $N_F+N_J\ge1$, without imposing the historical restriction $1\le N_J\le N_F$. In particular, we cover $N_J>N_F$, where more Jacobian queries than function queries are permitted. This extends the scope of the inequality; it does not assert attainability.

Our third principal result improves the barrier when evaluations are paired. If every point at which the Jacobian is queried is also a function-query point by the end of the cycle, then, for $N_F,N_J\ge1$, we prove the bound
\begin{equation*}
 p\le 2^{N_F-s+1}3^{s-1},
 \qquad s=\min\{N_F,N_J\}.
\end{equation*}
This inequality is sharper than \eqref{eq:intro-general} when $N_J>1$. We first prove the result for causally paired records and then show that pairing at termination suffices: function values may be queried in advance and stored without prematurely delivering them to the original rules or increasing the local budgets. The case without Jacobian queries is treated separately, with the bound $p\le2^{N_F-1}$ for $N_F\ge1$. Neither barrier, by itself, implies the existence of a method attaining the corresponding exponent for every budget.

The scalar construction is closely related to our earlier work on query barriers and the Kung--Traub conjecture \cite{CastroKT}. Both studies use indistinguishable perturbations, kernels of constraints, and procedures that produce a fixed witness. In the present article, we give a complete proof for queries of orders zero and one. The earlier scalar result also yields an alternative derivation of the general barrier, but is not a premise of our proof. We explain this relationship and the correspondence between the models in a separate section, making clear both the methodological continuity and the logical independence of the new argument.

The function class on which order is required is also important. We show that, if problems are restricted to polynomial fields while completely arbitrary rules are retained, a noncomputable algebraic decoding permits a locally exact procedure with few queries. This does not contradict our barriers: their assumptions require a class containing the analytic fields used to construct the fixed witness. Therefore, each result must be read together with its assumptions on the rules, function class, and notion of order.

Finally, we examine repeated queries, the distinction between nominal budget and new information, and the relationship between the local order of a cycle and the orders of iteration trajectories. A barrier for the former should not automatically be interpreted as a barrier for every notion of order of an iterative sequence; we provide an additional sufficient condition for establishing this relationship.

The article is organized as follows. In Sections~\ref{sec:model} and \ref{sec:main-results} we specify the model and state the results. Sections~\ref{sec:scalar-construction}--\ref{sec:entire} contain the scalar construction and fixed entire witness, and in Section~\ref{sec:general} we prove the vector barrier. Section~\ref{sec:kung-traub} explains the relationship with Kung--Traub. Sections~\ref{sec:paired-model} and \ref{sec:terminal-pairing} develop the paired refinement and remove the temporal precedence requirement. In Section~\ref{sec:scope} we study examples, limitations, and consequences. We end with the conclusions.

\section{Problems, queries, and local order}\label{sec:model}
We specify the objects received by a method, the information it may query, and the meaning of a universal local guarantee. The definitions separate budgets imposed problem by problem from any uniformity across fields. We then introduce the total scalar representation and the two notions of pairing used in the main results.
\subsection{Problems and function classes}
We use $\N=\{1,2,\ldots\}$, $\N_0=\N\cup\{0\}$, and the norm $\|\cdot\|_\infty$ on $\R^n$. For an open set $\Omega\subset\R^n$ and $F\in C^1(\Omega,\R^n)$, a root $\alpha$ is simple if
\(F(\alpha)=0,\qquad \det DF(\alpha)\ne0.\)
A problem is a triple $(\Omega,F,\alpha)$ with these properties. The only variable initial datum is $x$; neither $\alpha$ nor the expression for $F$ is provided to the method. The domain only determines which queries are available.

We define $\ER=\{f:\C\to\C:\ f\text{ entire},\ f(\R)\subseteq\R\}$ and
\[
 \mathcal E_*=
 \left\{f\in\ER:
 f'(u)>\tfrac12\ (-1\le u\le1),\quad
 \text{$f$ has a root in $(-\tfrac12,\tfrac12)$}
 \right\}.
\]
Since $f'>1/2$ on $[-1,1]$, each $f\in\mathcal E_*$ has a unique root $a_f$ there.

Put $\iota(s)=(s,0,\ldots,0)^T$, $I_0\in\{\R,(-1,1)\}$, and $\Omega_0=I_0\times\R^{n-1}$. We define $\Hclass(I_0)$ by
\begin{equation}\label{eq:hardclass}
 \left(\Omega_0,G_f,\alpha_f\right),\qquad
 G_f(z)=\bigl(f(z_1),z_2,\ldots,z_n\bigr)^T,
 \quad \alpha_f=(a_f,0,\ldots,0)^T,
 \quad f\in\mathcal E_*.
\end{equation}
Each field is understood to be restricted to $\Omega_0$. The same formula extends to an entire field on $\C^n$. For $n=1$, we interpret $\R^{n-1}$ as the singleton consisting of the empty tuple.

\begin{definition}[Admissible function class]\label{def:function-class}
A class $\Fclass$ is admissible for our barriers if $\Hclass(I_0)\subseteq\Fclass$ for some $I_0\in\{\R,(-1,1)\}$.
\end{definition}

In particular, the class of all real-analytic fields on open sets, with their simple roots, satisfies this condition. So does the class of all fields of class $C^r$, for any integer $r\ge1$, or of class $C^\infty$, provided the domains in \eqref{eq:hardclass} are admitted. Transferring the barrier uses only $\Hclass(I_0)\subseteq\Fclass$, not analyticity of every $F\in\Fclass$.

\subsection{Data and rules}
We define the set of queries and their response spaces by
\[
 \mathcal Q_n=\{0,1\}\times\R^n,\qquad
 \mathcal R_0=\R^n,\qquad \mathcal R_1=\R^{n\times n}.
\]
For $y\in\Omega$, the responses are $E_F(0,y)=F(y)$ and $E_F(1,y)=DF(y)$.
For $y\notin\Omega$, $E_F(\tau,y)$ is undefined.

A record of length $k$ is
\[
 R_k=\bigl((\tau_i,y_i,A_i)\bigr)_{i=1}^k,
 \qquad \tau_i\in\{0,1\},\quad y_i\in\R^n,
 \quad A_i\in\mathcal R_{\tau_i}.
\]
Put $R_0=()$, $h_k=(x,R_k)$, and $\mathcal H_k=\{(x,R_k):x\in\R^n,\ R_k\text{ has the preceding format}\}$.

Let $\ask$ and $\out$ be two distinct labels. A decision belongs to the disjoint union
\[
 \mathcal A_n=
 \bigl(\{\ask\}\times\mathcal Q_n\bigr)
 \sqcup
 \bigl(\{\out\}\times\R^n\bigr).
\]
Decisions are $(\ask,(\tau,y))$ (a query) and $(\out,v)$ (output $v$).

\begin{definition}[Query method]\label{def:method}
A deterministic, causal, memoryless method is a family of partial rules
\[
 \Phi_k:D_k\longrightarrow\mathcal A_n,
 \qquad D_k\subset\mathcal H_k,\qquad k\in\N_0,
\]
fixed independently of the field and its root. We impose no continuity, rationality, measurability, or computability on $\Phi_k$ or $D_k$.
\end{definition}

The method is \emph{total with budget $d\ge1$} if $D_k=\mathcal H_k$ ($0\le k\le d$) and $\Phi_d(\mathcal H_d)\subseteq\{\out\}\times\R^n$. Only these levels are used in an execution. For fields defined on all of $\R^n$, every query has a response and the execution ends after at most $d$ queries. On a proper domain, queried nodes must still belong to the domain.

Determinism: $\Phi_k$ is single-valued; causality: $\Phi_k$ depends only on $h_k$; absence of memory: the cycle begins with $x$ and fixed constants, without responses from previous cycles.

\begin{definition}[Valid execution]\label{def:execution}
A valid execution from $x$ on $F$ is a record $R_m$ and an output $v$ such that, for $0\le k<m$,
\[
\begin{aligned}
 (x,R_k)&\in D_k,\quad \Phi_k(x,R_k)=(\ask,(\tau_{k+1},y_{k+1})),\quad y_{k+1}\in\Omega,\\
 &A_{k+1}=E_F(\tau_{k+1},y_{k+1}).
\end{aligned}
\]
Finally, $(x,R_m)\in D_m$ and $\Phi_m(x,R_m)=(\out,v)$.
We write $\M_F(x)=v$. If no such finite execution exists, the output is undefined.
\end{definition}

Single-valuedness of every $\Phi_k(h_k)$ and induction on $k$ give uniqueness of the record and output whenever they exist.

We define the counters
\(n_F(R_m)=\#\{i:\tau_i=0\},\qquad n_J(R_m)=\#\{i:\tau_i=1\}.\)
Every request is counted once, even if repeated; reusing a response does not increase $n_F$ or $n_J$. Internal operations are not counted.

\begin{definition}[Local budget]\label{def:local-budget}
We say that $\M$ has local budget $(N_F,N_J)$ on $\Fclass$ if, for every $(\Omega,F,\alpha)\in\Fclass$, there exists $\delta_{F,\alpha}>0$ such that, whenever
\(0<\norm{x-\alpha}_\infty<\delta_{F,\alpha},\)
the execution is valid and $n_F(R_m)\le N_F$, $n_J(R_m)\le N_J$. Outside these neighborhoods we impose no budget.
\end{definition}

\begin{definition}[Universal local order]\label{def:order}
The method has universal local order at least $p>1$ on $\Fclass$ if, for every problem in that class, there exist $C_{F,\alpha},\rho_{F,\alpha}>0$ such that the execution is valid and
\begin{equation}\label{eq:order}
 \norm{\M_F(x)-\alpha}_\infty
 \le C_{F,\alpha}\norm{x-\alpha}_\infty^p
 \qquad
 (0<\norm{x-\alpha}_\infty<\rho_{F,\alpha}).
\end{equation}
\end{definition}

The exponent $p$ is common; $C_{F,\alpha}$ and $\rho_{F,\alpha}$ may depend on the problem. Intersecting neighborhoods allows us to require the budget and \eqref{eq:order} simultaneously.

\begin{example}[A Newton cycle and its data]\label{ex:model-newton}
Let $F$ be of class $C^2$ near a simple root $\alpha$. Starting from $x$, we first query $V=F(x)$ and then $A=DF(x)$. If $A$ is invertible, we return $x-A^{-1}V$; otherwise we may return $x$. On the invertible branch,
\(R_2=((0,x,V),(1,x,A)).\)
The rules use only $x,V,A$ and are fixed before choosing $F$. Since $\det DF(\alpha)\ne0$, both queries are admissible near $\alpha$ and $(N_F,N_J)=(1,1)$.

Write $e=x-\alpha$. In a neighborhood of the root, the inverses are bounded and $DF$ is Lipschitz continuous. Moreover,
\[
F(x)-DF(x)e
=\int_0^1\bigl(DF(\alpha+se)-DF(x)\bigr)e\,\dd s
=O(\norm e_\infty^2).
\]
Multiplying by $DF(x)^{-1}$ gives $\M_F(x)-\alpha=O(\norm e_\infty^2)$. The constant may depend on $F$ and $\alpha$, as allowed by Definition~\ref{def:order}. A second query of $F(x)$ would increase the nominal counter, even though its response was already available. None of these data are retained from one cycle to the next.
\end{example}

\subsection{Preservation of responses and internal procedures}
\begin{lemma}[Indistinguishability]\label{lem:indistinguishability}
Suppose that an execution of $\M$ on $F$ from $x$ is valid. If another field $\widetilde F$ is defined at every node of this execution and gives the same responses there, then the execution on $\widetilde F$ from $x$ has exactly the same record and output.
\end{lemma}
\begin{proof}
The initial histories coincide. Suppose they coincide after $k$ queries. The rule $\Phi_k$ is the same in both problems and receives the same history. Hence it is defined in both and produces the same decision. If it orders termination, the outputs coincide. If it orders a query, its type and node coincide; by hypothesis, so does the response. Thus the histories at level $k+1$ are equal. Induction reaches termination of the original execution. We have used no continuity of any rule.
\end{proof}

\begin{remark}[Operations between queries]
The definition allows us to group into a single rule all internal operations performed before the next query or output. To justify this, fix a deterministic procedure whose internal steps do not access the field. Starting from a history, execute these steps up to the next request or termination. If this segment fails to terminate or an operation is undefined, leave the rule undefined on that history. Otherwise, the next decision depends only on the known data. Induction on the number of queries preserves the record and output. No uniform bound on the number of internal operations is required, although each valid execution must reach every decision and terminate.
\end{remark}

\subsection{Total scalar representation}
In the scalar case, $(z,r)\in\R\times\{0,1\}$ returns $f^{(r)}(z)$, with $f^{(0)}=f$. Repeated requests count as queries; reusing a response does not.

\begin{definition}[Total scalar method]\label{au:def:total}
Let $d\ge1$. A total scalar method with budget $d$ is a family of mappings
\[
\Psi_k:\R^{k+1}\longrightarrow
\bigl(\{\ask\}\times(\R\times\{0,1\})\bigr)
\sqcup\bigl(\{\out\}\times\R\bigr),\qquad 0\le k\le d,
\]
such that $\Psi_d$ takes only output values. The mappings are fixed independently of $f$ and its root. We impose no continuity, measurability, rationality, or computability on them.

With $h_0=(s)$, the rules act according to
\[
\begin{aligned}
\Psi_k(h_k)=(\ask,(z,r))&\Longrightarrow\ \eta_{k+1}=f^{(r)}(z),\quad h_{k+1}=(h_k,\eta_{k+1}),\\
\Psi_k(h_k)=(\out,q)&\Longrightarrow\ T_f(s)=q.
\end{aligned}
\]
\end{definition}

The procedure ends after $J\le d$ queries, because the rules are total, the functions are globally defined, and level $d$ allows only outputs. The record of an execution is the list of its query pairs $(z_i,r_i)$ and responses. Once the method is fixed, the input and earlier responses also determine the queried pairs. Therefore, these need not be included as independent coordinates of the history.

\begin{lemma}[Indistinguishability]\label{au:lem:indist}
Fix a method from Definition~\ref{au:def:total} and an input $s$. Suppose its execution on $f$ has record $((z_i,r_i))_{i=1}^J$. If another real entire function $g$ satisfies
\(g^{(r_i)}(z_i)=f^{(r_i)}(z_i)\qquad(1\le i\le J),\)
then both executions have the same query pairs, the same responses, and the same output. The same assertion holds for any prefix of the record.
\end{lemma}
\begin{proof}
The initial histories coincide. Suppose they coincide after $k$ responses. Applying the same rule $\Psi_k$ to this history yields the same decision. If it is an output, the outputs coincide. If it is a query in the record under consideration, its node and type coincide, and the hypothesis gives equality of responses. Thus the histories at level $k+1$ coincide. Induction preserves the prefix and, for the complete record, reaches the same output.
\end{proof}

\subsection{Pairing and temporal order of queries}
The restrictions we study also depend on which responses are obtained at the same node. We distinguish a condition on each prefix from one that is required only at the end of the cycle.
\begin{definition}[Causal pairing]\label{def:pairing}
A record $R=((\tau_i,y_i,A_i))_{i=1}^m$ is causally paired if
\[
\forall i\ (\tau_i=1\ \Longrightarrow\ \exists j<i:\ \tau_j=0,\ y_j=y_i).
\]
A method is locally paired if this property holds in the neighborhoods where well-definedness and the budget are required.
\end{definition}
Let
\(R=\bigl((\tau_i,y_i,A_i)\bigr)_{i=1}^{m}\)
be the record of a completed valid execution. We define the sets of distinct nodes
\begin{equation}\label{term:eq:sets}
 \V(R)=\{y_i:\tau_i=0\},\qquad
 \J(R)=\{y_i:\tau_i=1\}.
\end{equation}
Since $n_F(R),n_J(R)$ count requests, including repetitions,
\(|\V(R)|\le n_F(R),\qquad |\J(R)|\le n_J(R).\)

\begin{definition}[Pairing at the end of an execution]\label{term:def:terminal}
We say that the record $R$ is \emph{paired without temporal precedence}, or simply \emph{paired at termination}, if
\begin{equation}\label{term:eq:terminal}
 \J(R)\subseteq\V(R).
\end{equation}
The causal notion is that of Definition~\ref{def:pairing}.
\end{definition}

Causal pairing implies \eqref{term:eq:terminal}. The converse is false as a statement about records: the list $DF(x),F(x)$ satisfies \eqref{term:eq:terminal} but is not causally paired. Under the new condition, a prefix may contain Jacobian nodes at which function values have not yet been requested. The inclusion is required at the end of the same cycle, not in a later iteration.

\begin{definition}[Local condition]\label{term:def:local}
Let $\mathscr F$ be a class of problems with simple roots. The method $\M$ is \emph{locally paired at termination} with budget $(N_F,N_J)$ if, for every $(\Omega,F,\alpha)\in\mathscr F$, there exists $\delta_{F,\alpha}>0$ such that, for $0<\|x-\alpha\|_\infty<\delta_{F,\alpha}$, the execution is valid and
\(n_F(R)\le N_F,\quad n_J(R)\le N_J,\quad \J(R)\subseteq\V(R).\)
\end{definition}

\begin{example}[The same set of responses in two orders]\label{ex:pairing-records}
Let $y$ be a point computed from the available data. The record $F(x),DF(x),F(y),DF(y)$ is causally paired. By contrast, $DF(x),DF(y),F(x),F(y)$ is paired at termination, but its initial prefixes are not. If $y$ depends on $DF(x)$, the second record still defines a causal procedure: the data needed to choose $y$ have already been received. We must not confuse causality of the rules with functional precedence within each pair.

In both records, the final sets satisfy $\J(R)\subseteq\V(R)$. If $F(y)$ is omitted from the second and $y\ne x$, the inclusion fails. A response that could be reconstructed by a special identity of the field is not counted as an actual function query.
\end{example}

\section{Main results}\label{sec:main-results}
We state the three contributions organizing this work: an independent proof of the general barrier, its extension to budgets with no ordering restriction, and the refinement under pairing. The proofs follow the necessary tools: the scalar witness is constructed in Section~\ref{sec:entire}, the vector transfer in Section~\ref{sec:general}, and the terminally paired barrier in Section~\ref{sec:terminal-pairing}. In every case, we distinguish an upper bound from an attainability statement.
\subsection{The scalar obstruction}
\begin{theorem}[Scalar barrier with a fixed witness]\label{au:thm:main}\label{thm:scalar-input}
Let $T$ be a total scalar method with budget $d\ge1$. There exist $f\in\Estar$, its root $a=a_f$, and a strictly decreasing sequence $s_m>a$, with $s_m\to a$, such that
\begin{equation}\label{au:eq:strong}
\frac{|T_f(s_m)-a|}{(s_m-a)^p}\longrightarrow+\infty
\qquad\text{for every }p>2^{d-1}.
\end{equation}
The objects $f,a,(s_m)$ do not depend on $p$. Moreover, the function may be chosen so that
\[
\max_{|z|\le3}|f(z)-z|\le2^{-10},\qquad
1-2^{-10}\le f'(u)\le1+2^{-10}\quad(-1\le u\le1).
\]
\end{theorem}
The proof turns perturbations constructed at each scale into a single witness, preserving the order of quantifiers in the statement.
\subsection{The general vector barrier}
\begin{theorem}[Vector barrier for arbitrary rules]\label{thm:general}
Let $n\ge1$, $N_F,N_J\in\N_0$, and $d=N_F+N_J\ge1$. Let $\M$ be a deterministic, causal, memoryless method with local budget $(N_F,N_J)$ on $\Hclass(I_0)$. Then there exist $f\in\mathcal E_*$, a fixed field $G_f$, its root $\alpha=\iota(a_f)$, and a sequence $X_m=\iota(s_m)$, with $s_m>a_f$ and $s_m\to a_f$, such that the original executions are valid for sufficiently large $m$ and
\begin{equation}\label{eq:vector-hard}
 \frac{\norm{\M_{G_f}(X_m)-\alpha}_\infty}
      {\norm{X_m-\alpha}_\infty^p}
 \longrightarrow+\infty
 \qquad\text{for every }p>2^{N_F+N_J-1}.
\end{equation}
The field, root, and sequence are chosen before the final exponent $p$.
\end{theorem}
A local estimate with exponent $p$ on a class containing $\Hclass(I_0)$ bounds the quotient for every fixed problem. Therefore, this theorem implies $p\le2^{N_F+N_J-1}$. We formalize this implication in Corollary~\ref{cor:CT}, after proving the simulation.

\begin{remark}[Extension of the budget range]\label{rem:extended-budgets}
The historical range is $1\le N_J\le N_F$; Theorem~\ref{thm:general} allows $N_F,N_J\in\N_0$, $N_F+N_J\ge1$, including $N_J>N_F$. This extension does not establish attainability of the exponent for every budget. The sharper bound under pairing is established separately in the following theorem.
\end{remark}

\subsection{The barrier under pairing}

For $N_F\ge1$ and $N_J\ge0$, write $s=\min\{N_F,N_J\}$ and define
\begin{equation}\label{eq:paired-budget}
 B(N_F,N_J)=
 \begin{cases}
  2^{N_F-1},&s=0,\\
  2^{N_F-s+1}3^{s-1},&s\ge1.
 \end{cases}
\end{equation}
\begin{theorem}[Paired barrier without temporal precedence]\label{term:thm:extension}
Let $n\ge1$, $N_F\ge1$, and $N_J\ge0$. Let $\M$ be a deterministic, causal, memoryless method with local budget $(N_F,N_J)$ on $\Hclass(I_0)$. Suppose that, in the neighborhoods where the budget is required, the completed record of every execution satisfies \eqref{term:eq:terminal}, without imposing an order on the two types of evaluation at the same point.

Then there exist $f\in\mathcal E_*$, a fixed field $G_f$, its simple root $\alpha=\iota(a_f)$, and a sequence $X_m=\iota(s_m)$, with $s_m>a_f$ and $s_m\to a_f$, such that the original executions are valid for sufficiently large $m$ and
\begin{equation}\label{term:eq:extended}
 \frac{\norm{\M_{G_f}(X_m)-\alpha}_\infty}
      {\norm{X_m-\alpha}_\infty^p}\longrightarrow+\infty
 \qquad\text{for every }p>B(N_F,N_J).
\end{equation}
The objects in the conclusion can be chosen before fixing $p$.

In particular, if $\M$ has universal local order at least $p>1$ on a class containing $\Hclass(I_0)$, then
\begin{equation}\label{term:eq:final}
 p\le B(N_F,N_J).
\end{equation}
\end{theorem}
If $1\le N_J\le N_F$, $B=2^{N_F-N_J+1}3^{N_J-1}$; if $N_J\ge N_F\ge1$, $B=2\,3^{N_F-1}$. The improvement over the general barrier is strict for $N_J>1$. These claims are justified in Subsection~\ref{subsec:terminal-consequences}; they are not obtained by assuming rationality of the rules.

\begin{example}[Three budgets and two classes of records]\label{ex:main-budgets}
For $(N_F,N_J)=(2,1),(2,2),(3,2)$, the pairs (general, paired) are $(4,4),(8,6),(16,12)$. The difference within each pair comes from an assumption on the queried nodes, not from counting the same requests differently. If only $F(x)$, $DF(x)$, $F(y)$, and $DF(z)$ are known, with $z$ different from the function-query nodes, the paired bound does not apply to that record. None of these calculations alone proves that the corresponding bound is attainable.
\end{example}

\section{Scalar construction: kernels and Wronskians}\label{sec:scalar-construction}
We construct polynomial perturbations preserving the responses of a record without requiring pairing. The problem is to control their value at the root with constants independent of the nodes. For this purpose, we study the kernels of the constraints, the zeros of their Wronskians, and a normalization that yields uniform bounds on disks.
\subsection{Kernels of the constraints of a record}\label{au:sec:kernels}
We use $\Pol_n(\mathbb K)=\{Q\in\mathbb K[z]:\deg Q\le n\}$ for $\mathbb K\in\{\R,\C\}$, with $\deg0=-\infty$. Fix $P\in\R[z]\setminus\{0\}$.

\begin{definition}[Kernels of a record]\label{au:def:kernels}
Let $\mathcal Z_j=((z_i,r_i))_{i=1}^j\in(\R\times\{0,1\})^j$. Define
\(L_i(Q)=(PQ)^{(r_i)}(z_i).\)
For $n\ge0$, put
\[
V_{j,n}=\{Q\in\Pol_n(\R):L_i(Q)=0\text{ for }1\le i\le j\},
\qquad V_{0,n}=\Pol_n(\R).
\]
Its complexification is
\(V_{j,n}^{\C}=\{Q_1+\mathrm iQ_2:Q_1,Q_2\in V_{j,n}\}\subset\Pol_n(\C).\)
\end{definition}

Each $L_i$ is linear in $Q$; we do not require independence of the constraints. If $Q\in V_{j,n}$, then for every real polynomial $g$ and every $\eta\in\R$,
\((g+\eta PQ)^{(r_i)}(z_i)=g^{(r_i)}(z_i)\qquad(i\le j).\)
Therefore, adding $\eta PQ$ preserves the responses of the record. For $r_i=1$ we require only $(PQ)'(z_i)=0$, not $(PQ)(z_i)=0$.

\begin{lemma}[Minimal degree and effective indices]\label{au:lem:min-degree}
Let $\mathcal Z_j$ be a record. Define
\(n_* =\min\{n\ge0:V_{j,n}\ne\{0\}\}.\)
Then $n_*\le j$, $\dim V_{j,n_*}=1$, and $\deg Q=n_*$ for every $Q\in V_{j,n_*}\setminus\{0\}$. If
\(\sigma_i=n_*+1-\dim V_{i,n_*}\qquad(0\le i\le j),\)
then
\(\sigma_0=0,\quad\sigma_j=n_*,\quad \sigma_i-\sigma_{i-1}\in\{0,1\}.\)
Thus, $\#\{i:\sigma_i-\sigma_{i-1}=1\}=n_*$, and for its indices $i_1<\cdots<i_{n_*}$ we have $\sigma_{i_s}=s$.
\end{lemma}
\begin{proof}
By rank-nullity, $\dim V_{j,j}\ge\dim\Pol_j(\R)-j=(j+1)-j=1$; hence $n_*\le j$. By minimality, no nonzero element of $V_{j,n_*}$ has degree less than $n_*$.

If there were two linearly independent elements $Q_1,Q_2$ in this kernel, both would have degree $n_*$. A nontrivial linear combination canceling their leading coefficients would be a nonzero kernel element of degree less than $n_*$. This contradicts minimality. Consequently, the kernel has dimension one.

By rank-nullity, $\sigma_i$ is the rank of the map $Q\mapsto(L_1(Q),\ldots,L_i(Q))$ on $\Pol_{n_*}(\R)$. Adding one linear functional cannot decrease the rank and increases it by at most one. Its initial value is zero and its final value is $n_*+1-1=n_*$. The claims about the indices follow by counting these increments.
\end{proof}

\subsection{Wronskians and concentration of zeros}\label{au:sec:wronskian}
\subsubsection{Wronskians of subspaces}
\begin{definition}[Wronskian]\label{au:def:wronskian}
For $(q_1,\ldots,q_k)\in\Pol_n(\C)^k$, put
\[
W(q_1,\ldots,q_k)(z)=\det\bigl[q_i^{(r)}(z)\bigr]_{r=0,\ldots,k-1;\ i=1,\ldots,k}.
\]
If $(q_i)$ is a basis of $V$, we write $W_V=W(q_1,\ldots,q_k)$. A change of basis multiplies $W_V$ by a nonzero scalar; $Z(W_V)$ and $\deg W_V$ are intrinsic.
\end{definition}

\begin{lemma}[Degree and change of variables]\label{au:lem:wronskian}
Let $V\subset\Pol_n(\C)$ have dimension $k\ge1$. There is a basis with degrees $d_1<\cdots<d_k\le n$, and
\begin{equation}\label{au:eq:Wdegree}
W_V\not\equiv0,\qquad
\deg W_V=\sum_{i=1}^k d_i-\frac{k(k-1)}2\le k(n-k+1).
\end{equation}
Moreover, $W_V$ is constant if and only if $V=\Pol_{k-1}(\C)$. If $a\in\C$, $\delta>0$, and
\(V_{a,\delta}=\{q(a+\delta\,\cdot):q\in V\},\)
then, for the corresponding bases,
\begin{equation}\label{au:eq:Wscale}
W_{V_{a,\delta}}(w)=\delta^{k(k-1)/2}W_V(a+\delta w).
\end{equation}
\end{lemma}
\begin{proof}
The dimension of $V\cap\Pol_s(\C)$ increases by at most one when $s$ increases by one. We choose a polynomial at each of the $k$ levels where an increase occurs. Their degrees are distinct and the polynomials form a basis. Write their leading terms as $c_i z^{d_i}$, with $c_i\ne0$.

The degree of every term in the Wronskian determinant is at most $\sum_i d_i-\sum_{r=0}^{k-1}r$. The coefficient at this degree is
\[
\left(\prod_{i=1}^k c_i\right)
\det[(d_i)_r]_{r=0,\ldots,k-1;\ i=1,\ldots,k},
\]
where $(d)_0=1$ and $(d)_r=d(d-1)\cdots(d-r+1)$. Since $(d)_r$ is a monic polynomial of degree $r$ in $d$, row operations turn the last determinant into the Vandermonde determinant. Its value is
\(\prod_{i<i'}(d_{i'}-d_i)\ne0.\)
This proves the degree formula and nonvanishing. The inequality in \eqref{au:eq:Wdegree} follows from $d_i\le n-k+i$.

If $\deg W_V=0$, then $\sum_i d_i=k(k-1)/2$; since $d_i\ge i-1$, it follows that $d_i=i-1$ ($1\le i\le k$). The $k$ basis polynomials belong to $\Pol_{k-1}(\C)$ and span it. Conversely, the Wronskian of $1,z,\ldots,z^{k-1}$ is the nonzero constant $\prod_{r=0}^{k-1}r!$.

Finally, the derivative of order $r$ of $q_i(a+\delta w)$ is $\delta^r q_i^{(r)}(a+\delta w)$. Factoring $\delta^r$ from each row gives \eqref{au:eq:Wscale}.
\end{proof}

\subsubsection{A uniform bound over all subspaces}
Put $\overline D(a,\delta)=\{z\in\C:|z-a|\le\delta\}$ and count zeros with multiplicity.

\begin{lemma}[Zeros and limits of polynomials]\label{au:lem:zeros-limit}
Let $h_\ell,h\in\Pol_n(\C)$, with $h_\ell\to h$ coefficientwise and $h\not\equiv0$.

If every $h_\ell$ has at least $k$ zeros in $\overline D(0,1)$, then $h$ also has at least $k$ zeros there. If instead there exist $R_\ell\to+\infty$ such that $h_\ell$ has no zeros in $\overline D(0,R_\ell)$, then $h$ is a nonzero constant.
\end{lemma}
\begin{proof}
Coefficientwise convergence implies uniform convergence on compact sets. For the first assertion, choose $\varepsilon>0$ so that $h$ has no zeros in $1<|z|\le1+\varepsilon$. This is possible because a nonzero polynomial has finitely many zeros. On $|z|=1+\varepsilon$, the minimum of $|h|$ is positive. For sufficiently large $\ell$, $|h_\ell-h|<|h|$ on this circle. By Rouché's theorem, both polynomials have the same number of zeros in $D(0,1+\varepsilon)$. That number is at least $k$, and all zeros of $h$ in this disk belong to $\overline D(0,1)$.

For the second assertion, suppose $h$ has a zero $u$. Choose a small circle around $u$ on which $h$ does not vanish. By uniform convergence and Rouché's theorem, for sufficiently large $\ell$ the polynomial $h_\ell$ has a zero inside it. But this disk lies in $\overline D(0,R_\ell)$ for large $\ell$, a contradiction. Thus $h$ has no complex zeros, and the fundamental theorem of algebra implies that it is a nonzero constant.
\end{proof}

\begin{lemma}[Uniform concentration of zeros]\label{au:lem:concentration}
Given $n\ge0$ and $1\le k\le n+1$, there exists $C_{n,k}\ge1$ such that, for every subspace $V\subset\Pol_n(\C)$ of dimension $k$, every $a\in\C$, and every $\delta>0$, the following implication holds: if some $h\in V\setminus\{0\}$ has at least $k$ zeros in $\overline D(a,\delta)$, then $W_V$ has a zero in $\overline D(a,C_{n,k}\delta)$ and is nonconstant.
\end{lemma}
\begin{proof}
The change of variables $z=a+\delta w$ and \eqref{au:eq:Wscale} reduce the problem to $a=0$, $\delta=1$.

Suppose that no uniform constant exists in that case. For each integer $\ell\ge1$, choose a subspace $V_\ell$ of dimension $k$ and $h_\ell\in V_\ell\setminus\{0\}$ with at least $k$ zeros in $\overline D(0,1)$, such that $W_{V_\ell}$ has no zeros in $\overline D(0,\ell)$.

Identify $\Pol_n(\C)$ with $\C^{n+1}$ through coefficients and use the Euclidean coefficient norm. Choose an orthonormal basis $q_1^{(\ell)},\ldots,q_k^{(\ell)}$ of $V_\ell$. Normalize $h_\ell$ to have norm one and write
\(h_\ell=\sum_{i=1}^k c_i^{(\ell)}q_i^{(\ell)},\qquad \sum_{i=1}^k|c_i^{(\ell)}|^2=1.\)
By compactness of spheres in finite dimension, there is a subsequence along which all the coefficient vectors $q_i^{(\ell)}$ and $c^{(\ell)}$ converge. Denote their limits by $q_i$ and $c$. The $q_i$ remain orthonormal. Therefore,
\(V_\infty=\operatorname{span}_{\C}\{q_1,\ldots,q_k\}\)
has dimension $k$, and $h_\infty=\sum_i c_iq_i$ belongs to it and has norm one. In particular, $h_\infty\not\equiv0$. By Lemma~\ref{au:lem:zeros-limit}, it has at least $k$ zeros in $\overline D(0,1)$.

The coefficients of a Wronskian are polynomials in the coefficients of its arguments. Thus
\(W(q_1^{(\ell)},\ldots,q_k^{(\ell)})\longrightarrow W(q_1,\ldots,q_k)\)
coefficientwise. The limit is nonzero by Lemma~\ref{au:lem:wronskian}. The preceding Wronskians have no zeros in disks whose radii tend to infinity, even after passing to the subsequence. Lemma~\ref{au:lem:zeros-limit} implies that the limiting Wronskian is constant. By Lemma~\ref{au:lem:wronskian}, $V_\infty=\Pol_{k-1}(\C)$. This is impossible, since it contains the nonzero polynomial $h_\infty$ with at least $k$ zeros in a disk.

The contradiction yields $C_{n,k}$. Moreover, if the Wronskian of the original space were constant, that space would be $\Pol_{k-1}(\C)$ and could contain no nonzero polynomial with $k$ zeros. This also proves the final assertion.
\end{proof}

Having fixed the budget $d$, choose once and for all
\begin{equation*}
C_d=\max_{\substack{0\le n\le d\\1\le k\le n+1}}C_{n,k},\qquad
N_d=\sum_{n=0}^d\left\lfloor\frac{(n+1)^2}{4}\right\rfloor.
\end{equation*}
The constants $C_d,N_d$ depend only on $d$, not on $P$ or the nodes; we do not need an explicit expression for $C_d$.

\subsection{Exceptional sets and quantitative perturbations}\label{au:sec:perturbations}
\begin{definition}[Exceptional set of a record]\label{au:def:exceptional}
For a record $\mathcal Z_i$ and the fixed polynomial $P$, define
\begin{equation*}
Z_i=\bigcup_{\substack{0\le n\le d,\ V_{i,n}\ne\{0\}\\
W_{V_{i,n}^{\C}}\text{ nonconstant}}}
\{\zeta\in\C:W_{V_{i,n}^{\C}}(\zeta)=0\}.
\end{equation*}
Adoptamos $\dist(b,\varnothing)=+\infty$. 
\end{definition}

By \eqref{au:eq:Wdegree}, a Wronskian corresponding to a space of dimension $k$ has at most $k(n-k+1)\le\lfloor(n+1)^2/4\rfloor$ distinct zeros. Therefore,
\begin{equation*}
\#Z_i\le N_d.
\end{equation*}
It is important that $Z_i$ depends only on $P$ and the pairs $(z_1,r_1),\ldots,(z_i,r_i)$. It does not depend on the response values, the current polynomial, or future queries. Once those pairs are known, the set remains fixed throughout the stage.

\subsection{An exact example and its degeneracies}\label{au:sec:example}
The following example explains why we need the sets $Z_i$ and why we must avoid a rigid choice of kernel degree.

\begin{example}[One function query and one derivative query]\label{au:ex:two}
Let $P=1$, reference budget $d=3$, and consider the record consisting of the pairs for $f(x)$ and $f'(v)$, with $x\ne v$. The kernel of minimal degree is
\[
V_{2,2}=\operatorname{span}\{Q_2\},\qquad
Q_2(z)=(z-x)(z-(2v-x)).
\]
Indeed, $Q_2(x)=0$ and $Q_2'(v)=0$. In $\Pol_1(\R)$, $Q(x)=Q'(v)=0$ implies $Q=0$. In degree three, a basis of the kernel is
\(Q_2(z),\qquad R_3(z)=(z-x)(z-v)^2.\)
The polynomials are independent because their degrees differ, and the two independent constraints leave dimension two in $\Pol_3(\R)$. Computing the Wronskian gives
\begin{equation}\label{au:eq:example-W}
W(Q_2,R_3)(z)=(z-x)^2(z-v)(z-(3v-2x)).
\end{equation}
To verify the identity, put $u=z-x$ and $h=v-x$. Then $Q_2=u(u-2h)$ and $R_3=u(u-h)^2$. Therefore,
\begin{align*}
Q_2R_3'-Q_2'R_3
&=u(u-2h)(u-h)(3u-h)-2(u-h)u(u-h)^2\\
&=u(u-h)\bigl((u-2h)(3u-h)-2(u-h)^2\bigr)\\
&=u^2(u-h)(u-3h).
\end{align*}
Returning to $z$, we obtain \eqref{au:eq:example-W}. Consequently,
\(Z_2=\{x,v,2v-x,3v-2x\}.\)
The exceptional points are not limited to the queried nodes.
\end{example}

\begin{example}[A third function query]
Add the query $f(w)$. Initially suppose that $x,w,v$ are distinct and define
\(C=(v-x)(v-w),\qquad D=2v-x-w.\)
The polynomial
\begin{equation*}
\widetilde Q_3(z)=(z-x)(z-w)\bigl(C-D(z-v)\bigr)
\end{equation*}
satisfies all three constraints. Vanishing at $x,w$ is immediate. If $A(z)=(z-x)(z-w)$, then $A(v)=C$ and $A'(v)=D$, so that
\(\widetilde Q_3'(v)=A'(v)C-A(v)D=DC-CD=0.\)
For $D\ne0$, $\deg\widetilde Q_3=3$ and its third zero is $v+C/D$. For $D=0$ ($v=(x+w)/2$), $\widetilde Q_3=C(z-x)(z-w)$ and $\deg\widetilde Q_3=2$.

For $D\ne0$, the zero $v+C/D$ may belong to $Z_2$. For example, with $x=0$, $v=1$, $w=4$ is $5/2$, which does not belong to $\{0,1,2,3\}$. By contrast, with $x=0$, $v=1$, $w=3/2$ is zero, which does belong to that set. There is no universal exclusion.

If the nodes coincide, we work directly with the constraints. For example, $x=w\ne v$ repeats a function value and does not increase the rank; if $x=v=w$, the kernel requires a double zero at $x$. Proposition~\ref{au:prop:perturbation} includes these cases through minimal degree and effective indices, without dividing by $v-x$, $v-w$, or $D$.
\end{example}

\subsection{Zero counting and normalized perturbations}\label{au:sec:quantitative}
We now use separation from exceptional sets to control zeros of a kernel element. Centered normalization turns this count into a lower bound for the value of the perturbation at the root.
\begin{lemma}[Counting zeros from a separation]\label{au:lem:count}
Let $\rho>0$, $b\in\R$, $\dist(b,Z_i)>\rho$, and $V_{i,n}\ne\{0\}$, with $n\le d$. Every $Q\in V_{i,n}\setminus\{0\}$ has at most $\dim V_{i,n}-1$ zeros in $\overline D(b,\rho/C_d)$.
\end{lemma}
\begin{proof}
Write $k=\dim V_{i,n}=\dim_{\C}V_{i,n}^{\C}$. If its Wronskian is constant, the complexified space is $\Pol_{k-1}(\C)$ and no nonzero polynomial in it has more than $k-1$ zeros in the whole plane.

If the Wronskian is nonconstant and $Q$ had $k$ zeros in the indicated disk, Lemma~\ref{au:lem:concentration} would yield a Wronskian zero at distance at most
\(C_{n,k}\frac{\rho}{C_d}\le\rho\)
from $b$. That zero belongs to $Z_i$, contradicting the hypothesis.
\end{proof}

\begin{definition}[Centered normalization]\label{au:def:normalization}
For $Q\in\Pol_n(\C)\setminus\{0\}$ and $b\in\C$, write
\(Q(b+w)=\sum_{\nu=0}^n q_\nu w^\nu.\)
We say that $Q$ is normalized at $b$ when $\max_\nu|q_\nu|=1$. Dividing by the positive maximum of these absolute values always yields this normalization. If $Q,b$ are real, normalization preserves real coefficients.
\end{definition}

For $R>0$ and a continuous function on $\overline D(0,R)$, write
\(\norm h_R=\max_{|z|\le R}|h(z)|.\)
The coefficient norm in Subsection~\ref{au:sec:wronskian} is used only in the compactness proof; $\norm h_R$ denotes the uniform norm.

\begin{lemma}[Value of a normalized polynomial]\label{au:lem:normalized-value}
Let $Q\in\Pol_n(\C)$ be normalized at $b$, of degree $D$, and let $\zeta_1,\ldots,\zeta_D$ be its zeros, counted with multiplicity. Then
\begin{equation}\label{au:eq:normalized-lower}
|Q(b)|\ge2^{-n}\prod_{i=1}^D\min\{1,|b-\zeta_i|\}.
\end{equation}
If $|b|\le1/2$ and $R\ge1$, moreover,
\begin{equation}\label{au:eq:normalized-upper}
|Q(b)|\le1,\qquad \norm Q_R\le(n+1)(R+1)^n.
\end{equation}
\end{lemma}
\begin{proof}
For $D=0$, $|Q|=1$ and the empty product is one. Suppose $D\ge1$ and put $\eta_i=\zeta_i-b$. The factorization is
\[
Q(b+w)=q_D\prod_{i=1}^D(w-\eta_i).
\]
Each coefficient of the product consists of a sum of products of the $\eta_i$. The sum of the absolute values of all such products is bounded by $\prod_i(1+|\eta_i|)$. Consequently,
\[
1=\max_\nu|q_\nu|\le |q_D|\prod_{i=1}^D(1+|\eta_i|).
\]
Multiplying the resulting inequality for $|q_D|$ by $\prod_i|\eta_i|$, we obtain
\[
|Q(b)|\ge\prod_{i=1}^D\frac{|\eta_i|}{1+|\eta_i|}
\ge2^{-D}\prod_{i=1}^D\min\{1,|\eta_i|\}.
\]
The inequality $D\le n$ proves \eqref{au:eq:normalized-lower}. The argument also applies when some of the $\eta_i$ vanish, since then the lower bound is zero.

By normalization, $|q_0|\le1$ and $|q_\nu|\le1$. If $|z|\le R$ and $|b|\le1/2$, then
\[
|Q(z)|\le\sum_{\nu=0}^n|z-b|^\nu
\le(n+1)(R+1)^n,
\]
because $R+1\ge1$. This proves \eqref{au:eq:normalized-upper}.
\end{proof}

\begin{proposition}[Perturbation preserving a record]\label{au:prop:perturbation}
Let $0\le j\le d$. Suppose $b\in\R$, $|b|\le1/2$, $P(b)\ne0$, and
\(0<\rho_i\le1,\qquad \dist(b,Z_i)>\rho_i\qquad(1\le i\le j),\)
where $Z_i$ corresponds to the prefix of length $i$. There exists a real polynomial $Q$ of degree at most $j$, normalized at $b$, satisfying all constraints $L_i(Q)=0$ and
\begin{equation}\label{au:eq:product-lower}
|Q(b)|\ge(4C_d)^{-j}\prod_{i=1}^j\rho_i.
\end{equation}
In particular, $H=PQ$ preserves the responses of the record and satisfies
\begin{align*}
\norm H_R&\le\norm P_R(j+1)(R+1)^j\qquad(R\ge1),\\
|P(b)|(4C_d)^{-j}\prod_{i=1}^j\rho_i
&\le |H(b)|\le|P(b)|.
\end{align*}
\end{proposition}
\begin{proof}
For $j=0$, take $Q=1$. Let $j\ge1$. By Lemma~\ref{au:lem:min-degree}, the minimal-degree kernel $V_{j,n_*}$ is one-dimensional and has a real generator of degree $n_*\le j$. Choose that generator and normalize it at $b$. Denote the result by $Q$.

Let $i_1<\cdots<i_{n_*}$ be the effective indices from the same lemma. If $n_*=0$, the polynomial $Q$ is constant of modulus one, and \eqref{au:eq:product-lower} is immediate. Suppose $n_*\ge1$. For $1\le s\le n_*$, we have
\(\dim V_{i_s,n_*}=n_*+1-s,\qquad Q\in V_{i_s,n_*}.\)
Lemma~\ref{au:lem:count} shows that $Q$ has at most $n_*-s$ zeros in $\overline D(b,\rho_{i_s}/C_d)$.

Order its zeros so that
\(|\zeta_1-b|\le\cdots\le|\zeta_{n_*}-b|.\)
For $1\le m\le n_*$, put $s=n_*+1-m$. If $|\zeta_m-b|\le\rho_{i_s}/C_d$, the first $m=n_*-s+1$ zeros would lie in that disk, contrary to the preceding count. Hence
\[
|\zeta_m-b|>\frac{\rho_{i_{n_*+1-m}}}{C_d}\qquad(1\le m\le n_*).
\]
Applying $\rho_i\le1$, $C_d\ge1$, and Lemma~\ref{au:lem:normalized-value}, we obtain
\[
|Q(b)|\ge(2C_d)^{-n_*}\prod_{s=1}^{n_*}\rho_{i_s}
\ge(2C_d)^{-n_*}\prod_{i=1}^j\rho_i
\ge(4C_d)^{-j}\prod_{i=1}^j\rho_i.
\]
In the second inequality we add factors no greater than one; in the third we use $n_*\le j$. This proves \eqref{au:eq:product-lower}.

By construction $Q\in V_{j,n_*}$, so $H^{(r_i)}(z_i)=L_i(Q)=0$. The remaining bounds follow from \eqref{au:eq:normalized-upper}, $n_*\le j$, and $\norm{PQ}_R\le\norm P_R\norm Q_R$.
\end{proof}

\begin{remark}[Dependent constraints and confluent nodes]
The proof does not assume that the rank increases at every query. Ineffective indices contribute only factors $\rho_i\le1$ to the lower bound. Nor do we divide by differences between nodes. Thus repeated queries and coincident function and derivative nodes require no additional cases.
\end{remark}

\begin{example}[Effective constraints, repetition, and normalization]\label{ex:kernel-normalized}
Take $P=1$ and the record $f(0),f'(1),f(0)$. The third query repeats the first. The minimal degree is two, and the kernel is generated by $Q(z)=z(z-2)$. In $\Pol_2$ the successive ranks are $0,1,2,2$; the effective indices are one and two.

If the current root is $b=1/2$, then
\(Q(b+w)=w^2-w-\frac34.\)
Since $\max\{1,1,3/4\}=1$, $Q$ is already normalized at $b$. Its zeros are at distances $1/2$ and $3/2$ from $b$, and the bound in Lemma~\ref{au:lem:normalized-value} gives
\(|Q(b)|=\frac34\ge2^{-2}\min\{1,1/2\}\min\{1,3/2\}=\frac18.\)
On the other hand, $Q'(1)=0$ but $Q(1)=-1$. Hence adding $\eta Q$ preserves the derivative response at one without necessarily preserving the function value there, which was not queried. The example brings together the distinction between constraint types, minimal degree, effective indices, and centered normalization.
\end{example}

\section{Perturbations and the finite stage}\label{sec:finite-general}
We use the perturbations of the preceding section to separate roots from exceptional sets without changing previously received responses. We first fix a neighborhood of simple roots and then the amplitudes of all modifications. This order will provide an obstruction at every scale with constants chosen before the procedure's nodes are known.
\subsection{Simple roots and selection among perturbations}\label{au:sec:roots}
\begin{lemma}[Cauchy bound on the interval]\label{au:lem:cauchy}
If $h$ is holomorphic in a neighborhood of $\overline D(0,R)$ and $R\ge3$, then
\(|h'(u)|\le\norm h_R\qquad(-1\le u\le1).\)
\end{lemma}
\begin{proof}
The circle centered at $u$ with radius one is contained in $\overline D(0,2)$. Cauchy's integral formula gives
\[
h'(u)=\frac1{2\pi\mathrm i}\int_{|z-u|=1}\frac{h(z)}{(z-u)^2}\,dz.
\]
Its length is $2\pi$, and the denominator has modulus one. Taking absolute values gives the bound.
\end{proof}

\begin{lemma}[Uniform neighborhood with a simple root]\label{au:lem:root-neighborhood}
Let $f,P\in\R[z]$, $a\in(-1/2,1/2)$, $f(a)=0$, $f'(a)>0$, and $P(a)\ne0$. Given $R\ge3$ and $\varepsilon>0$, there exist $r,\delta,\mu,M>0$, with $M>1$ and $\delta<\varepsilon$, such that
\(I=[a-r,a+r]\subset(-1/2,1/2),\qquad \min_I|P|>0,\)
and every real polynomial $g$ with $\norm{g-f}_R<\delta$ has a unique root in $I$, lying in its interior, and satisfies
\begin{equation}\label{au:eq:root-derivative}
0<\mu\le g'(u)\le M\qquad(u\in I).
\end{equation}
\end{lemma}
\begin{proof}
By continuity, choose $r>0$ with $I\subset(-1/2,1/2)$, $\min_I f'>0$, and $\min_I|P|>0$. Define
\(\mu=\frac12\min_I f'>0,\qquad M=\max_I f'+1>1.\)
Take $0<\delta<\min\{\varepsilon,1,\mu,\mu r\}$. By Lemma~\ref{au:lem:cauchy}, $|g'-f'|<\delta$ on $I$, which implies \eqref{au:eq:root-derivative}.

Moreover, $f(a+r)\ge2\mu r$ and $f(a-r)\le-2\mu r$. The bound $|g-f|<\delta<\mu r$ preserves the signs at these endpoints. The intermediate value theorem yields an interior root. The inequality $g'\ge\mu$ shows that this root is unique and simple.
\end{proof}

\begin{lemma}[Two indistinguishable roots]\label{au:lem:two-roots}
Suppose $g$ and $g\pm\eta H$, with $\eta>0$, have roots $b,b_+,b_-$ in $I$ and all three functions satisfy \eqref{au:eq:root-derivative}. If $H(b)\ne0$, then $b_+$ and $b_-$ lie on opposite sides of $b$ and
\[
\frac{\eta|H(b)|}{M}\le|b_\pm-b|\le\frac{\eta|H(b)|}{\mu}.
\]
For every $q\in\R$, at least one of these two roots satisfies
\begin{equation*}
|b_\pm-q|\ge\frac{\eta|H(b)|}{M}.
\end{equation*}
\end{lemma}
\begin{proof}
Since $(g\pm\eta H)(b)=\pm\eta H(b)$, the mean value theorem gives
\[
b_\pm-b=\frac{\mp\eta H(b)}{(g\pm\eta H)'(\xi_\pm)}
\]
for intermediate points $\xi_\pm\in I$. Since $\mu\le(g\pm\eta H)'(\xi_\pm)\le M$, we obtain the bounds and opposite signs. In particular,
\(|b_+-b_-|\ge\frac{2\eta|H(b)|}{M}.\)
The triangle inequality forces at least one of the distances to $q$ to be no less than half this separation.
\end{proof}

\begin{lemma}[Chain of roots avoiding a finite set]\label{au:lem:chain}
Use $f,I,\delta,\mu,M$ from Lemma~\ref{au:lem:root-neighborhood}. Let $N\ge1$, $K>0$, $\eta>0$, and real polynomials $g,H$ such that
\[
\norm H_R\le K,\quad \norm{g-f}_R+N\eta K<\delta,
\quad N\eta K\le\mu/2.
\]
Let $b$ be the root of $g$ in $I$, and suppose $H(b)\ne0$. For $0\le\ell\le N$, denote by $b_\ell$ the root in $I$ of $g+\ell\eta H$. Then
\begin{equation}\label{au:eq:chain-upper}
|b_\ell-b|\le\frac{N\eta|H(b)|}{\mu},
\end{equation}
the sequence $b_0,\ldots,b_N$ is strictly monotone and
\begin{equation}\label{au:eq:chain-gap}
|b_{\ell+1}-b_\ell|\ge\frac{\eta|H(b)|}{2M}.
\end{equation}
If $Z\subset\C$ has at most $N-1$ points, there exists $\ell\in\{1,\ldots,N\}$ such that
\begin{equation}\label{au:eq:chain-avoid}
\dist(b_\ell,Z)\ge\frac{\eta|H(b)|}{4M}.
\end{equation}
\end{lemma}
\begin{proof}
For $0\le\ell\le N$, $\norm{g+\ell\eta H-f}_R<\delta$. By Lemma~\ref{au:lem:root-neighborhood}, their roots and derivative bounds exist. The mean value theorem between $b$ and $b_\ell$ gives
\[
|b_\ell-b|=\frac{\ell\eta|H(b)|}{(g+\ell\eta H)'(\xi)}
\le\frac{N\eta|H(b)|}{\mu}.
\]
By Lemma~\ref{au:lem:cauchy}, $|H'|\le K$ on $I$, so that
\[
|H(b_\ell)-H(b)|\le K|b_\ell-b|
\le\frac{KN\eta}{\mu}|H(b)|\le\frac12|H(b)|.
\]
Consequently, $H(b_\ell)$ has the sign of $H(b)$ and modulus at least $|H(b)|/2$. Now,
\((g+(\ell+1)\eta H)(b_\ell)=\eta H(b_\ell),\)
and another application of the mean value theorem shows that
\[
b_{\ell+1}-b_\ell
=-\frac{\eta H(b_\ell)}{(g+(\ell+1)\eta H)'(\xi_\ell)}.
\]
All increments have the same sign and modulus at least $\eta|H(b)|/(2M)$. This proves monotonicity and \eqref{au:eq:chain-gap}.

Put $\Delta=\eta|H(b)|/(2M)$. By monotonicity, distinct roots in the chain are separated by at least $\Delta$. Each point of $Z$ can lie at distance strictly less than $\Delta/2$ from at most one of the $N$ roots $b_1,\ldots,b_N$: otherwise the triangle inequality would give separation less than $\Delta$. Since $N-1$ has at most $Z$ points, one of those roots is at distance at least $\Delta/2$ from all of them. This is \eqref{au:eq:chain-avoid}.
\end{proof}

\subsection{Geometric amplitudes and preservation of separations}\label{au:sec:amplitudes}
Define
\begin{equation}\label{au:eq:exponents}
e_1=1,\qquad e_j=2^{j-2}\ (j\ge2),\qquad
E_0=0,\qquad E_j=\sum_{i=1}^j e_i=2^{j-1}\ (j\ge1).
\end{equation}
In particular, $e_{j+1}=\max\{1,E_j\}$ and $(e_j)$ is nondecreasing.

\begin{lemma}[Reserve for future displacements]\label{au:lem:tail-amplitudes}
Let $\lambda_j,\gamma_j>0$, $\mu>0$, and $N\ge1$ satisfy
\begin{equation}\label{au:eq:geometric-condition}
\lambda_{j+1}\le\frac{\lambda_j}{2},\qquad
\lambda_{j+1}\le\frac{\mu\gamma_j}{16N}.
\end{equation}
For $0<t<1$, suppose a finite sequence of real points $b_0,\ldots,b_L$ satisfies
\[
|b_j-b_{j-1}|\le\frac{N\lambda_j}{\mu}t^{e_j}
\qquad(1\le j\le L).
\]
Then, for $1\le i<\ell\le L$,
\begin{equation*}
\sum_{j=i+1}^{\ell}|b_j-b_{j-1}|
\le\frac{\gamma_i}{8}t^{e_i}.
\end{equation*}
If, in addition, $\dist(b_i,Z)\ge\gamma_i t^{e_i}$, it follows that
\begin{equation}\label{au:eq:preserved-distance}
\dist(b_\ell,Z)\ge\frac78\gamma_i t^{e_i}>
\frac12\gamma_i t^{e_i}.
\end{equation}
\end{lemma}
\begin{proof}
Since $e_j\ge e_i$ and $t<1$, we have $t^{e_j}\le t^{e_i}$ for $j>i$. The first inequality of \eqref{au:eq:geometric-condition} implies
\[
\sum_{j=i+1}^{\ell}\lambda_j
\le\lambda_{i+1}\sum_{r=0}^{\infty}2^{-r}
=2\lambda_{i+1}.
\]
Consequently,
\[
\sum_{j=i+1}^{\ell}|b_j-b_{j-1}|
\le\frac{2N\lambda_{i+1}}{\mu}t^{e_i}
\le\frac{\gamma_i}{8}t^{e_i}.
\]
Distance to a nonempty set is a 1-Lipschitz function. Applying this between $b_i$ and $b_\ell$ gives \eqref{au:eq:preserved-distance}. For the empty set, the statement is interpreted according to our convention and follows immediately.
\end{proof}

\begin{example}[Stage exponent and reserve for separations]\label{ex:stage-exponents}
For a budget $d=3$, the exponents for the three queries are $e_1=1$, $e_2=1$, and $e_3=2$, so that $E_3=4$. In the final modification, $e_4=4$. A separation created at the second query has magnitude at least $\gamma_2t$. Although the root is subsequently modified twice, Lemma~\ref{au:lem:tail-amplitudes} bounds the sum of these displacements by $\gamma_2t/8$ and preserves separation at least $7\gamma_2t/8$.

These exponents belong to the proof: they do not presuppose power expansions in $t$ for the method's nodes.
\end{example}
\subsection{The finite stage for unpaired queries}\label{au:sec:finite}
\begin{proposition}[Obstruction at a fixed scale]\label{au:prop:finite}
Let $T$ be a total scalar method with budget $d\ge1$. Let $f,P\in\R[z]$ and $a\in(-1/2,1/2)$, with
\(f(a)=0,\qquad f'(a)>0,\qquad P(a)\ne0.\)
Given $R\ge3$ and $\varepsilon>0$, there exist $c\in(0,1/8]$ and $\tau>0$ such that, for every $0<t<\tau$, with $x=a+t$, there exists a real polynomial $g$ with a simple root $b\in(-1/2,1/2)$ satisfying
\begin{align*}
\norm{g-f}_R&<\varepsilon,\qquad P\mid(g-f),\qquad
\deg(g-f)\le\deg P+d,\\
|b-a|&<t/8,\qquad |T_g(x)-b|\ge c\,t^{2^{d-1}}.
\end{align*}
Moreover, $P(b)\ne0$ and $b$ does not coincide with any node of the execution on $g$ from $x$.

The constants $c,\tau$ depend only on $f,P,a,R,\varepsilon,d$. In particular, they are chosen before $t$ and all execution nodes, and do not depend on the particular rules of $T$ beyond budget $d$.
\end{proposition}
\begin{proof}
\par\medskip\noindent\textit{Choice of constants}\par\nopagebreak
Apply Lemma~\ref{au:lem:root-neighborhood}. This fixes $I=[a-r,a+r]$, $\delta,\mu,M$, with $\delta<\varepsilon$ and $M>1$. Define
\[
p_-=\min_I|P|>0,\qquad
K=\norm P_R(d+1)(R+1)^d>0,\qquad N=N_d+2.
\]
Choose
\begin{equation}\label{au:eq:lambda-initial}
\lambda_1=\min\left\{1,\frac{\mu}{4NK},\frac{\mu}{32N},
\frac{\delta}{8NK}\right\}>0.
\end{equation}
Next, for $1\le j\le d+1$, define successively
\begin{equation*}
\kappa_j=p_-(4C_d)^{-(j-1)}\prod_{i=1}^{j-1}\frac{\gamma_i}{2},
\qquad
\gamma_j=\frac{\lambda_j\min\{\kappa_j,1\}}{4M},
\end{equation*}
and, when $j\le d$,
\begin{equation*}
\lambda_{j+1}=\min\left\{\frac{\lambda_j}{2},
\frac{\mu\gamma_j}{16N}\right\}.
\end{equation*}
The empty product is one. This recurrence is not circular: first $\lambda_j$ is known, then $\kappa_j$, using only earlier $\gamma_i$, then $\gamma_j$, and finally $\lambda_{j+1}$.

All constants are positive, and
\[
0<\gamma_j\le\frac{\lambda_j}{4M}\le\frac14,
\qquad
\sum_{j=1}^{d+1}\lambda_j\le2\lambda_1.
\]
Ponemos
\begin{equation*}
c=\min\left\{\frac18,
\min_{1\le j\le d+1}\frac{\lambda_j\min\{\kappa_j,1\}}{M}\right\}>0,
\qquad \tau=\min\{1/4,r/2\}>0.
\end{equation*}
The constants $c,\tau$ are fixed before $t$ and all adaptive queries.

\par\medskip\noindent\textit{Bounds guaranteeing admissibility}\par\nopagebreak
Fix $0<t<\tau$ and $x=a+t$. Begin with the polynomial $g_0=f$ and its root $b_0=a$. Every modification at step $j$ will have the form
\(g\longmapsto g+\ell\eta H, \qquad |\ell|\le N,\qquad 0<\eta\le\lambda_j,\)
where $H=PQ$, $\deg Q\le d$, and $Q$ is normalized at the current root. Proposition~\ref{au:prop:perturbation} supplies $\norm H_R\le K$.

For any choice of candidates in the first steps, the norm of the difference from $f$ is bounded by
\begin{equation}\label{au:eq:all-candidates}
NK\sum_j\lambda_j\le2NK\lambda_1\le\delta/4<\delta.
\end{equation}
This bound includes the candidate being considered before its selection. Consequently, all intermediate polynomials and all their candidates have a unique root in $I$ and satisfy \eqref{au:eq:root-derivative}. Likewise,
\(N\eta K\le N\lambda_1K\le\mu/4\le\mu/2.\)
The assumptions of the root lemmas are therefore fulfilled whenever $H$ does not vanish at the current root. This nonvanishing will follow from Proposition~\ref{au:prop:perturbation}.

Every increment is divisible by $P$ and has degree at most $\deg P+d$. The same holds for any sum of increments. This accounts for the algebraic properties that we shall maintain.

\par\medskip\noindent\textit{Invariants of the construction}\par\nopagebreak
After responding to $j$ queries we have a polynomial $g_j$, a root $b_j\in I$, and a record of length $j$. All responses in that record are exactly those of $g_j$. For each answered query define
\(\widehat Z_i=Z_i\cup\{z_i\},\qquad \#\widehat Z_i\le N_d+1=N-1.\)
When choosing $b_i$, require
\begin{equation}\label{au:eq:created-separation}
\dist(b_i,\widehat Z_i)\ge\gamma_i t^{e_i}.
\end{equation}
Every subsequent modification, including the final one, will satisfy
\begin{equation}\label{au:eq:step-displacement}
|b_j-b_{j-1}|\le\frac{N\lambda_j}{\mu}t^{e_j}.
\end{equation}
Lemma~\ref{au:lem:tail-amplitudes} then implies, for all previously constructed prefixes,
\begin{equation}\label{au:eq:prefix-separation}
\dist(b_j,\widehat Z_i)\ge\frac78\gamma_i t^{e_i}
>\frac12\gamma_i t^{e_i}\qquad(i\le j).
\end{equation}
For $i=j$, the bound follows directly from \eqref{au:eq:created-separation}. For $i<j$, only the sum of displacements of already constructed steps is used. We need not assume that future steps exist.

\par\medskip\noindent\textit{Before answering a query}\par\nopagebreak
Suppose a prefix of length $j-1$ has been constructed, where $1\le j\le d$. If the method orders termination, go to the final step. Otherwise, let $(z_j,r_j)$ be the requested query. This pair is known before responding: it depends only on $x$ and the responses already received.

From that pair, form the new record of pairs and conceptually compute $Z_j$ and $\widehat Z_j$. Definition~\ref{au:def:exceptional} shows that the new response need not be known in order to fix these sets.

Let $b=b_{j-1}$ be the current root. For $i<j$, define
\(\rho_i=\frac{\gamma_i}{2}t^{e_i}\in(0,1].\)
The bound \eqref{au:eq:prefix-separation} guarantees $\dist(b,Z_i)>\rho_i$. Since $b\in I$, we have $|b|<1/2$ and $P(b)\ne0$. Apply Proposition~\ref{au:prop:perturbation} to the prefix of length $j-1$ to obtain $H_j=PQ_{j-1}$, with $\deg Q_{j-1}\le j-1$, preserving all previous responses and satisfying
\begin{equation}\label{au:eq:H-lower-stage}
|H_j(b)|\ge p_-(4C_d)^{-(j-1)}
\prod_{i=1}^{j-1}\left(\frac{\gamma_i}{2}t^{e_i}\right)
=\kappa_j t^{E_{j-1}}>0.
\end{equation}
For $j=1$ we use $Q_0=1$ and the empty product.

Choose the amplitude
\begin{equation}\label{au:eq:amplitude}
\eta_j=\lambda_j\min\left\{1,\frac{t^{e_j}}{|H_j(b)|}\right\}>0.
\end{equation}
It follows that
\begin{equation}\label{au:eq:amplitude-product}
\lambda_j\min\{\kappa_j,1\}t^{e_j}
\le\eta_j|H_j(b)|\le\lambda_j t^{e_j}.
\end{equation}
Indeed, $\eta_j|H_j(b)|=\lambda_j\min\{|H_j(b)|,t^{e_j}\}$. For $j\ge2$, we have $E_{j-1}=e_j$ and apply \eqref{au:eq:H-lower-stage}. For $j=1$, $|H_1(b)|\ge\kappa_1$ and $t<1$, giving the same inequality.

Consider the $N$ candidates
\(g_{j-1}+\ell\eta_j H_j,\qquad 1\le\ell\le N.\)
All preserve the response prefix of length $j-1$. By Lemma~\ref{au:lem:indist}, they all produce the same next query $(z_j,r_j)$. Since $\#\widehat Z_j\le N-1$, Lemma~\ref{au:lem:chain} lets us choose a candidate whose root satisfies
\[
\dist(b_j,\widehat Z_j)
\ge\frac{\eta_j|H_j(b)|}{4M}
\ge\gamma_j t^{e_j}.
\]
Denote the chosen candidate by $g_j$ and now respond with the number $g_j^{(r_j)}(z_j)$. The record is valid for $g_j$: previous responses are preserved, and the new one is $g_j^{(r_j)}(z_j)$.

Finally, \eqref{au:eq:chain-upper} and \eqref{au:eq:amplitude-product} give
\[
|b_j-b_{j-1}|\le\frac{N\eta_j|H_j(b)|}{\mu}
\le\frac{N\lambda_j}{\mu}t^{e_j}.
\]
This verifies \eqref{au:eq:step-displacement}. Applying Lemma~\ref{au:lem:tail-amplitudes} to the sequence constructed thus far preserves previous separations. This completes the induction step.

\par\medskip\noindent\textit{After the last query}\par\nopagebreak
By totality and the budget, the method terminates after $J\le d$ queries. Denote by $q$ its output on the current polynomial $g_J$, and by $b_J$ the current root. We also allow $J=0$.

Apply Proposition~\ref{au:prop:perturbation} to the complete record, with the earlier separations, to obtain $H=PQ_J$ such that
\(\deg Q_J\le J\le d,\qquad |H(b_J)|\ge\kappa_{J+1}t^{E_J}>0.\)
Define $\eta$ by \eqref{au:eq:amplitude}, with $j=J+1$ and $H$ in place of $H_j$. The two polynomials $g_J\pm\eta H$ preserve the complete record and have the same output $q$.

By Lemma~\ref{au:lem:two-roots}, one of their roots, which we denote by $b_{J+1}$, satisfies
\begin{align*}
|q-b_{J+1}|
&\ge\frac{\eta|H(b_J)|}{M}\\
&\ge\frac{\lambda_{J+1}\min\{\kappa_{J+1},1\}}{M}
 t^{e_{J+1}}\ge c\,t^{e_{J+1}}.
\end{align*}
Moreover, its displacement satisfies \eqref{au:eq:step-displacement}, because the upper bound in the same lemma does not contain the factor $N$ and $N\ge1$. Denote the chosen polynomial by $g$ and its root by $b=b_{J+1}$.

By \eqref{au:eq:exponents}, $e_{J+1}\le2^{d-1}$; since $0<t<1$, we obtain
\(|T_g(x)-b|=|q-b|\ge c\,t^{2^{d-1}}.\)
If $J=0$, $e_{J+1}=e_1=1$.

\par\medskip\noindent\textit{Verification of all conclusions}\par\nopagebreak
By \eqref{au:eq:all-candidates}, $\norm{g-f}_R\le\delta/4<\varepsilon$. Divisibility and degree are preserved by adding the increments. The total displacement satisfies
\[
|b-a|\le\frac{N}{\mu}\sum_{j=1}^{J+1}\lambda_j t^{e_j}
\le\frac{2N\lambda_1}{\mu}t\le\frac{t}{16}<\frac t8,
\]
where we have used $e_j\ge1$ and \eqref{au:eq:lambda-initial}.

Since $b\in I$, we have $g'(b)\ge\mu>0$ and $P(b)\ne0$. By \eqref{au:eq:prefix-separation}, including the final modification,
\[
\dist(b,\widehat Z_i)\ge\frac78\gamma_i t^{e_i}>0
\qquad(1\le i\le J).
\]
Since $z_i\in\widehat Z_i$, the final root does not coincide with any queried node. Finally, every modification preserved all preceding responses, and the last one also preserved the output. Lemma~\ref{au:lem:indist} proves that the constructed record is exactly the record of the method's execution on $g$ from $x$.
\end{proof}

\begin{remark}[Order of choices and absence of circularity]\label{au:rem:choice}
The order is: data $f,P,a,R,\varepsilon,d$; constants $C_d,N_d,I,\delta,\mu,M,K$; finite sequences $\lambda_j,\kappa_j,\gamma_j$; constants $c,\tau$; scale $t$; and only then adaptive queries and perturbations. Preservation of a separation uses the upper displacement bound fixed by \eqref{au:eq:amplitude}, not its lower bound. The latter is used to obtain a new separation. The geometric reserve shows that both requirements are compatible without choosing constants after the nodes are known.
\end{remark}

\section{A single entire witness}\label{sec:entire}\label{au:sec:entire}
We prove the passage to the limit in a formulation that does not mention pairing. This allows us to use the same construction for both the general barrier and another finite stage yielding a different exponent.

\begin{proposition}[Principle of exact preservation between stages]\label{au:prop:entire}
Let $T$ be a total scalar method with a finite budget. Let $B\ge1$. Suppose it has the following property: for all data $f,P,a,R,\varepsilon$ in Proposition~\ref{au:prop:finite}, there exist $c\in(0,1/8]$ and $\tau>0$ such that, for each $0<t<\tau$, there exists a real polynomial $g$ with a simple root $b\in(-1/2,1/2)$ satisfying
\begin{equation}\label{au:eq:abstract-stage}
\norm{g-f}_R<\varepsilon,\quad P\mid(g-f),\quad
|b-a|<t/8,\quad |T_g(a+t)-b|\ge c t^B,
\end{equation}
and moreover $P(b)\ne0$ and $b$ does not coincide with any node of that execution.

Then there exist $f_\infty\in\Estar$, its root $a$, and a strictly decreasing sequence $x_m>a$, $x_m\to a$, such that
\[
\frac{|T_{f_\infty}(x_m)-a|}{(x_m-a)^p}\longrightarrow+\infty
\qquad\text{for every }p>B.
\]
These three objects are fixed before choosing $p$. Moreover,
\begin{equation}\label{au:eq:limit-near-id}
\norm{f_\infty-z}_3\le2^{-10},\qquad
1-2^{-10}\le f_\infty'(u)\le1+2^{-10}\quad(-1\le u\le1).
\end{equation}
\end{proposition}
\begin{proof}
\par\medskip\noindent\textit{Initial data and induction hypothesis}\par\nopagebreak
Begin with
\(f_0(z)=z,\qquad a_0=0,\qquad S_0=\varnothing.\)
$S_m$ will be the finite set of nodes queried up to stage $m$, without repetitions.

At stage $m$, we construct $f_m$, its root $a_m$, a scale $t_m$, the input $x_m$, the set $S_m$, and a margin $h_m>0$. We maintain $a_m\notin S_m$, that $f_m'>1/2$ on $[-1,1]$, and that $|a_m|<2^{-10}$. These conditions hold at the initial stage, interpreting the last strict inequality as $|a_0|=0$.

\par\medskip\noindent\textit{Choosing a new stage}\par\nopagebreak
Suppose the first $m-1$ stages have been constructed. Set
\begin{equation}\label{au:eq:eps-m}
R_m=m+3,\qquad
\varepsilon_m=2^{-m-10}\min\{1,h_1,\ldots,h_{m-1}\},
\end{equation}
where for $m=1$ the minimum is one. Define
\begin{equation*}
P_{m-1}(z)=\prod_{u\in S_{m-1}}(z-u)^2.
\end{equation*}
If $S_{m-1}$ is empty, the product is one. Since $a_{m-1}\notin S_{m-1}$, we have $P_{m-1}(a_{m-1})\ne0$.

Apply the property in the statement to $f_{m-1},P_{m-1},a_{m-1},R_m,\varepsilon_m$. We obtain $c_m,\tau_m>0$, with $c_m\le1/8$, before choosing the scale. Put
\(p_m=B+\frac1m.\)
For $m=1$ define $\vartheta_m=1$; for $m\ge2$, define $\vartheta_m=t_{m-1}/4$. Choose
\begin{equation*}
0<t_m<\min\left\{\tau_m,2^{-m},\frac14,\vartheta_m,
\left(\frac{c_m}{2^{p_m+1}m}\right)^{1/(p_m-B)}\right\}.
\end{equation*}
All terms in the minimum are positive. The last condition implies
\begin{equation}\label{au:eq:large-ratio}
\frac{c_m}{2^{p_m+1}}t_m^{B-p_m}>m.
\end{equation}
We require only $c_m>0$, without $\inf_m c_m>0$.

Put $x_m=a_{m-1}+t_m$. The finite stage provides $f_m$ and $a_m$ such that
\begin{align}
\norm{f_m-f_{m-1}}_{R_m}&<\varepsilon_m,
\qquad P_{m-1}\mid(f_m-f_{m-1}),\label{au:eq:fm-difference}\\
|a_m-a_{m-1}|&<t_m/8,
\qquad |T_{f_m}(x_m)-a_m|\ge A_m,
\qquad A_m=c_m t_m^B>0.\label{au:eq:fm-error}
\end{align}
$a_m\notin S_{m-1}$ and $a_m\ne u$ for every node $u$ of the current stage. Add these nodes to $S_{m-1}$ to form $S_m$. Then $S_m$ is finite, $a_m\notin S_m$, and we can choose the margin
\begin{equation}\label{au:eq:h-m}
h_m=\min\{1,A_m,t_m,\dist(a_m,S_m)\}>0.
\end{equation}
Our convention for distance to the empty set makes this formula valid even if all executions considered have terminated without queries.

\par\medskip\noindent\textit{The induction can continue}\par\nopagebreak
Since $R_j\ge3$ and $\varepsilon_j\le2^{-j-10}$,
\begin{equation}\label{au:eq:near-id}
\norm{f_m-z}_3
\le\sum_{j=1}^m\norm{f_j-f_{j-1}}_3
<\sum_{j=1}^m2^{-j-10}<2^{-10}.
\end{equation}
Applying Lemma~\ref{au:lem:cauchy} to $f_m-z$ gives
\(|f_m'(u)-1|<2^{-10}\qquad(-1\le u\le1).\)
In particular, $f_m'>1/2$ there. The finite stage places $a_m$ in $(-1/2,1/2)$; by $f_m(a_m)=0$ and \eqref{au:eq:near-id},
\(|a_m|=|a_m-f_m(a_m)|\le\norm{f_m-z}_3<2^{-10}.\)
Finally, $P_m(a_m)\ne0$ because $a_m\notin S_m$. All assumptions for the next stage have been verified.

\par\medskip\noindent\textit{Convergence on all complex disks}\par\nopagebreak
Fix $R>0$. When $m$ is sufficiently large that $R_{m+1}\ge R$, for every $\ell>m$ we have
\[
\norm{f_\ell-f_m}_R
\le\sum_{j=m+1}^{\ell}\norm{f_j-f_{j-1}}_{R_j}
<\sum_{j>m}\varepsilon_j
\le\sum_{j>m}2^{-j-10}.
\]
The final tail tends to zero. By completeness of the space of continuous functions on a compact disk, a uniform limit exists on each disk. These limits agree on intersections by uniqueness of pointwise limits. They define a function $f_\infty$ on $\C$.

We show that it is entire. Fix a circle and a point in its interior. Cauchy's integral formula for $f_m$ passes to the limit by uniform convergence on the circle. The resulting integral is holomorphic in the interior and agrees there with $f_\infty$. Since every point of the plane lies inside a circle of this kind, $f_\infty$ is entire. The same formula, differentiated under the integral on a smaller interior disk, gives uniform convergence of $f_m^{(r)}$ to $f_\infty^{(r)}$ on compact sets, for each fixed $r$.

$f_m(\R)\subseteq\R$ for every $m$; passing to the limit, $f_\infty\in\ER$. Passing to the limit in \eqref{au:eq:near-id} and the derivative bounds gives the estimates \eqref{au:eq:limit-near-id}.

\par\medskip\noindent\textit{Convergence of roots and tail control}\par\nopagebreak
The mean value theorem applied to $f_j$ between $a_{j-1}$ and $a_j$, using $f_j'\ge1/2$, gives
\[
|a_j-a_{j-1}|
\le2|f_j(a_{j-1})|
=2|(f_j-f_{j-1})(a_{j-1})|
<2\varepsilon_j.
\]
Thus $(a_m)$ is Cauchy; let $a=\lim_{m\to\infty}a_m$. For every $m\ge1$ and every $j>m$, the definition \eqref{au:eq:eps-m} includes $h_m$ in the minimum. Hence
\begin{equation}\label{au:eq:root-tail}
|a-a_m|\le2\sum_{j>m}\varepsilon_j
\le2h_m\sum_{j>m}2^{-j-10}
=2^{-m-9}h_m<\frac{h_m}{4}.
\end{equation}
We have $|a|\le2^{-10}<1/2$. Moreover,
\[
|f_\infty(a_m)|=|f_\infty(a_m)-f_m(a_m)|
\le\norm{f_\infty-f_m}_1\longrightarrow0.
\]
Since $a_m\to a$ and $f_\infty$ is continuous, $f_\infty(a)=0$. By \eqref{au:eq:limit-near-id}, the root is simple and $f_\infty\in\Estar$.

\par\medskip\noindent\textit{Exact preservation of each execution}\par\nopagebreak
Fix a stage $m$ and a node $u$ in its record. For every $j>m$, the factor $P_{j-1}$ contains $(z-u)^2$. Divisibility in \eqref{au:eq:fm-difference} implies
\(f_j(u)=f_{j-1}(u),\qquad f_j'(u)=f_{j-1}'(u).\)
By induction, these values are those of $f_m$ for every $j>m$. Now choose a disk containing this fixed node $u$. Convergence of the functions and their first derivatives on this disk allows passage to the limit, giving
\begin{equation}\label{au:eq:exact-limit}
f_\infty(u)=f_m(u),\qquad f_\infty'(u)=f_m'(u).
\end{equation}
We do not require $|u|\le R$ for a $R$ uniform across stages. We use only that, once the stage is fixed, each of its nodes is a finite real number and belongs to some disk on which convergence holds.

Applying Lemma~\ref{au:lem:indist} to the finite record of stage $m$ gives the exact equality
\begin{equation}\label{au:eq:output-preserved}
T_{f_\infty}(x_m)=T_{f_m}(x_m).
\end{equation}
Equality of the outputs follows from \eqref{au:eq:exact-limit} and Lemma~\ref{au:lem:indist}, without continuity of $T$.

Moreover, by \eqref{au:eq:root-tail} and \eqref{au:eq:h-m}, the limiting root belongs to no $S_m$: if $S_m\ne\varnothing$,
\[
\dist(a,S_m)\ge\dist(a_m,S_m)-|a-a_m|
>\frac34\dist(a_m,S_m)>0.
\]

\par\medskip\noindent\textit{Output error and initial error}\par\nopagebreak
From \eqref{au:eq:fm-error}, \eqref{au:eq:root-tail}, and \eqref{au:eq:output-preserved},
\begin{equation}\label{au:eq:limit-error}
|T_{f_\infty}(x_m)-a|
\ge |T_{f_m}(x_m)-a_m|-|a_m-a|
> A_m-\frac{A_m}{4}\ge\frac{A_m}{2}.
\end{equation}
On the other hand,
\(x_m-a=t_m+(a_{m-1}-a_m)+(a_m-a).\)
The last two differences have moduli less than $t_m/8$ and $t_m/4$, respectively. Consequently,
\begin{equation}\label{au:eq:input-error}
\frac58t_m<x_m-a<\frac{11}{8}t_m<2t_m.
\end{equation}
Thus $x_m>a$ and $x_m\to a$. The choice $t_{m+1}<t_m/4$ also implies
\[
x_{m+1}-a<\frac{11}{8}t_{m+1}
<\frac{11}{32}t_m<\frac58t_m<x_m-a.
\]
The sequence is strictly decreasing.

\par\medskip\noindent\textit{A single witness for all exponents}\par\nopagebreak
Using \eqref{au:eq:large-ratio}, \eqref{au:eq:limit-error}, and \eqref{au:eq:input-error}, we obtain
\[
\frac{|T_{f_\infty}(x_m)-a|}{(x_m-a)^{p_m}}
\ge\frac{c_m t_m^B}{2(2t_m)^{p_m}}
=\frac{c_m}{2^{p_m+1}}t_m^{B-p_m}>m.
\]
Fix now, and only now, an exponent $p>B$. For sufficiently large $m$ we have $p_m\le p$ and $0<x_m-a<1$. Therefore,
\[
\frac{|T_{f_\infty}(x_m)-a|}{(x_m-a)^p}
\ge\frac{|T_{f_\infty}(x_m)-a|}{(x_m-a)^{p_m}}>m.
\]
The quotient diverges to infinity. $f_\infty$, $a$, and $(x_m)$ were fixed before choosing $p>B$.
\end{proof}

\begin{proof}[Proof of Theorem~\ref{au:thm:main}]
Proposition~\ref{au:prop:finite} verifies all hypotheses of Proposition~\ref{au:prop:entire} with $B=2^{d-1}$. Apply the latter and rename $f_\infty$ as $f$ and $x_m$ as $s_m$.
\end{proof}

\begin{remark}[Why a scale-dependent family does not suffice]
Proposition~\ref{au:prop:finite} alone does not contradict a local estimate whose constants and radius depend on the polynomial. The contradiction with universal order comes after constructing $f_\infty$. For that single problem, the method's local constants would be fixed, whereas the quotients along the constructed sequence are unbounded. This is the decisive use of exact preservation between stages.
\end{remark}

\begin{example}[Preserving an execution does not mean approximating its rules]\label{ex:entire-preservation}
Suppose a stage queries $f(0)$ and $f'(1)$. To preserve this stage in all subsequent stages, we may choose
\(P(z)=z^2(z-1)^2.\)
Every modification $P(z)Q(z)$ preserves $f(0)$ and $f'(1)$; it also preserves $f'(0)$ and $f(1)$, although these were not queried responses. If the resulting functions converge together with their first derivatives on compact disks, both record responses remain exactly the same in the limit. By induction on the rules, the procedure reproduces the same execution from the same input, even if it makes discontinuous decisions from those responses.

The additional preservation is imposed between stages. By contrast, within the unpaired stage we impose only the actually queried constraints when constructing the kernels. These are two different roles of the perturbations and must not be confused in the count determining $2^{d-1}$.
\end{example}

\section{The general vector barrier}\label{sec:general}\label{sec:transfer}
We prove that full queries by a vector procedure on diagonal fields can be reproduced by a scalar query of the same type. First, we extend partial rules without altering locally valid executions. The scalar barrier already proved will then provide a fixed vector field that prevents an improvement of the general exponent.
\subsection{Extension of rules with budget control}
\begin{lemma}[Truncation and total extension]\label{lem:totalization}
Let $N_F,N_J\in\N_0$ and $d=N_F+N_J\ge1$. From the rules of $\M$ we construct a method $\widehat\M$ that, on fields defined on all of $\R^n$, terminates after at most $d$ queries and respects both counters. It preserves every valid execution of $\M$ respecting these counters.
\end{lemma}
\begin{proof}
For $k<d$, we preserve an original output decision whenever the rule is defined. We preserve a query decision only if the rule is defined and making that query does not exceed the corresponding counter. In every other case, we order termination with the zero vector. At level $d$, we preserve a defined original output and return the zero vector in all other cases. On histories formally incompatible with earlier decisions, we may also define the output to be zero.

The rules are total and independent of $F$; on global fields, $n_F\le N_F$, $n_J\le N_J$, and execution terminates no later than level $d$.

Now consider a valid original execution whose counters do not exceed the budgets. The original rule is defined at every prefix. No query exceeds its counter, and the final output is reached at a level no greater than $d$. By induction, every decision, response, and the output are preserved. Membership in the domains of the rules need not be algorithmically decidable.
\end{proof}

The extension also preserves a valid execution on a field defined on an open set if we compare it with an extension of the field agreeing at all its nodes. Only the responses at these nodes are used.

\subsection{Reconstruction of a vector response}
Recall the embedding $\iota(s)=(s,0,\ldots,0)^T$ from Section~\ref{sec:model}.
Given $Y\in\R^n$ and a real number $\eta$, define
\begin{equation}\label{eq:reconstruction}
 \mathcal R_0(Y,\eta)=(\eta,Y_2,\ldots,Y_n)^T,
 \qquad
 \mathcal R_1(Y,\eta)=\diag(\eta,1,\ldots,1).
\end{equation}
For every $f$, the reconstruction identities are
\(G_f(Y)=\mathcal R_0(Y,f(Y_1)),\qquad DG_f(Y)=\mathcal R_1(Y,f'(Y_1)).\)

\begin{theorem}[Exact simulation without additional queries]\label{thm:simulation}
Let $\M$ be a vector method with arbitrary rules. There exists a scalar method $T$ querying $f$ and $f'$ that reproduces the first component of the execution on $G_f$ from $\iota(s)$. For every valid execution,
\begin{equation}\label{eq:simulation}
 T_f(s)=\bigl(\M_{G_f}(\iota(s))\bigr)_1.
\end{equation}
Every vector function query corresponds to one query of $f$, and every Jacobian query to one query of $f'$. The counters are preserved exactly, even for repeated queries. The assertion holds both on $\R^n$ and on a cylinder $I\times\R^{n-1}$.
\end{theorem}
\begin{proof}
The initial scalar datum is $s$. From it we construct $\iota(s)$ and apply the first vector rule. If it orders termination with $v$, the scalar rule terminates with $v_1$. If it requests $(\tau,Y)$, we query $(Y_1,\tau)$ and receive $f^{(\tau)}(Y_1)$. Upon receiving $\eta$, we reconstruct the vector response by \eqref{eq:reconstruction} and append it to the vector record.

Suppose that, after $k$ scalar queries, we have reconstructed the vector history through that level. The next vector rule depends only on this history. Hence the next query type and node, or the output, can be determined from $s$ and the $k$ earlier scalar responses. Repeat the construction. If a vector rule is undefined, the corresponding scalar rule is also undefined.

This defines scalar rules depending only on the initial datum and previous responses. The vector nodes, their coordinates, and their types can be reconstructed by induction; they are not supplied as external data. On histories incompatible with this reconstruction, leave the rule undefined or, after totalization, return zero.

Applying these rules to a function $f$, we again prove by induction that every reconstructed response is the true response of $G_f$. For type $0$, every component except the first is a coordinate of the known node. For type $1$, every matrix entry except the first diagonal entry is a known zero or one. Thus only one unknown number is queried in either case. The entire execution is preserved and \eqref{eq:simulation} follows.

If $f$ is defined on $I$, then $Y\in I\times\R^{n-1}\iff Y_1\in I$. Thus query availability also coincides. No mechanism for answering queries outside the domain is used.

In particular, $n_F^{\rm esc}=n_F^{\rm vec}$ and $n_J^{\rm esc}=n_J^{\rm vec}$. We have not assumed that $Y_2=\cdots=Y_n=0$ at internal nodes or that the rules constructing them have any regularity.
\end{proof}

\begin{corollary}[Transfer of a local estimate]\label{cor:order-transfer}
If $\M$ satisfies \eqref{eq:order} on a diagonal problem with root $\alpha=\iota(a)$, then the induced scalar method satisfies
\(|T_f(s)-a|\le C_{G_f,\alpha}|s-a|^p\)
for $s$ sufficiently close to $a$. The local budget is also preserved.
\end{corollary}
\begin{proof}
We have $\norm{\iota(s)-\iota(a)}_\infty=|s-a|$ and
\(|T_f(s)-a| \le\norm{\M_{G_f}(\iota(s))-\iota(a)}_\infty.\)
Apply the vector estimate and shrink the neighborhood so that the budgets also hold.
\end{proof}

\subsection{Scalar rules on response vectors}
\begin{lemma}[Compatibility with the scalar query model]\label{lem:response-rules}
Let $\widehat\M$ be a total vector method with budget $d$, whose rules at level $d$ produce only outputs. There exist total mappings $\Psi_k$ ($0\le k\le d$), with the domain, codomain, and terminal condition of Definition~\ref{au:def:total}, defining a total scalar method. Only derivatives of orders zero and one are queried. For every real entire function $f$ and every $s\in\R$, the scalar method preserves the number and types of queries of $\widehat\M$ on $G_f$ from $\iota(s)$, and returns the first component of its output. The rules depend only on $\widehat\M$ and $d$, not on $f$ or its root.
\end{lemma}
\begin{proof}
At level $k$,
\(h_k=(s,\eta_1,\ldots,\eta_k)\in\R^{k+1}.\)
We must define the rules also on numerical vectors that occur in no execution. To do so, introduce a label $\bot$ distinct from those used earlier, and consider the disjoint union
\(\mathcal S_k=\mathcal H_k\sqcup(\{\bot\}\times\R^n).\)
An element of the second part indicates that reconstruction has already ended with a vector output, which we retain without using additional responses.

Define total mappings $\rho_k:\R^{k+1}\to\mathcal S_k$ by induction. At the initial level, put
\(\rho_0(s)=(\iota(s),()).\)
Suppose $\rho_k$ has been defined for $k<d$. For $\rho_k(h_k)=(\bot,v)$,
\(\rho_{k+1}(h_k,\eta)=(\bot,v)\qquad(\eta\in\R).\)
If $\rho_k(h_k)=(\iota(s),R_k)\in\mathcal H_k$, apply the total rule $\widehat\Phi_k$ to that history. When it returns $(\out,v)$, define $\rho_{k+1}(h_k,\eta)=(\bot,v)$ again. When it returns $(\ask,(\tau,Y))$, define
\[
 \rho_{k+1}(h_k,\eta)
 =\bigl(\iota(s),R_k\mathbin{\|}(\tau,Y,\mathcal R_\tau(Y,\eta))\bigr),
\]
where $R_k\mathbin{\|}(\tau,Y,A)$ denotes the list obtained by appending this triple to $R_k$. The mapping $\mathcal R_\tau$ is the one in \eqref{eq:reconstruction}. Its value belongs to the appropriate response space, whatever $\eta$ may be.

Now define $\Psi_k$. If $\rho_k(h_k)=(\bot,v)$, take $\Psi_k(h_k)=(\out,v_1)$. Otherwise compute $\widehat\Phi_k(\rho_k(h_k))$. An output $(\out,v)$ is transformed into $(\out,v_1)$, and a query $(\ask,(\tau,Y))$ is transformed into
\(\Psi_k(h_k)=(\ask,(Y_1,\tau)).\)
In the last expression, we write the node first and then the derivative order, as in Definition~\ref{au:def:total}. Since $\rho_k$ and $\widehat\Phi_k$ are total, so is $\Psi_k$. At level $d$ only an output can occur, whether reconstruction ended earlier or reached a history of level $d$.

Now fix $f\in\ER$ and $s\in\R$. Before the scalar execution terminates, we prove by induction that $\rho_k(h_k)$ coincides with the true vector history after $k$ queries. This holds for $k=0$. If the histories coincide and $(Y_1,\tau)$ is requested, the response is $\eta=f^{(\tau)}(Y_1)$. By \eqref{eq:reconstruction}, the appended triple contains exactly the response of $G_f$ or $DG_f$ at $Y$. Thus coincidence is preserved. If an output occurs, both executions end at that level and the scalar output is its first component.

States labeled with $\bot$ serve only to define the rules on response vectors that artificially prolong a terminated execution. They add no queries to any actual execution. We conclude that the counters, types, and output are preserved, and that the constructed mappings satisfy exactly the domains and terminal rule required by the scalar model.
\end{proof}

\begin{example}[Queries off the test axis]\label{ex:simulation-off-axis}
In dimension two, $G_f(u,v)=(f(u),v)^T$ and $DG_f(u,v)=\operatorname{diag}(f'(u),1)$. Suppose a rule first requests the value at $Y=(s+1,7)^T$ and, after receiving it, requests the Jacobian at $Z=(\theta,9)^T$, where $\theta$ is an arbitrary function of the available data. The simulation queries $f(s+1)$, reconstructs $(f(s+1),7)^T$, computes the same $\theta$, and queries $f'(\theta)$. The matrix is reconstructed as $\operatorname{diag}(f'(\theta),1)$.

The nodes do not lie on the first axis, and the rule defining $\theta$ need not be continuous. We use only that it receives the same reconstructed history. This example does not specify a convergent method; it shows preservation of the information and adaptivity required by Theorem~\ref{thm:simulation}.
\end{example}
\subsection{Proof of the barrier and its consequences}
\begin{proof}[Proof of Theorem~\ref{thm:general}]
Apply Lemma~\ref{lem:totalization} to the rules of $\M$. We obtain $\widehat\M$, which terminates within the stated budgets on every global field. Its rules do not depend on $I_0$ or $f$.

By Theorem~\ref{thm:simulation} and Lemma~\ref{lem:response-rules}, execution of $\widehat\M$ on the global fields $G_f$ induces a scalar method whose rules are defined on all of $\R^{k+1}$ for $0\le k\le d$, whose rule at level $d$ produces only outputs, and which makes at most $d$ queries. Its rules are independent of $f$; the only queried derivative orders are zero and one. All hypotheses of Theorem~\ref{thm:scalar-input} therefore hold. That theorem provides $f$, $a=a_f$, and $(s_m)$ satisfying \eqref{au:eq:strong}.

Now consider the problem in \eqref{eq:hardclass}, with domain $\Omega_0$ fixed in the statement. By the local budget assumption, for sufficiently large $m$ the original execution from $X_m=\iota(s_m)$ is valid, uses only nodes in $\Omega_0$, and respects both counters. Execution on the global field with the extended rules preserves exactly its responses: the extension and restriction agree at all these nodes. Lemma~\ref{lem:totalization} and induction by indistinguishability give
\(\M_{G_f}(X_m)=\widehat\M_{G_f}(X_m)\)
for all these indices. On the right we use the global field, and on the left its restriction when $I_0=(-1,1)$.

The first component on the right is $T_f(s_m)$. Moreover,
\[
 \norm{X_m-\alpha}_\infty=s_m-a,
 \qquad
 \norm{\M_{G_f}(X_m)-\alpha}_\infty\ge |T_f(s_m)-a|.
\]
Dividing by $(s_m-a)^p$ and using \eqref{au:eq:strong} gives \eqref{eq:vector-hard}. The index after which the executions coincide depends only on the problem and its budget neighborhood, not on $p$. Simultaneity from the scalar theorem is preserved.
\end{proof}
\begin{corollary}[Cordero--Torregrosa bound]\label{cor:CT}
Under the preceding hypotheses, if $\M$ has universal local order at least $p>1$ on a class $\Fclass\supset\Hclass(I_0)$, then the bound \eqref{eq:intro-general} holds.
\end{corollary}
\begin{proof}
If $p$ were larger, \eqref{eq:order} applied to the single problem constructed in Theorem~\ref{thm:general} would bound the quotient in \eqref{eq:vector-hard} by a fixed constant. This contradicts its divergence. The inclusion of classes ensures that this problem belongs to the class on which the estimate is required.
\end{proof}

\begin{remark}[Assumptions not used]
The proof does not prescribe the first node or the order of queries. It does not require starting with $F(x)$ and $DF(x)$, pairing the Jacobians, separating nodes, maintaining an exact query profile, or imposing a uniform bound on the number of internal operations. New information must nevertheless come exclusively from the full queries being counted. Neither external memory nor access to the expression of the field is allowed.
\end{remark}

\begin{remark}[The function class cannot be arbitrarily small]
The theorem uses entire test fields. With arbitrary rules, it does not establish a barrier on the class of polynomials of one fixed degree alone. Nor can high orders be forbidden on every subclass: for affine fields with invertible Jacobian, one Newton step solves the problem exactly. Generality of the rules and breadth of the function class are different assumptions.
\end{remark}

\begin{corollary}[Only a local scalar budget]\label{cor:scalar-local}
Let $T$ be a scalar procedure in the partial-rule model with local total budget $d\ge1$ on the problems of $\mathcal E_*$. Then there exist a fixed function and a sequence of inputs that simultaneously exclude all exponents greater than $2^{d-1}$.
\end{corollary}
\begin{proof}
Apply the construction of Lemma~\ref{lem:totalization} with a single total counter: before exceeding $d$, or when the original rule is undefined, return zero. The extension is total and queries only values or first derivatives. Theorem~\ref{thm:scalar-input} supplies a fixed witness. For that problem, the sequence eventually enters the neighborhood where the original execution is valid and uses at most $d$ queries. No substitution is then activated, so the outputs coincide. Divergence of the quotients is preserved.
\end{proof}

\section{Relationship with Kung--Traub and independence of the proof}\label{sec:kung-traub}
We clarify two distinct relationships with our manuscript on Kung--Traub \cite{CastroKT}. The first is a deduction: its scalar barrier, applied to the procedure just simulated, also yields the Cordero--Torregrosa inequality. The second is methodological: both constructions preserve records through kernels of constraints and produce a fixed entire witness. We explain why this continuity does not contradict the logical independence of the present proof.

\subsection{The deduction using the earlier scalar result}
The scalar barrier in \cite[Theorem 7.5 and Proposition 7.7]{CastroKT} allows total rules on response vectors and queries of derivatives of different orders. Its simultaneous version provides, for a budget $d\ge1$, a fixed real entire function, a fixed root, and a sequence of inputs excluding all exponents greater than $2^{d-1}$. We recall this assertion only to describe an alternative route. We have not used it in Sections~\ref{sec:scalar-construction}--\ref{sec:general}.

\begin{proposition}[Alternative deduction from the earlier scalar barrier]\label{prop:alternative-KT}
Assume the scalar result of \cite[Proposition 7.7]{CastroKT}, with its test-function class and hypotheses on the rules. Then the simulation in Section~\ref{sec:transfer} implies Theorem~\ref{thm:general} and the bound \eqref{eq:intro-general}.
This deduction does not require paired function and Jacobian queries or query nodes lying on the first coordinate axis.
\end{proposition}
\begin{proof}
Let $\M$, $N_F$, $N_J$, and $I_0$ be as in Theorem~\ref{thm:general}, and put $d=N_F+N_J$. Apply Lemma~\ref{lem:totalization} to obtain truncated total rules preserving the original valid executions within the budgets. By Theorem~\ref{thm:simulation}, a function query at $Y$ is simulated by querying $f(Y_1)$ and a Jacobian query by querying $f'(Y_1)$. The remaining components are reconstructed using \eqref{eq:reconstruction}.

Lemma~\ref{lem:response-rules} writes the induced decisions as total mappings of the initial point and earlier scalar responses. At level $d$, only an output is permitted. The rules are independent of $f$ and the root. The only derivative orders queried are zero and one, so the method belongs to the model of the assumed scalar result. No higher-derivative information is used to construct the simulation.

That result provides $f\in\mathcal E_*$, $a=a_f$, and $s_m>a$, with $s_m\to a$, for which
\[
\frac{|T_f(s_m)-a|}{(s_m-a)^p}\longrightarrow+\infty
\qquad(p>2^{d-1}).
\]
Once these objects are fixed, the local hypothesis on $\M$ implies that its original executions from $X_m=\iota(s_m)$ are valid and respect the budgets for sufficiently large $m$. No truncation is activated in those executions. If $I_0=(-1,1)$, the responses of the entire extension agree with those of the restriction at every node of these executions. Therefore,
\[
\bigl(\M_{G_f}(X_m)\bigr)_1=T_f(s_m),\qquad
\norm{X_m-\iota(a)}_\infty=s_m-a.
\]
The norm of the vector error is no less than the modulus of its first component. Scalar divergence therefore implies \eqref{eq:vector-hard}. If universal local order $p$ is guaranteed, the corresponding estimate on that fixed problem excludes $p>2^{d-1}$. Since $d=N_F+N_J$, we obtain the stated bound.
\end{proof}

\subsection{Methodological continuity and specialization to two query types}
\begin{remark}[Methodological origin and self-contained proof]\label{rem:methodological-continuity}
The proof in this article belongs to the same line of construction developed in \cite{CastroKT}. This relationship must be distinguished from using a theorem as a premise. We explain the correspondences and differences involved in the argument.

\textit{The constraints describe exactly the information received.}
If a record contains the pairs $(z_i,r_i)$, a modification $H$ preserves their responses when $H^{(r_i)}(z_i)=0$ for all indices. Writing $H=PQ$ separates preservation of earlier stages, handled by $P$, from constraints of the current stage, handled by $Q$. Here $r_i\in\{0,1\}$; in \cite{CastroKT}, higher orders are allowed. This restriction does not remove the difficulty of unpaired nodes: knowing $f'(v)$ is not equivalent to knowing $f(v)$.

\textit{The kernel is not replaced by a multiplicity that was not queried.}
The kernels are $V_{j,n}=\bigcap_{i=1}^j\ker(L_i|_{\Pol_n(\R)})$. Their minimal degree and the indices at which the rank increases also permit dependent constraints. We use no quotients $1/(z_i-z_j)$. In particular, coincident nodes and repeated queries are covered by the same construction. For the paired case, we later use additional structure: a value and a derivative at the same node can be preserved by a double factor. This simplification does not apply at a node where only the derivative was queried.

\textit{Exceptional sets arise from the perturbation space.}
Zeros of Wronskians are not limited to the method's nodes. Example~\ref{au:ex:two} shows that two queries can produce further exceptional points. Uniform concentration of zeros relates separation of the root from these sets to a lower bound on the value of a normalized kernel element. This step turns a geometric condition into a quantitative estimate for the root displacement.

\textit{We preserve exact prefixes without regularizing the rules.}
Before answering a query, its node and type are determined by previous responses. The corresponding exceptional set can be fixed at that point because it depends on the queried pairs, not on the new response. We then choose a perturbation preserving the prefix and respond only afterward. The procedure may change discontinuously as the scale varies; the proof needs neither a continuous branch for these changes nor an asymptotic expansion.

\textit{The amplitude reserve controls the entire continuation of a stage.}
In this presentation, we use Lemma~\ref{au:lem:tail-amplitudes} to choose the amplitudes before the scale. The inequality
\(\sum_{j=i+1}^{\ell}|b_j-b_{j-1}| \le\frac{\gamma_i}{8}t^{e_i}\)
reserves a fixed fraction of each separation for all later steps, including the final modification. It is enough that $e_{j+1}\ge e_j$. This choice lets us prove preservation of separations by a single tail estimate. It does not change the admissible information types or add any regularity assumption on the rules.

\textit{The conclusion is obtained on a fixed problem.}
The finite stage produces polynomials that may depend on the scale. The principle in Section~\ref{sec:entire} then constructs an entire function through summable modifications on expanding disks. Double factors preserve values and first derivatives at all nodes of completed stages. Outputs remain exactly equal by indistinguishability; we take no limit inside the rules. Choosing exponents $B+1/m$ yields simultaneous divergence for every $p>B$. This preservation-and-limit architecture is common to the general and paired barriers in the present work.

Thus independence of the proof means that all lemmas required for the case studied are proved here. It does not mean attributing a different methodological origin to this specialization than to the earlier work. Conversely, citing that work introduces no logical dependence when its theorem is not among the premises used. Presenting both relationships allows us to see how one architecture adapts to different query types while maintaining the assumptions and proof of each result.
\end{remark}

\subsection{What follows and what requires another proof}
Specialization to queries of orders zero and one proves the barrier needed for vector transfer. It does not by itself prove a barrier when higher derivatives are allowed: admitting more information enlarges the class of procedures and requires verifying the constraints, degrees, and displacements of the construction again. Nor does it allow us to deduce the paired bound simply by substituting a smaller number of queries into the general result. The paired improvement rests on different accounting, which we develop in the next section.

\section{The barrier for paired queries}\label{sec:paired-model}
We now study the information obtained when every first-derivative query is preceded by a function-value query at the same node. This property allows us to preserve responses using simple and double factors and improve the exponent of the finite stage. We develop this construction and verify its hypotheses before applying the fixed entire witness principle already proved.
\subsection{Pairing and repetition}
We use causal pairing as defined in Definition~\ref{def:pairing}. The function value must precede the derivative at the same node, but the two queries need not be consecutive.

Repeated queries of $F(y)$ or $DF(y)$ increase the nominal counter but not the distinct constraints in the record. This distinction will matter when constructing perturbations.

\begin{lemma}[Extension preserving pairing]\label{lem:paired-totalization}
If $\M$ has local budget $(N_F,N_J)$ and is locally paired, a truncated total extension can be chosen that respects the counters and whose records are paired in every global execution. It preserves the original valid executions respecting both the budget and pairing.
\end{lemma}
\begin{proof}
Add to Lemma~\ref{lem:totalization} one condition before preserving a query of type $1$: the history must contain an earlier query of type $0$ at the same node. Otherwise, return zero rather than making the query. This check uses only the record and equality of two known points. Every preserved execution satisfies the condition. By induction, the extension never makes an unpaired derivative query. On an original execution with the required property, no substitution is activated and all its data are preserved.
\end{proof}

The simulation in Theorem~\ref{thm:simulation} preserves pairing. Indeed, if $DF(Y)$ is queried after $F(Y)$, then $f'(Y_1)$ is requested after $f(Y_1)$. Projection can cause distinct vector nodes to coincide; this may increase scalar repetitions but cannot destroy pairing.

\subsection{Integers associated with a record}
For each distinct function node $u_i$ of a paired scalar record, assign $e_i\in\N$ and $\mu_i\in\{1,2\}$, with $\mu_i=1$ if only $f(u_i)$ is known and $\mu_i=2$ after the first query of $f'(u_i)$.

Start with an empty list and $E=0$. On adding a new function node, define
\begin{equation}\label{eq:exponent-update}
 e_{\mathrm{nuevo}}=\max\{1,E\},\qquad
 \mu_{\mathrm{nuevo}}=1,
 \qquad E^+=E+e_{\mathrm{nuevo}}.
\end{equation}
On querying the derivative for the first time at a node of multiplicity $1$, put
\begin{equation}\label{eq:derivative-update}
 \mu_i^+=2,\qquad E^+=E+e_i.
\end{equation}
Repetitions leave $e_i,\mu_i,E$ unchanged. We maintain
\(E=\sum_i\mu_i e_i.\)
These integers are used in the proof. We do not assume that the nodes depend regularly on a parameter or have a prescribed asymptotic order.

We use the function $B(N_F,N_J)$ defined in \eqref{eq:paired-budget}. Its value bounds the integer process in the next assertion.

\begin{lemma}[Growth of the construction integers]\label{lem:exponent-growth}
For every paired record with at most $N_F$ function queries and $N_J$ Jacobian queries, the preceding integers satisfy
\begin{equation}\label{eq:exponent-bound}
 \max\{1,E\}\le B(N_F,N_J).
\end{equation}
\end{lemma}
\begin{proof}
If there are no nodes, $E=0$ and the inequality is immediate. Suppose there are $v\ge1$ distinct function nodes and $r$ nodes at which a derivative has been queried. By pairing,
\(v\le N_F,\qquad r\le\min\{v,N_J\}.\)
The first node gives $(e_1,E)=(1,1)$; every subsequent distinct function node satisfies $e_{\rm nuevo}=E$ and $E^+=2E$.

Consider a first derivative query at a node other than the first. Let $e_i$ be its exponent. Before adding that node, the integer was $e_i$; after its function query, it became $2e_i$. Since the integer does not decrease, before querying the derivative we have $E\ge2e_i$. Consequently,
\(E^+=E+e_i\le\frac32E.\)
For the first node we need only $e_1=1\le E$, hence $E^+\le2E$. This last contribution can occur at most once.

If $r=0$, we obtain $E\le2^{v-1}$. If $r\ge1$ and the derivative has been queried at the first node, then
\[
 E\le2^{v-1}\,2\left(\frac32\right)^{r-1}
   =2^{v-r+1}3^{r-1}.
\]
If that derivative has not been queried, all $r$ derivative contributions multiply by a factor no greater than $3/2$, and
\[
 E\le2^{v-1}\left(\frac32\right)^r
   =\frac34\,2^{v-r+1}3^{r-1}
   \le2^{v-r+1}3^{r-1}.
\]
Repeated queries leave $E$ unchanged.

The bound $2^{v-1}$ increases with $v$. Passing from zero derivatives to one derivative multiplies the corresponding bound by $2$, and every additional distinct derivative multiplies it by $3/2$. For a fixed number of derivatives, adding one function node multiplies the expression by $2$. Thus the maximum allowed by $v\le N_F$ and $r\le\min\{v,N_J\}$ is exactly \eqref{eq:paired-budget}. This proves \eqref{eq:exponent-bound}.
\end{proof}

\begin{remark}
The maximum of this integer process is attained by beginning with a function node, querying its derivative when $s\ge1$, then making $s-1$ new value-and-derivative pairs, and adding the remaining $N_F-s$ values. Each later pair multiplies by $3$ and each remaining value by $2$. This observation concerns only the integer process; it does not assert the existence of a vector method attaining the barrier for every budget.
\end{remark}

\begin{example}[Integers of a causal record and the effect of repetition]\label{ex:paired-integers}
Consider the requests $f(u_1),f'(u_1),f(u_2),f'(u_2),f(u_3)$, at distinct nodes. The successive values of $E$ are
\(0\longmapsto1\longmapsto2\longmapsto4\longmapsto6\longmapsto12.\)
Thus $(e_1,e_2,e_3)=(1,2,6)$ and $(\mu_1,\mu_2,\mu_3)=(2,2,1)$. Hence $E=2\cdot1+2\cdot2+1\cdot6=12=B(3,2)$. If we repeat $f'(u_1)$ before the last request, the derivative counter increases, but multiplicities and $E$ remain unchanged. $E$ quantifies the constraints preserved by the perturbation, not the cost of operations.
\end{example}
\subsection{Polynomial tools for preserving responses}\label{sec:polynomial-tools}
\subsubsection{Norms on disks and normalized factors}
We use the uniform norm $\norm h_R$ and the bound $|h'(u)|\le\norm h_R$ for $R\ge3$ and $u\in[-1,1]$, already proved in Lemma~\ref{au:lem:cauchy}.

For a real node $u$, define the factor
\begin{equation*}
 L_u(z)=\frac{z-u}{1+|u|}.
\end{equation*}
$L_u\in\R[z]$ and $L_u(u)=0$. If $R\ge1$, then
\begin{equation*}
 \norm{L_u}_R\le R+1,
 \qquad
 |L_u(b)|\le1\quad\text{if }|b|\le\tfrac12.
\end{equation*}
Indeed,
\[
 |L_u(z)|\le\frac{|z|+|u|}{1+|u|}\le R+1,
 \qquad
 |L_u(b)|\le\frac{1/2+|u|}{1+|u|}\le1.
\]
The bounds are uniform for $u\in\R$, including $|u|\to\infty$.

Let $P$ be a fixed polynomial and $I\subset(-1/2,1/2)$ a compact interval on which $P$ does not vanish. Consider
\begin{equation}\label{eq:H-product}
 H(z)=P(z)\prod_{i=1}^v L_{u_i}(z)^{\mu_i},
 \qquad \mu_i\in\{1,2\},\quad \sum_i\mu_i\le d.
\end{equation}
Define
\[
 p_-:=\min_I|P|>0,\qquad
 K_0:=\max\{1,\max_I|P|\},\qquad
 K_R:=\max\{1,\norm{P}_R(R+1)^d\}.
\]
Then
\begin{equation*}
 \norm{H}_R\le K_R,\qquad
 |H(b)|\le K_0\quad(b\in I).
\end{equation*}
Moreover, for every node present and every $b\in I$,
\begin{equation}\label{eq:H-distance}
 |H(b)|\le K_0|b-u_i|.
\end{equation}
To prove \eqref{eq:H-distance}, extract a factor $L_{u_i}$ from \eqref{eq:H-product}. The remaining factors have modulus at most one, and $1/(1+|u_i|)\le1$.

\begin{lemma}[Preservation by simple and double factors]\label{lem:product-preserves}
Suppose a scalar record is paired. If every node at which only the function value is known has multiplicity $1$, and every node where the derivative is also known has multiplicity $2$, the polynomial $H$ in \eqref{eq:H-product} satisfies
\[
 H(u)=0\quad\text{at every function query},
 \qquad
 H'(v)=0\quad\text{at every derivative query}.
\]
Therefore, adding a multiple of $H$ preserves the record exactly.
\end{lemma}
\begin{proof}
At a function node $u$, $L_u(u)=0$; if $f'(u)$ is also queried, the factor $L_u^2$ gives $H(u)=H'(u)=0$. An additional factor $P$ does not alter these properties. Apply Lemma~\ref{lem:indistinguishability} in dimension one.
\end{proof}

\subsubsection{Using the root estimates already proved}
We apply Lemmas~\ref{au:lem:root-neighborhood} and \ref{au:lem:two-roots}. To distinguish the derivative lower bound from multiplicities $\mu_i$, we denote it by $m:=\mu$ in \eqref{au:eq:root-derivative}. Thus, for the roots of $g$ and $g\pm\eta H$,
\begin{equation}\label{eq:two-roots-shift}
\frac{\eta|H(b)|}{M}\le|b_\pm-b|\le\frac{\eta|H(b)|}{m}.
\end{equation}
Since $(b_+-b)(b_--b)<0$, for every $y\in\R$ one sign satisfies
\begin{equation}\label{eq:avoid-point}
|b_\pm-y|\ge\frac{\eta|H(b)|}{M}.
\end{equation}
In each application, we verify the hypotheses through uniform control of the candidates' norm.

\subsubsection{An amplitude controlling the root without losing the exponent}
Given $\lambda>0$ and $0<t<1$, for $H(b)\ne0$ define
\begin{equation}\label{eq:capped-amplitude}
 \eta=\lambda\min\left\{1,\frac{t}{|H(b)|}\right\}.
\end{equation}
Then
\begin{equation}\label{eq:capped-value}
 \eta|H(b)|=\lambda\min\{|H(b)|,t\},
 \qquad 0<\eta\le\lambda.
\end{equation}
The constant $\lambda$ is chosen before $t$. The amplitude $\eta$ may depend on data already constructed, since it belongs to the perturbation used in the proof, not to the method's rules.

\begin{lemma}[Preservation of separations and a new node]\label{lem:separation-step}
Suppose $H$ has the form \eqref{eq:H-product}, $b\in I$ does not coincide with any of its nodes, and $\eta$ is chosen by \eqref{eq:capped-amplitude}. Suppose also that Lemma~\ref{au:lem:two-roots} is applied and $\lambda K_0/m\le1/2$.

For either root $b_\pm$ and every previous node,
\begin{equation}\label{eq:old-distance}
 \frac{|b_\pm-u_i|}{1+|u_i|}
 \ge\frac12\frac{|b-u_i|}{1+|u_i|}.
\end{equation}
For every new node $y$, one of the two roots satisfies
\begin{equation}\label{eq:new-distance}
 \frac{|b_\pm-y|}{1+|y|}
 \ge\kappa\min\{|H(b)|,t\},
 \qquad
 \kappa=\min\left\{\frac{\lambda}{3M},\frac12\right\}.
\end{equation}
Moreover, $|b_\pm-b|\le\lambda t/m$.
\end{lemma}
\begin{proof}
By \eqref{eq:two-roots-shift}, \eqref{eq:H-distance}, and $\eta\le\lambda$,
\(|b_\pm-b|\le\frac{\lambda K_0}{m}|b-u_i| \le\frac12|b-u_i|.\)
The reverse triangle inequality, followed by division by $1+|u_i|$, proves \eqref{eq:old-distance}.

Choose the sign in \eqref{eq:avoid-point} for $y$. If $|y|\le2$, then $1+|y|\le3$ and
\[
 \frac{|b_\pm-y|}{1+|y|}
 \ge\frac{\lambda}{3M}\min\{|H(b)|,t\}.
\]
If $|y|>2$, since $|b_\pm|<1/2$,
\(\frac{|b_\pm-y|}{1+|y|} \ge\frac{|y|-1/2}{1+|y|}>\frac12.\)
The right-hand side is at least $\tfrac12\min\{|H(b)|,t\}$ because $t<1$. This proves \eqref{eq:new-distance}. The last assertion follows from \eqref{eq:capped-value} and the upper bound in \eqref{eq:two-roots-shift}.
\end{proof}

\begin{lemma}[Lower control of the product]\label{lem:product-lower}
Let $0<c\le1$ and $0<t<1$. If
\[
 \frac{|b-u_i|}{1+|u_i|}\ge c\,t^{e_i}
 \quad\text{for all nodes present},
\]
and $E=\sum_i\mu_i e_i$, then
\begin{equation}\label{eq:lower-H}
 |H(b)|\ge p_-c^d t^E,
 \qquad
 \min\{|H(b)|,t\}
 \ge\min\{p_-c^d,1\}\,t^{\max\{1,E\}}.
\end{equation}
The assertion includes the empty list, for which $E=0$.
\end{lemma}
\begin{proof}
By the product and $\sum_i\mu_i\le d$,
\[
 |H(b)|=|P(b)|\prod_i
 \left(\frac{|b-u_i|}{1+|u_i|}\right)^{\mu_i}
 \ge p_-c^{\sum_i\mu_i}t^E
 \ge p_-c^dt^E.
\]
If $E\ge1$, both $p_-c^dt^E$ and $t$ are at least $\min\{p_-c^d,1\}t^E$. If $E=0$, both $p_-c^d$ and $t$ are at least $\min\{p_-c^d,1\}t$. This proves the second inequality. For an empty list, the product equals one and the first inequality still holds.
\end{proof}

\subsection{Finite stage for arbitrary paired queries}\label{sec:finite-stage}
\begin{proposition}[Constants fixed before the scale]\label{prop:finite-stage}
Let $T$ be a total paired scalar method querying only $f$ and $f'$, with budgets $N_F\ge1$ and $N_J\ge0$. Put $d=N_F+N_J$ and $B=B(N_F,N_J)$. Let $f,P\in\R[z]$ and $a\in(-1/2,1/2)$ be such that
\(f(a)=0,\qquad f'(a)>0,\qquad P(a)\ne0.\)
For every $R\ge3$ and $\varepsilon>0$, there exist $c,\tau>0$, fixed before $t$, such that for all $0<t<\tau$, with $x=a+t$, we can construct a real polynomial $g$ and a simple root $b$ satisfying
\begin{equation}\label{eq:stage-norm}
 \norm{g-f}_R<\varepsilon,
 \qquad P\mid(g-f),
 \qquad \deg(g-f)\le\deg P+d,
\end{equation}
\begin{equation}\label{eq:stage-root}
 |b-a|<t/8,
 \qquad |T_g(x)-b|\ge c\,t^B.
\end{equation}
Moreover, $P(b)\ne0$ and $b$ differs from all nodes queried in this execution. The constants can be chosen in terms of $f,P,a,R,\varepsilon,N_F,N_J$, without depending on the nodes the method will query.
\end{proposition}
\begin{proof}
Apply Lemma~\ref{au:lem:root-neighborhood} and fix $I,r,\delta,m,M$, writing $m=\mu$ as indicated above. Also fix $p_-,K_0,K_R$ as in Subsection~\ref{sec:polynomial-tools}, with integer $d$. Choose
\begin{equation*}
 0<\lambda\le
 \min\left\{
 1,\frac{m}{8K_0},\frac{m}{32(d+1)},
 \frac{\delta}{4(d+1)K_R}
 \right\},
 \qquad
 0<\tau<\min\{1/4,r/2\}.
\end{equation*}
These choices precede $t$ and every query.

Define $\kappa=\min\{\lambda/(3M),1/2\}$, $c_0=1$, and, for $0\le j\le d$,
\begin{equation}\label{eq:c-recursion}
 c_{j+1}=\min\left\{
 \frac{c_j}{2},\ 
 \kappa\min\{p_-c_j^d,1\}
 \right\}.
\end{equation}
All these numbers are positive, at most one, and nonincreasing. Let $c_*=c_{d+1}$ and put
\begin{equation}\label{eq:c-stage}
 c=\min\left\{\frac18,
 \frac{\lambda}{M}\min\{p_-c_*^d,1\}\right\}>0.
\end{equation}

Now fix $0<t<\tau$ and $x=a+t$. Initially use the polynomial $f$ and its root $a$, with no acquired responses. We construct an execution of $T$ and, simultaneously, a current polynomial $g_0$, a current root $b_0$, and the list of distinct function nodes. We use subscripts only when necessary; at a given step we simply write $g,b$.

During the construction, we maintain the following properties. The current polynomial preserves all previously provided responses. Its root belongs to $I$, coincides with no function node in the list, and does not annihilate $P$. After $j$ polynomial modifications, each listed node satisfies
\begin{equation}\label{eq:stage-invariant}
 \frac{|b-u_i|}{1+|u_i|}\ge c_j t^{e_i}.
\end{equation}
Multiplicities and integers are updated by \eqref{eq:exponent-update} and \eqref{eq:derivative-update}.

We first check that all perturbations to be used are admissible in the neighborhood of Lemma~\ref{au:lem:root-neighborhood}. Every modification has the form $g\mapsto g\pm\eta H$, where $H$ is given by \eqref{eq:H-product} for the current record and $\eta$ is \eqref{eq:capped-amplitude}. Then
\(\norm{\eta H}_R\le\lambda K_R \le\frac{\delta}{4(d+1)}.\)
We make at most one modification for each new function node and one final modification; hence no more than $N_F+1\le d+1$. All intermediate polynomials and both candidates at each choice stay within distance $\delta/4<\delta$ of $f$. Thus they have a unique simple root in $I$ and satisfy the same derivative bounds $m,M$.

Moreover, each root displacement is at most $\lambda t/m$, by Lemma~\ref{lem:separation-step}. The sum of all displacements is at most
\begin{equation}\label{eq:stage-total-shift}
 (d+1)\frac{\lambda t}{m}\le\frac{t}{32}<\frac t8.
\end{equation}
Existence of the candidates was established earlier from the norm bound, not from this estimate; there is no circular dependence.

Consider a decision of the method from the current history. If it requests a function query at a node already in the list, respond with the value of $g$ at that node. This agrees with the earlier response, and we do not change the polynomial or the integers.

Suppose it requests a function query at a new node $y$. The node is determined by the already fixed history. Form $H$ by \eqref{eq:H-product}. Since $P(b)\ne0$ and $b$ is not a previous node, we have $H(b)\ne0$. Choose $\eta$ by \eqref{eq:capped-amplitude} and consider the two candidates $g\pm\eta H$.

Both candidates preserve the history exactly by Lemma~\ref{lem:product-preserves}. Hence the method requests the same node $y$ on both. Choose the sign given by Lemma~\ref{lem:separation-step}, replace the current polynomial and root by that candidate and its root, and only then answer the new function query.

For preceding nodes, \eqref{eq:old-distance} preserves at least half the separation. For the new node, Lemmas~\ref{lem:separation-step} and \ref{lem:product-lower} give
\[
 \frac{|b_{\mathrm{nuevo}}-y|}{1+|y|}
 \ge \kappa\min\{p_-c_j^d,1\}
          t^{\max\{1,E\}}.
\]
Assign $e_{\mathrm{nuevo}}=\max\{1,E\}$ and multiplicity $1$. By \eqref{eq:c-recursion}, \eqref{eq:stage-invariant} holds with $j+1$ for every node. The new root differs from $y$ and remains distinct from earlier nodes. Also, $P$ does not vanish at $I$.

Now suppose the next query is a derivative query. By pairing, its node $u_i$ belongs to the function-node list. Answer with $g'(u_i)$, leaving the polynomial and root unchanged. If the derivative has not been queried before, change the multiplicity from $1$ to $2$ and add $e_i$ to $E$. If it is repeated, leave the list unchanged. From then on, the double factor will also preserve the derivative response in all later modifications. The separation estimate is unchanged.

The budget and totality guarantee that this procedure reaches an output $q\in\R$ after at most $d$ queries. Each history is a true history for the current polynomial because all earlier responses were preserved exactly. Thus arbitrary decisions of the method introduce no additional difficulty.

Now form $H$ for the complete record and make one final modification. The two polynomials $g\pm\eta H$ produce the same record and output $q$. By Lemma~\ref{au:lem:two-roots}, choose a sign for which the new root, denoted by $b$, satisfies
\[
 |q-b|\ge\frac{\lambda}{M}\min\{|H(b_{\mathrm{anterior}})|,t\}.
\]
Before this modification, all nodes have separation at least $c_*t^{e_i}$, because the number of modifications does not exceed $d+1$ and the $c_j$ are nonincreasing. Lemma~\ref{lem:product-lower} implies
\[
 |q-b|\ge
 \frac{\lambda}{M}\min\{p_-c_*^d,1\}
       t^{\max\{1,E\}}.
\]
By Lemma~\ref{lem:exponent-growth}, $\max\{1,E\}\le B$. Since $0<t<1$, we conclude that $|q-b|\ge ct^B$.

The final polynomial preserves the output $q=T_g(x)$. All its corrections are multiples of $P$ and have degree at most $\deg P+d$, because the sum of the current multiplicities does not exceed the number of queries made. Its total norm is less than $\delta\le\varepsilon$. Estimate \eqref{eq:stage-total-shift} proves control of the root. Finally, \eqref{eq:old-distance} preserves separation from all nodes through the last modification, and $P(b)\ne0$ because $b\in I$. This proves every assertion.
\end{proof}

\begin{remark}[Uniformity and arbitrary rules]
The record may vary discontinuously with $t$; the integers $e_i$ encode only lower bounds and are not valuations of the nodes. The constants $\lambda,c_j,c,\tau$ are fixed before the scale, whereas the effective amplitude may depend on the constructed product and satisfies $0<\eta\le\lambda$. Identities \eqref{eq:capped-value} and \eqref{eq:lower-H} show that this choice preserves the required exponent. We do not need a stable rational branch.
\end{remark}

\subsection{A fixed entire witness for the paired barrier}
The preceding stage has the same form as the hypothesis of the preservation principle in Section~\ref{sec:entire}. We verify all its conditions to apply it with exponent $B(N_F,N_J)$, without repeating the limit construction.
\begin{theorem}[Paired scalar barrier with a fixed witness]\label{thm:scalar-paired}
Let $T$ be a total, deterministic, causal, memoryless scalar method that queries only values and first derivatives, always paired. Suppose its counters do not exceed $N_F\ge1$ and $N_J\ge0$. Let $B=B(N_F,N_J)$.

There exist $f\in\mathcal E_*$, its root $a$, and a sequence $x_m>a$, $x_m\to a$, such that
\begin{equation}\label{eq:paired-fixed-witness}
 \frac{|T_f(x_m)-a|}{(x_m-a)^p}\longrightarrow+\infty
 \qquad\text{for every }p>B.
\end{equation}
The function, root, and sequence do not depend on the final exponent.
\end{theorem}
\begin{proof}
The total budget $d=N_F+N_J$ is finite and positive, since $N_F\ge1$. The procedure $T$ is total, and its only queries are values and first derivatives. Expression \eqref{eq:paired-budget} satisfies $B(N_F,N_J)\ge1$.

Fix any real polynomial $f$, any polynomial $P$, a root $a\in(-1/2,1/2)$ with $f'(a)>0$ and $P(a)\ne0$, and numbers $R\ge3$ and $\varepsilon>0$. Proposition~\ref{prop:finite-stage} supplies $c,\tau>0$ before choosing $t$. Its definition \eqref{eq:c-stage} ensures $c\le1/8$. For every $0<t<\tau$ it produces a real polynomial $g$ and a simple root $b\in I\subset(-1/2,1/2)$ satisfying \eqref{eq:stage-norm} and \eqref{eq:stage-root}. Moreover, $P(b)\ne0$ and $b$ does not coincide with any node of that execution. These are precisely the conditions in \eqref{au:eq:abstract-stage} required by Proposition~\ref{au:prop:entire}, with $B=B(N_F,N_J)$.

Apply that principle. We obtain a function $f_\infty\in\mathcal E_*$, its root $a$, and a strictly decreasing sequence $x_m>a$, with $x_m\to a$, for which the quotients diverge simultaneously for all $p>B$. Renaming $f_\infty$ as $f$ gives \eqref{eq:paired-fixed-witness}.

In particular, the double factors used between stages preserve both values and first derivatives of every earlier record. Equality of outputs follows from Lemma~\ref{au:lem:indist}, not from continuity of the method. The same application retains $\norm{f-z}_3\le2^{-10}$ and the derivative bounds from the preservation principle. No result of \cite{CastroKT} is used as a premise.
\end{proof}
\subsection{Transfer of the paired barrier to systems}\label{sec:paired-vector}
The simulation preserves a value-and-derivative pair at each node. We therefore apply the paired scalar witness to the totalized procedure and then return to the original local executions.
\begin{theorem}[Paired barrier for arbitrary rules]\label{thm:paired-general}
Let $n\ge1$, $N_F\ge1$, and $N_J\ge0$. Let $\M$ be a deterministic, causal, memoryless method with local budget $(N_F,N_J)$ on $\Hclass(I_0)$. Suppose its records are paired in these neighborhoods.

Then there exist a fixed field $G_f$ with $f\in\mathcal E_*$, a simple root $\alpha$, and a sequence $X_m\to\alpha$ such that the original executions are valid for sufficiently large $m$ and divergence in \eqref{term:eq:extended} holds for every $p>B(N_F,N_J)$.
In particular, if $\M$ has universal local order at least $p>1$ on a class containing $\Hclass(I_0)$, then \eqref{term:eq:final} holds.
We do not prescribe the first node, require initialization with $F(x),DF(x)$, or impose regularity on the rules.
\end{theorem}
\begin{proof}
By Lemma~\ref{lem:paired-totalization}, the original rules admit a truncated total extension with paired records. Theorem~\ref{thm:simulation} yields a total scalar method with the same counters. Every query of $f'$ has a preceding query of $f$ at its node, because the simulation preserves pairing.

Apply Theorem~\ref{thm:scalar-paired}. We obtain $f\in\mathcal E_*$, $a=a_f$, and $x_m>a$, $x_m\to a$, satisfying scalar divergence for all exponents greater than $B(N_F,N_J)$. Put
\(\alpha=\iota(a),\qquad X_m=\iota(x_m).\)
For sufficiently large $m$, the original execution on the diagonal problem is valid, respects the budgets, and is paired. By preservation of executions under the extension and the simulation,
\(\bigl(\M_{G_f}(X_m)\bigr)_1=T_f(x_m).\)
If $I_0=(-1,1)$, compare the global extension with the restricted field exactly as in the proof of Theorem~\ref{thm:general}: every node of a valid original execution belongs to the domain and the responses agree there.

Finally,
\[
 \norm{X_m-\alpha}_\infty=x_m-a,
 \qquad
 \norm{\M_{G_f}(X_m)-\alpha}_\infty\ge|T_f(x_m)-a|.
\]
Scalar divergence implies \eqref{term:eq:extended}. If universal order greater than $B$ existed, the local estimate for this single field would bound the quotient, which is impossible.
\end{proof}

\begin{corollary}[Refinement of the Cordero--Torregrosa bound]\label{cor:paired-CT}
Under the hypotheses of Theorem~\ref{thm:paired-general}, if $1\le N_J\le N_F$, then
\begin{equation*}
 p\le2^{N_F-N_J+1}3^{N_J-1}
   =2^{N_F+N_J-1}\left(\frac34\right)^{N_J-1}.
\end{equation*}
The improvement is strict when $N_J>1$.
\end{corollary}
\begin{proof}
In this range, $s=\min\{N_F,N_J\}=N_J$. Substitute in \eqref{eq:paired-budget} and factor. If $N_J>1$, then $(3/4)^{N_J-1}<1$.
\end{proof}

\begin{corollary}[Budgets with repeated derivatives]\label{cor:saturation}
If $N_J\ge N_F\ge1$ and all derivatives are paired, then
\(p\le2\,3^{N_F-1}.\)
\end{corollary}
\begin{proof}
There are at most $N_F$ distinct function nodes. At each, the first derivative can contribute new information only once. Therefore $s=N_F$ in \eqref{eq:paired-budget}, and the expression equals $2\,3^{N_F-1}$.
\end{proof}

\begin{remark}[Analytic classes and larger differentiability classes]
Theorems~\ref{thm:general} and \ref{thm:paired-general} produce fields extending to entire fields. Thus their obstructions hold for the class of all real-analytic fields with simple roots. If the same method claims a common order on a class $C^r$ or $C^\infty$ containing these problems, it must also satisfy the corresponding bounds. No nonanalytic counterexample needs to be constructed.
\end{remark}

\begin{remark}[An error exponent need not be an integer]
The definition of order allows any real $p>1$. The integers $e_i$ and $E$ occur only in constructing a lower bound. The final contradiction uses $t^{B-p}\to+\infty$ when $p>B$ and thus does not require $p$ to be an integer or a Taylor expansion of the output map to exist.
\end{remark}

\subsubsection{Independent consequences without Jacobian queries}
\begin{corollary}[Methods querying only function values]\label{cor:values-only}
Let $N_F\ge1$ and let $T$ be a deterministic, causal, memoryless scalar method querying only function values and having local budget $N_F$ on the problems of $\mathcal E_*$, globally defined or restricted to $I_0$. Then there exist a fixed test function $f\in\mathcal E_*$, with simple root $a$, and a sequence $x_m>a$, $x_m\to a$, for which the original executions are eventually valid and
\[
 \frac{|T_f(x_m)-a|}{(x_m-a)^p}\longrightarrow+\infty
 \qquad\text{for every }p>2^{N_F-1}.
\]
In particular, universal local order greater than $2^{N_F-1}$ cannot be guaranteed.
\end{corollary}
\begin{proof}
Every record consisting only of values trivially satisfies pairing. Take $n=1$ and $N_J=0$ in Theorem~\ref{thm:paired-general}. Since $s=\min\{N_F,0\}=0$, formula \eqref{eq:paired-budget} gives $B(N_F,0)=2^{N_F-1}$. The construction and local treatment of domains are those of that theorem. This deduction uses the paired finite stage of Section~\ref{sec:finite-stage} and the preservation principle of Section~\ref{sec:entire}. It also agrees with the general scalar barrier, which has the same exponent for this budget.
\end{proof}

\begin{corollary}[Vector barrier without Jacobian queries]\label{cor:values-only-vector}
Let $N_F\ge1$ and let $\M$ be a deterministic, causal, memoryless vector method with local budget $(N_F,0)$ on $\Hclass(I_0)$. Then a fixed analytic field exists that simultaneously excludes all local exponents greater than $2^{N_F-1}$. In particular, if the method guarantees universal local order $p>1$ on a class containing these problems, then
\(p\le2^{N_F-1}.\)
\end{corollary}
\begin{proof}
There are no Jacobian queries, so every record satisfies pairing. Apply Theorem~\ref{thm:paired-general} directly with $N_J=0$ and use $B(N_F,0)=2^{N_F-1}$. The entire function and sequence come from the paired stage and the principle in Section~\ref{sec:entire}.
\end{proof}

\subsubsection{Pairing on the test fields}
\begin{definition}[Pairing in the first coordinate]\label{def:projected-pairing}
We say that $\M$ is locally paired in the first coordinate on $\Hclass(I_0)$ if, in the neighborhoods of well-definedness and budget, every Jacobian query $(1,Y)$ is preceded by some function query $(0,Y')$ with $Y'_1=Y_1$.
\end{definition}

\begin{proposition}[Refinement under projected pairing]\label{prop:projected-pairing}
Let $N_F\ge1$, $N_J\ge0$, and $\M$ be a method with local budget $(N_F,N_J)$ on $\Hclass(I_0)$. If Definition~\ref{def:projected-pairing} holds, then obstruction \eqref{term:eq:extended} and bound \eqref{term:eq:final} of Theorem~\ref{thm:paired-general} remain valid.
\end{proposition}
\begin{proof}
Construct the truncated extension of Lemma~\ref{lem:totalization}. Before preserving a Jacobian query $(1,Y)$, check that the record contains an earlier function query $(0,Y')$ with $Y'_1=Y_1$. Otherwise, return zero. These decisions depend on available data and use no additional information about the field. The extension is total, respects both counters, and all its global executions are paired in the first coordinate.

We show that it preserves the original executions specified in the statement. At every prefix of one of these executions, the original rule is defined, the next query respects the budget, and, if it is a Jacobian query, an earlier function value with the same first coordinate exists. Thus no substitution by the zero output is activated. Induction preserves all queries, responses, and the output.

Now apply Theorem~\ref{thm:simulation} to the extension. Every query of $DG_f(Y)$ becomes a query of $f'(Y_1)$. The earlier query of $G_f(Y')$, with $Y'_1=Y_1$, became a query of $f(Y_1)$. Hence the scalar method is paired in the usual sense and has the same counters. Theorem~\ref{thm:scalar-paired} provides $f\in\mathcal E_*$, a root $a$, and a sequence $s_m>a$, $s_m\to a$, simultaneously excluding every exponent greater than $B(N_F,N_J)$.

If $I_0=(-1,1)$, compare the global extension with the restricted field: every node of a valid original execution belongs to the cylinder and the responses agree there.

For this fixed problem, points $\iota(s_m)$ eventually lie in the neighborhood in the hypothesis. At those points the original outputs coincide with the outputs of the extension, whose first component is $T_f(s_m)$. Since
\[
 \norm{\iota(s_m)-\iota(a)}_\infty=s_m-a,\qquad
 \norm{\M_{G_f}(\iota(s_m))-\iota(a)}_\infty\ge|T_f(s_m)-a|,
\]
we obtain \eqref{term:eq:extended}. A universal estimate with a higher exponent contradicts this divergence, proving \eqref{term:eq:final}.
\end{proof}

\begin{remark}[Difference between the two notions of pairing]\label{rem:projected-pairing}
Equality of full vector nodes implies equality of their first coordinates, but the converse is false. For example, for $n\ge2$, nodes $Y$ and $Y+e_2$ have the same first coordinate but are distinct. On $G_f$, a query of $G_f(Y)$ followed by $DG_f(Y+e_2)$ provides $f(Y_1)$ and $f'(Y_1)$, although the evaluations are not paired at the same vector point.

The proposition concerns the projection and diagonal family fixed for the proof. It does not assert that $F(Y)$ determines $F(Y')$ for an arbitrary field, or that any pair of different nodes gives pairing. We retain the full-node condition in Theorem~\ref{thm:paired-general}; projected pairing is a separate, explicit extension.
\end{remark}

\section{Pairing without temporal precedence}\label{sec:terminal-pairing}
In this section we prove causal normalization of queries.
The idea is to preserve the decisions of the original procedure
even when we query a function value in advance. To this end, we distinguish
the records of queries actually made from those received by the
original rules. We verify, without assuming rationality or continuity,
that both procedures produce the same output and preserve
the local budgets.

Recall the sets $\V(R)$ and $\J(R)$ in \eqref{term:eq:sets} and condition $\J(R)\subseteq\V(R)$ of Definition~\ref{term:def:terminal}. The counters $n_F,n_J$ include repetitions; $\V(R),\J(R)$ are sets of distinct nodes. This distinction allows queries to be advanced without increasing the local budgets.

\subsection{Causal normalization by anticipation and reuse}\label{subsec:causal-normalization}

The proof represents the same procedure with a different order of response acquisition. We keep separate the record seen by the original procedure and the record of queries actually made by the representation.

\begin{lemma}[Truncation without imposing pairing]\label{term:lem:truncate}
Let $N_F,N_J\in\N_0$ and $d=N_F+N_J\ge1$. From a method $\M$, we construct a method $\A$ whose rules are total through level $d$, whose rule at level $d$ produces only outputs, and which respects both counters on global fields. The method $\A$ preserves all valid executions of $\M$ that respect these counters, without requiring pairing in their prefixes.
\end{lemma}
\begin{proof}
Use the truncation of Lemma~\ref{lem:totalization}. At a level $k<d$, preserve a defined original output. If the original rule requests a query, preserve it only when it does not exceed the corresponding counter. If the rule is undefined or the counter would be exceeded, return zero. At level $d$, preserve a defined original output and replace every other decision by zero.

Complete the rules with zero output on formally incompatible histories. They are total mappings fixed without knowing the field. In a valid original execution within budget, no replacement occurs: induction over prefixes preserves every query, its response, and the output. The construction also preserves valid executions on an open set, since their nodes remain in that set. At no point do we check whether a derivative has an earlier function value.
\end{proof}

\begin{proposition}[Exact causal normalization]\label{term:prop:normalization}
Let $\A$ be a truncated total method as in Lemma~\ref{term:lem:truncate}. There exists a method $\mathcal C$, determined solely by the rules of $\A$, with the following properties.
\begin{enumerate}[label=\textup{(\roman*)},leftmargin=2.2em]
\item Sus registros son causalmente emparejados.
\item For every valid execution of $\A$ from $x$ on $F$, the execution of $\mathcal C$ is also valid and
\begin{equation}\label{term:eq:same-output}
 \mathcal C_F(x)=\A_F(x).
\end{equation}
\item If $R$ is the original record and $\widetilde R$ the record of the representation, then
\begin{equation}\label{term:eq:counts}
 n_F(\widetilde R)=|\V(R)\cup\J(R)|,\qquad
 n_J(\widetilde R)=n_J(R).
\end{equation}
\item All effective query nodes belong to the set of nodes of $R$. In particular, no additional domain-availability condition is imposed.
\item On global fields, $\mathcal C$ is total with total budget at most $d+N_J$. If $R$ is paired at termination, then
\begin{equation}\label{term:eq:no-increase}
 n_F(\widetilde R)\le n_F(R),\qquad
 n_J(\widetilde R)=n_J(R).
\end{equation}
\end{enumerate}
\end{proposition}
\begin{proof}
\textit{Construction of the representation.}
At the beginning of the cycle, fix the same initial point $x$, an empty simulated record $R^{\rm sim}=()$, and an empty table $\mathcal C_0$. This table stores pairs $(y,F(y))$. Its domain is always a finite set of points already queried effectively. $R^{\rm sim}$ contains only the responses delivered to $\A$, not those obtained in advance.

Suppose $\A$, applied to the simulated history, produces an output $v$. The representation returns this same vector. If it requests a function query at $y$, distinguish two cases. If the value at $y$ is already stored, append the triple $(0,y,\mathcal C_0(y))$ to the simulated record without another field query. If it is not stored, make the effective query $F(y)$, store its response, and append the same triple to the simulated record.

Now suppose $\A$ requests $DF(y)$. If the function value at $y$ has not yet been stored, first make the effective query $F(y)$ and store its response. During this operation the simulated record remains unchanged. Then make the effective query $DF(y)$ and append only its response to the simulated record. If the function value was already stored, make the Jacobian query directly and proceed in the same way.

In particular, when we query $F(y)$ in advance, we do not prematurely deliver that datum to the original rules. The decision requesting $DF(y)$ remains pending until the Jacobian response arrives. This distinction prevents the advanced value from affecting the decisions of the procedure being reproduced.

\textit{Preservation of decisions and output.}
Fix a valid execution of $\A$ on $F$. We prove by induction on the number $j$ of answered original requests that the simulated record coincides with its prefix of length $j$. For $j=0$, the initial points coincide and both records are empty.

Suppose the claim holds for $j$. The original rule receives the same history, so it makes the same query or returns the same output. If it requests an unstored function value, the effective response is exactly $F(y)$. If the value is stored, it is also $F(y)$, since it was previously obtained from the same fixed field at the same point. If it requests the Jacobian, the response delivered to the simulated record is exactly $DF(y)$, regardless of whether the representation queried the function value first. In all cases the next original prefix is preserved. Induction reaches the output and proves \eqref{term:eq:same-output}.

The induction holds for adaptive nodes: equality of $R_j^{\rm sim}$ preserves the next decision of $\A$.

\textit{Availability of queries.}
An advanced function evaluation is made only at a point $y$ where the original procedure is requesting $DF(y)$. A valid execution satisfies $y\in\Omega$, hence $F(y)$ is defined. All other queries are made at their original nodes. Thus the representation is valid whenever the execution it reproduces is valid; no extension of the field outside its domain is required.

\textit{Accounting and pairing.}
Every point first appearing as a function or Jacobian node gives rise to exactly one effective query of $F$ and is stored. Its effective function query is never repeated. At termination, the table's domain is exactly $\V(R)\cup\J(R)$. This proves the first identity in \eqref{term:eq:counts}.

Every original request of $DF(y)$ corresponds to exactly one effective request of $DF(y)$. We do not reuse matrices here: repeated Jacobian queries are preserved and counted. This proves the second identity. Moreover, before each effective Jacobian query, $F$ has already been queried at that point, so the effective record is causally paired.

Since $|\V(R)\cup\J(R)|\le n_F(R)+n_J(R)\le d$, the total number of effective queries is at most $d+N_J$. When $\J(R)\subseteq\V(R)$, the union equals $\V(R)$, giving \eqref{term:eq:no-increase}.

\textit{Rules on effective histories.}
It remains to show that the preceding description defines a method of the same model, rather than a procedure receiving an external state. Given $x$ and an effective record of length $k$, reconstruct from the beginning the simulated record, the table, and, when needed, the pending Jacobian query. Use the effective-record responses in the order specified by the construction. If $F(y)$ belongs to the table, answer the simulated query without querying. Other operations consume the corresponding effective response or, if all $k$ responses have already been consumed, determine the next effective query.

If the original rule returns an output, that determines the next decision. If a supplied history is incompatible with the types or nodes required by the reconstruction, define the output there to be zero. The same can be done when it contains additional responses following termination. These histories are unreachable in valid executions.

The reconstruction contains no infinite internal loop. Every answered original query increases its simulated level, and there are at most $d$ of them. Advancing a function value introduces at most one additional phase before each Jacobian query. Thus a finite number of internal transitions leads to the next effective query or output. At level $d+N_J$ no further query can be needed in a compatible reconstruction; return zero on the other histories.

Thus the effective decision is a mapping of the history, independent of $F$. The table is reconstructed from data in the same cycle and starts empty in the next. The new method is deterministic, causal, and without external memory. We impose neither continuity nor computability on the rules.
\end{proof}

\begin{remark}[The future is not queried]\label{term:rem:future}
The representation does not check whether the original procedure will later request $F(y)$. It queries that value in advance whenever it must answer a first Jacobian query at $y$. Final inclusion \eqref{term:eq:terminal} is used only to prove that the function budget did not increase. Outside executions satisfying that inclusion, the function cost can increase up to the size of the union of node sets. The representation remains valid; \eqref{term:eq:no-increase} is restricted to $\J(R)\subseteq\V(R)$.
\end{remark}

\begin{corollary}[Normalization with local budgets]\label{term:cor:local-normal}
If $\M$ satisfies Definition~\ref{term:def:local} on $\mathscr F$, there exists a causally paired method $\mathcal C$ with the same local budget $(N_F,N_J)$ and the same output map in those neighborhoods. If a global extension with these counters is needed, $\mathcal C$ can be truncated without altering the specified local executions.
\end{corollary}
\begin{proof}
First apply Lemma~\ref{term:lem:truncate} without requiring causal pairing, obtaining $\A$. In the neighborhoods of Definition~\ref{term:def:local}, its execution coincides with that of $\M$, including the complete record and final inclusion of nodes.

Normalize $\A$ using Proposition~\ref{term:prop:normalization}. In those same neighborhoods,
\[
 \mathcal C_F(x)=\A_F(x)=\M_F(x),\qquad
 n_F(\mathcal C,F,x)\le N_F,\quad n_J(\mathcal C,F,x)\le N_J.
\]
Every Jacobian query of $\mathcal C$ has an earlier function value. It therefore satisfies the local-budget and causal-pairing assumptions required by the preceding theorem of this manuscript.

If a global extension with the original counters is desired, apply truncation again to the already normalized method. The preceding local executions remain unchanged, since their counters do not exceed the limits. Causal pairing is preserved under early termination: each Jacobian query already made had its earlier function value. There is no need to complete any pending query.
\end{proof}

\begin{remark}[Why truncation with the old condition is insufficient]
We cannot apply directly to $\M$ the rule that returns zero at the first Jacobian lacking an earlier function value. A Newton step querying $DF(x)$ before $F(x)$ satisfies the final condition and would be removed by that rule. The correct order is: truncation without a temporal condition, causal normalization, and only then truncation with the desired budgets. This construction avoids deciding, from a prefix, whether its Jacobian nodes will receive function values in the future.
\end{remark}

\subsection{Extension of the paired barrier}\label{subsec:terminal-consequences}

We apply the causal barrier of Theorem~\ref{thm:paired-general} to the normalized method, with the classes and budget already defined.

\begin{proof}[Proof of Theorem~\ref{term:thm:extension}]
By Corollary~\ref{term:cor:local-normal}, there exists a causally paired method $\mathcal C$ with local budget $(N_F,N_J)$ on the same test class. For every $f\in\mathcal E_*$, there is a punctured neighborhood of $\alpha_f$ in which
\begin{equation}\label{term:eq:local-equality}
 \mathcal C_{G_f}(x)=\M_{G_f}(x).
\end{equation}
The neighborhood may depend on $f$, as with the original budget. The rules of $\mathcal C$ were fixed before choosing $f$.

Apply the causal barrier \eqref{term:eq:extended} to $\mathcal C$. We obtain a fixed field $G_f$, its root $\alpha_f$, and a sequence $X_m\to\alpha_f$ with the stated simultaneous divergence. Once these objects are fixed, points $X_m$ belong to the neighborhood of \eqref{term:eq:local-equality} for all sufficiently large $m$. At these indices outputs coincide and both executions are valid. Consequently, equal outputs have identical quotients and satisfy \eqref{term:eq:extended} for every exponent greater than $B$. The index of coincidence depends on $f$ and its local neighborhood, not on $p$.

Finally, if universal local order $p>B$ held on a class containing the constructed problem, there would exist $C,\rho>0$, depending on this fixed problem, such that
\[
 \norm{\M_{G_f}(x)-\alpha_f}_\infty
 \le C\norm{x-\alpha_f}_\infty^p
 \quad\text{for }0<\norm{x-\alpha_f}_\infty<\rho.
\]
Evaluating at $X_m$ for large $m$ would bound the quotient in \eqref{term:eq:extended} by $C$, contradicting its divergence. Therefore \eqref{term:eq:final} holds.
\end{proof}

\begin{corollary}[The three budget regimes]
Under the hypotheses of Theorem~\ref{term:thm:extension}, the bounds are
\[
 \begin{array}{ll}
 p\le2^{N_F-1},&N_J=0,\\[0.3em]
 p\le2^{N_F-N_J+1}3^{N_J-1},&1\le N_J\le N_F,\\[0.3em]
 p\le2\,3^{N_F-1},&N_J\ge N_F\ge1.
 \end{array}
\]
For $N_J>1$, $B(N_F,N_J)<2^{N_F+N_J-1}$.
\end{corollary}
\begin{proof}
All three expressions follow by substituting $s=\min\{N_F,N_J\}$ into \eqref{eq:paired-budget}. If $1\le N_J\le N_F$, the ratio of the paired bound to the general one is $(3/4)^{N_J-1}$. If $N_J>N_F$, this ratio is
\[
 \frac{2\,3^{N_F-1}}{2^{N_F+N_J-1}}
 =\left(\frac34\right)^{N_F-1}2^{N_F-N_J}<1.
\]
This includes $N_F=1$, because the first factor then equals one and the second is less than one.
\end{proof}

\subsection{Relationship with the multiplicity construction}

The preceding normalization preserves the proof of the causal barrier without modification. We can also verify directly that the proposed exponent accounting for advanced derivatives is compatible with the same expression $B$.

Now consider all distinct nodes occurring in a prefix, whatever the type of their first query. Assign to each node $u_i$ an integer $e_i\ge1$ and a multiplicity $\mu_i$. The multiplicity is $1$ if only the value has been queried and $2$ from the moment the derivative is queried, even if the value has not yet been requested. Put $E=\sum_i\mu_i e_i$ and start with $E=0$.

If the first query at a new node is functional, use update \eqref{eq:exponent-update}. If it is Jacobian, define
\begin{equation}\label{term:eq:reverse-update}
 e_{\rm nuevo}=\max\{1,E\},\qquad \mu_{\rm nuevo}=2,
 \qquad E^+=E+2e_{\rm nuevo}.
\end{equation}
A first derivative at a node of multiplicity $1$ adds $e_i$ to the total integer and changes the multiplicity to $2$. A later function query at a node that began with a derivative does not change $E$. Neither do repeated queries.

\begin{proposition}[The combinatorial bound is preserved]\label{term:prop:integers}
In every completed record satisfying \eqref{term:eq:terminal} and budgets $(N_F,N_J)$, the preceding process satisfies \eqref{eq:exponent-bound}. The same bound holds for each prefix.
\end{proposition}
\begin{proof}
Replace, only in the abstract accounting, every first Jacobian query at a new node by two consecutive operations: introduction of a node of multiplicity $1$ followed immediately by a derivative query at that node. Both use the same exponent $e=\max\{1,E\}$ computed before introducing the node. The total increment is $e+e=2e$, agreeing with \eqref{term:eq:reverse-update}. The later function query at that node is treated as a repetition and does not alter the process.

This produces an abstract causal process with $v=|\V(R)\cup\J(R)|=|\V(R)|\le N_F$ node introductions and $r=|\J(R)|\le\min\{v,N_J\}$ first derivatives. The first node contributes $1$. Each later new value doubles $E$. A first derivative at a node other than the first adds $e_i$ when we already have $E\ge2e_i$, and therefore multiplies $E$ by a factor no greater than $3/2$. The derivative at the first node multiplies by a factor no greater than $2$, and this contribution occurs at most once.

Consequently, if $r=0$, $E\le2^{v-1}$. If $r\ge1$, $E\le2^{v-r+1}3^{r-1}$; if the first node receives no derivative, an even smaller bound results. The expressions increase as the available resources increase, and their maximum is \eqref{eq:paired-budget}. The case with no queries gives $E=0$ and also satisfies the assertion. Since the original process does not decrease, every prefix has an integer no larger than the final one.
\end{proof}

\begin{remark}[Overpreservation and budget]
In a direct perturbation proof, assigning multiplicity $2$ to the first Jacobian at $u_i$ preserves both value and derivative. Although the method does not yet know the value, the perturbation may satisfy this additional condition. What cannot be done is to identify, at each prefix, the sum of multiplicities with the number of queries already made: a first derivative query contributes multiplicity $2$ but counts as only one query at that moment. The correct bound comes from the completed record, since
\(\sum_i\mu_i=v+r\le N_F+N_J.\)
The normalization in Subsection~\ref{subsec:causal-normalization} avoids transferring these changes to the entire analytic construction: it produces a causally paired method to which the previously proved theorem applies.
\end{remark}

\subsection{Example and scope of the extension}

\begin{example}[Two Newton steps in reversed order]
Consider the procedure making the following requests, in this order,
\[
 DF(x),\quad F(x),\quad DF(y),\quad F(y),
 \qquad y=x-DF(x)^{-1}F(x),
\]
and return
\[
 \M_F(x)=y-DF(y)^{-1}F(y).
\]
Near a simple root of a field $C^2$, both matrices are invertible and each Newton step has quadratic error. By the calculation in Example~\ref{ex:model-newton}, the first correction satisfies $y-\alpha=O(\norm{x-\alpha}^2)$; applying it at $y$, we obtain $\M_F(x)-\alpha=O(\norm{x-\alpha}^4)$.

The procedure has budget $(2,2)$ and is paired at termination, but not causally. Its representation first queries $F(x)$, then $DF(x)$, and uses the stored value when the original procedure requests it. It does the same at $y$. If the nodes are distinct, the effective queries are $F(x),DF(x),F(y),DF(y)$; if they coincide, the value is reused and the budget only decreases. The output does not change. Theorem~\ref{term:thm:extension} allows us to apply the bound $B(2,2)=6$ to the procedure originally written in reversed order.
\end{example}

The same argument allows $DF(x),DF(y),F(x),F(y)$ with $y=\Theta(x,DF(x))$: the representation queries $F(x)$ and $F(y)$ in advance but does not deliver them to rule $\Theta$ or other rules before the original requests, without requiring regularity of $\Theta$.

\begin{remark}[Scope of normalization]
Pairing at termination permits the same local output maps as causal pairing, with the same maximum budgets and bound $B(N_F,N_J)$. The simulation does not preserve the physical order of requests: it preserves the original decisions and output. We do not assert equality of arithmetic or storage cost, which are outside the model, and the value table belongs only to the current cycle.

If a Jacobian node does not receive a function query during the cycle, normalization with these budgets does not apply. Nor does it imply a statement about the order $Q$ of all trajectories.
\end{remark}

\section{Examples, scope, and limitations}\label{sec:scope}
We specify the scope of the barriers through reuse of responses, a counterexample on polynomial fields, and the distinction between the local order of one cycle and the order of trajectories. These complementary results are not premises of the main barriers.
\subsection{Effective budget and reuse of responses}\label{sec:effective}
The diagonal reduction distinguishes the number of requests from the number of new scalar responses. We establish a bound based on the latter number and then show why two distinct vector nodes may provide the same scalar response.

\begin{proposition}[Simulation with exact reuse]\label{prop:cache}
Consider a vector method already truncated to a finite number $d$ of queries. In its diagonal simulation, each scalar response can be stored under the key $(\tau,Y_1)$. If a key appears again, the vector response is reconstructed without a new scalar query. The output remains unchanged, and the number of scalar queries is exactly the number of distinct keys in the diagonal record.
\end{proposition}
\begin{proof}
When a key $(\tau,Y_1)$ appears for the first time, we query $f^{(\tau)}(Y_1)$ and store the response. If the key reappears, we use the stored value. In the function case, the remaining components of the vector response are reconstructed from the coordinates of the current node, not those of the node where the key first appeared. In the Jacobian case, the remaining entries are always zeros and ones.

Hence each reconstructed vector response coincides exactly with the original one. Induction over the requests shows that the decisions and the output are preserved. Between two new scalar queries, at most $d$ vector requests are reproduced, so the internal processing terminates. Each distinct key gives rise to exactly one scalar query, proving the asserted count.
\end{proof}

\begin{corollary}[Bound by the number of distinct keys]\label{cor:effective}
Let $\M$ be a method with local budget $(N_F,N_J)$ on $\Hclass(I_0)$, where $d=N_F+N_J\ge1$. Suppose there exists $d_{\mathrm{ef}}\ge1$ such that, for each problem in this class, the record contains at most $d_{\mathrm{ef}}$ distinct keys $(\tau,Y_1)$ in a punctured neighborhood of its root. Then every universal local order $p>1$ on a class containing $\Hclass(I_0)$ satisfies
\(p\le2^{d_{\mathrm{ef}}-1}.\)
\end{corollary}
\begin{proof}
Apply Lemma~\ref{lem:totalization} and obtain a procedure $\widehat\M$ with finite global budget $d$, whose executions agree with the original ones near each root whenever the local budgets are respected. In those same neighborhoods, the set of distinct keys is unchanged.

We construct the scalar procedure $T'$. Given an initial point $s$ and a finite list of scalar responses, we simulate $\widehat\M$ from $\iota(s)$, retain the keys already queried, and reconstruct every vector response through \eqref{eq:reconstruction}. When a known key appears, we use its stored response without consuming another scalar input. When a new key appears, we request $f^{(\tau)}(u)$ and append the response to the list. Decisions use only the data available. After consuming a given list of responses, the rules may continue processing known keys until reaching a new query or an output. This process makes at most $d$ vector requests, so it terminates. On lists that do not correspond to a valid execution, we complete the rules with conventional outputs. Thus we obtain total scalar rules that are causal and independent of $f$.

By Proposition~\ref{prop:cache}, on valid original executions the output of $T'_f(s)$ equals the first component of $\M_{G_f}(\iota(s))$, and the effective scalar budget does not exceed $d_{\mathrm{ef}}$. We also totalize $T'$ with this last budget: if another new query were attempted after reaching it, we return zero; an earlier output is preserved. The truncation does not intervene in the neighborhoods under consideration. Denote the resulting procedure by $\widehat T'$.

Theorem~\ref{thm:scalar-input}, applied to $\widehat T'$, provides a fixed function $f\in\mathcal E_*$, its root $a$, and a sequence $s_m>a$, with $s_m\to a$, for which
\[
 \frac{|\widehat T'_f(s_m)-a|}{(s_m-a)^p}\longrightarrow+\infty
 \qquad\text{for every }p>2^{d_{\mathrm{ef}}-1}.
\]
For sufficiently large $m$, the points $\iota(s_m)$ lie in the neighborhood where both totalizations preserve the executions. Consequently,
\(\widehat T'_f(s_m)=T'_f(s_m) =\bigl(\M_{G_f}(\iota(s_m))\bigr)_1.\)
If the vector method had universal local order $p>2^{d_{\mathrm{ef}}-1}$, we would obtain
\[
 |\widehat T'_f(s_m)-a|
 \le\norm{\M_{G_f}(\iota(s_m))-\iota(a)}_\infty
 \le C_{G_f,a}(s_m-a)^p,
\]
contradicting the divergence. This proves the corollary.
\end{proof}

\begin{remark}[A request that adds no information]
A cycle that queries $F(x)$ twice, queries $DF(x)$, and returns $x-DF(x)^{-1}F(x)$ has nominal budget $(2,1)$ but only two distinct scalar keys. Its effective-budget barrier is two, rather than the nominal barrier four. Repeating the request adds no information, even though it counts as a request.
\end{remark}

\begin{example}[Distinct nodes with the same projected information]\label{ex:distinct-nodes-cache}
Let $n\ge2$ and let $e_2$ be the second canonical basis vector. From $x$, we query $V=F(x)$ and define
\(h=\sum_{i=1}^n V_i^2,\qquad Y=x+he_2.\)
We query $W=F(Y)$ and $A=DF(x)$; if $A$ is invertible, we return
\begin{equation}\label{eq:distinct-nodes-method}
 \M_F(x)=x-A^{-1}V+W-V-hAe_2.
\end{equation}
The budget is $(2,1)$. Near a simple root $\alpha$, if $x\ne\alpha$, we have $V\ne0$, so $h>0$ and $Y\ne x$.

If $F\in C^2$ near $\alpha$, then $V=O(\norm{x-\alpha}_\infty)$ and $h=O(\norm{x-\alpha}_\infty^2)$. Shrinking the neighborhood, $DF(x)^{-1}$ is uniformly bounded and the segment joining $x$ and $Y$ stays inside the domain. By the integral formula,
\[
 W-V-hAe_2
 =\int_0^1\bigl(DF(x+she_2)-DF(x)\bigr)he_2\,\dd s
 =O(h^2)=O(\norm{x-\alpha}_\infty^4).
\]
The calculation in Example~\ref{ex:model-newton}, applied to the same field $F\in C^2$ with simple root $\alpha$, gives
$ x-DF(x)^{-1}F(x)-\alpha=O(\norm{x-\alpha}_\infty^2)$.
Thus \eqref{eq:distinct-nodes-method} has local order at least two.

On the diagonal family, $Y_1=x_1$ and
\(G_f(Y)=G_f(x)+he_2,\qquad DG_f(x)e_2=e_2.\)
Hence the additional correction in \eqref{eq:distinct-nodes-method} vanishes, and the output is a Newton step. The three requests project onto the keys
\((0,x_1),\qquad(0,x_1),\qquad(1,x_1).\)
The function-valued vector responses are different, since the second component increases by $h$, but the second is reconstructed from the same scalar value and the current node $Y$. Only two distinct keys are required. Corollary~\ref{cor:effective} yields the barrier two, compared with the nominal barrier four. Finally, with $f(s)=s+s^2$ and $x=\iota(t)$, the first component of the output is $t^2/(1+2t)$; its nonzero leading term shows that exponent two is attained in this example.
\end{example}

\subsection{Nonregular rules and comparison of budgets}\label{subsec:scope-comparison}
\subsubsection{A discontinuous rule admitted by the model}
Recall that a finite rational program uses a finite tree of rational operations, sign decisions, and queries. On a fixed field, its data depend rationally on the input along each branch. We do not impose this restriction on the methods of our theorems.

\begin{example}[A discontinuous rule of order at least two]\label{ex:nonrational}
Let $\chi_{\Q}:\R\to\{0,1\}$ be the indicator function of the rational numbers. We query $F(x)$ and $DF(x)$, define $d=DF(x)^{-1}F(x)$ when the matrix is invertible, and return
\begin{equation}\label{eq:discontinuous-method}
 \M_F(x)=x-d+\chi_{\Q}(x_1)\norm{d}_2^2d.
\end{equation}
This rule is deterministic, causal, and without memory, with budget $(1,1)$. We do not assume that the decision $x_1\in\Q$ is computable.

Let $F\in C^2$ near a simple root $\alpha$. If $e=x-\alpha$, then $d=O(\norm e)$. By the Newton estimate in Example~\ref{ex:model-newton}, $x-d-\alpha=O(\norm e^2)$. The last summand of \eqref{eq:discontinuous-method} is $O(\norm e^3)$, since $|\chi_{\Q}|\le1$. Therefore, the method has local order at least two.

We show that its output map cannot be represented on all fields by a finite rational program. Take $F(x)=x$ and $x=\iota(t)$ with $t>0$. The first component of \eqref{eq:discontinuous-method} is $\chi_{\Q}(t)t^3$, which is discontinuous at every $t>0$. By contrast, on this field the query responses are rational expressions in the nodes. Induction over the tree shows that, along each branch, all data are rational functions of $t$. Since the tree is finite, the zeros of the nonidentically zero numerators and denominators appearing in its instructions form a finite set $Z$. Identically zero comparisons retain the equality branch, while all others retain their sign on each component of $(0,\infty)\setminus Z$. Thus, on one of those intervals a fixed branch is executed and its output is rational and continuous. This contradicts the discontinuity of $\chi_{\Q}(t)t^3$ on every interval. Nonrepresentability is proved.

\end{example}

\subsubsection{Values of the general and paired barriers}
Table~\ref{tab:bounds} includes the regime $N_J>N_F$, which is absent from the historical range.
\begin{table}[htbp]
\centering
\caption{General barrier and refinement under pairing at termination.}\label{tab:bounds}
\begin{tabular}{@{}cccc@{}}
\toprule
$N_F$ & $N_J$ & $2^{N_F+N_J-1}$ & $B(N_F,N_J)$\\
\midrule
1&1&2&2\\
1&2&4&2\\
1&3&8&2\\
2&1&4&4\\
2&2&8&6\\
2&3&16&6\\
3&2&16&12\\
3&3&32&18\\
4&3&64&36\\
4&4&128&54\\
\bottomrule
\end{tabular}
\end{table}

\subsection{An obstruction to extending the result to polynomial fields for arbitrary rules}\label{sec:polynomial-counterexample}
We now consider the class $\mathscr P_n$ of polynomial fields with simple roots, retaining arbitrary rules and problem-dependent validity radii. We construct a locally exact method that recovers the field algebraically from real responses. The construction uses the axiom of choice; it is not a computable numerical algorithm.

\begin{definition}[Local exactness for each problem]\label{neg:def:exact}
A method is locally exact on $\mathscr P_n$ if, for every $(\R^n,F,\alpha)\in\mathscr P_n$, there exists $\delta_{F,\alpha}>0$ such that
\(\M_F(x)=\alpha\qquad\text{when }\norm{x-\alpha}<\delta_{F,\alpha}.\)
The method must be the same for all problems. Only the radius may depend on them.
\end{definition}

Local exactness implies an estimate of universal local order at least $p$ for each $p>1$: take $C_{F,\alpha}=1$ in the same neighborhood. This does not assign an infinite exact $Q$ order to a terminating trajectory. These are estimates for the error of one cycle.

\subsubsection{Subfields and evaluations determining a polynomial}
\paragraph{An algebraic exhaustion of the real numbers}
A set $S\subset\R$ is algebraically independent over $L\subset\R$ if, for any distinct $s_1,\ldots,s_k\in S$ and any $P\in L[Z_1,\ldots,Z_k]\setminus\{0\}$, we have $P(s_1,\ldots,s_k)\ne0$. A transcendence basis $\mathcal T$ of $\R/\Q$ is algebraically independent and $\R/\Q(\mathcal T)$ is algebraic.

\begin{lemma}[A sequence of subfields with independent probes]\label{neg:lem:fields}
For fixed $n\ge1$, there exist subfields
\[
 K_1\subseteq K_2\subseteq\cdots\subset\R,
 \qquad \bigcup_{m\ge1}K_m=\R,
\]
and vectors $\theta_m=(\theta_{m,1},\ldots,\theta_{m,n})\in(0,1)^n$ such that the coordinates of $\theta_m$ are algebraically independent over $K_m$. The $K_m,\theta_m$ are independent of $F$.
\end{lemma}
\begin{proof}
First choose a countable family
\(\{\theta_{m,j}:m\in\N,\ 1\le j\le n\}\subset(0,1)\)
algebraically independent over $\Q$. To construct it, enumerate the pairs $(m,j)$. At each step, the field generated by the elements already chosen is countable, since only finitely many elements have been selected. The real numbers algebraic over a countable field form a countable set: there are countably many polynomials, and each has finitely many roots. We can therefore choose the next element of $(0,1)$ outside that set.

Extend this family to a maximal algebraically independent family $\mathcal T\subset\R$. Zorn's lemma gives such an extension: the union of a chain is still algebraically independent, since each polynomial relation involves only finitely many elements. Maximality implies that $\R$ is algebraic over $\Q(\mathcal T)$. Indeed, if some real number were transcendental over that field, we could adjoin it to $\mathcal T$.

For each $m\ge1$, define
\[
 S_m=\{\theta_{k,j}:k\ge m,\ 1\le j\le n\},
 \qquad L_m=\Q(\mathcal T\setminus S_m)
\]
and
\(K_m=\{a\in\R:a\text{ is algebraic over }L_m\}.\)
$K_m$ is a subfield of $\R$. If $a,b$ are algebraic over $L_m$, the field $L_m(a,b)$ is a finite extension of $L_m$ and contains their sums, products, and defined quotients, which are therefore algebraic. $S_{m+1}\subset S_m\Longrightarrow L_m\subseteq L_{m+1}\Longrightarrow K_m\subseteq K_{m+1}$.

We show that the union of the $K_m$ is $\R$. Let $a\in\R$. Since it is algebraic over $\Q(\mathcal T)$, it satisfies a nonzero polynomial equation whose coefficients are rational functions of finitely many elements of $\mathcal T$. Thus there is a finite subset $A\subset\mathcal T$ such that $a$ is algebraic over $\Q(A)$. If $m$ is larger than every index of the blocks occurring in $A$, then $A\cap S_m=\varnothing$. It follows that $\Q(A)\subset L_m$ and hence $a\in K_m$. The same inclusion holds for every later index.

Finally, the coordinates of $\theta_m$ are algebraically independent over $L_m$. A relation with coefficients in $L_m$, after clearing denominators, would be a nonzero polynomial relation between distinct elements of $\mathcal T$. We show that this independence is preserved when passing to $K_m$. For one coordinate, being algebraic over $K_m$ would imply, by transitivity of algebraic extensions, being algebraic over $L_m$. For several coordinates, argue by induction. If the first $j-1$ coordinates are independent over $K_m$ and a relation involves coordinate $j$, then that coordinate would be algebraic over
\(K_m(\theta_{m,1},\ldots,\theta_{m,j-1}).\)
This field is algebraic over $L_m(\theta_{m,1},\ldots,\theta_{m,j-1})$. By transitivity we would obtain a dependence over $L_m$, a contradiction. This proves the required independence.
\end{proof}

\begin{corollary}[Coefficients of a fixed field enter the same subfield]\label{neg:cor:coefficients}
For every real polynomial field $F:\R^n\to\R^n$, there exists $m_F\in\N$ such that all its coefficients belong to $K_m$ for every $m\ge m_F$.
\end{corollary}
\begin{proof}
If $F$ has coefficients $a_1,\ldots,a_r$ and $a_i\in K_{m_i}$, take $m_F=\max_i m_i$: by monotonicity, $a_i\in K_m$ for every $m\ge m_F$.
\end{proof}

\paragraph{A fixed rule for recovering coefficients}
\begin{lemma}[Injective evaluation]\label{neg:lem:decoder}
For every $m$, the map
\[
 \operatorname{ev}_m:K_m[Z_1,\ldots,Z_n]\longrightarrow\R,
 \qquad P\longmapsto P(\theta_m)
\]
is injective. Consequently, there exists a total map
\[
 \mathcal D_m:\R\longrightarrow K_m[Z_1,\ldots,Z_n]\sqcup\{\bad\}
\]
that returns the unique polynomial $P$ satisfying $P(\theta_m)=a$, when one exists, and returns the marker $\bad$ otherwise.
\end{lemma}
\begin{proof}
If $P(\theta_m)=Q(\theta_m)$, then $(P-Q)(\theta_m)=0$. Algebraic independence of the coordinates implies $P-Q=0$. Thus the inverse on the image is well defined. Defining it to equal $\bad$ extends it to the rest of $\R$.
\end{proof}

The map $\mathcal D_m$ does not query the unknown field. It is a mapping of a real number already received. Its definition uses only the fixed objects in Lemma~\ref{neg:lem:fields}. We do not claim it is computable or that membership in the image of $\operatorname{ev}_m$ can be decided algorithmically.

\begin{example}[What it means to recover a polynomial from one response]\label{ex:algebraic-decoder}
In one dimension, fix an index $m$ and consider $P(Z)=Z^2-1/16$. $\Q\subseteq K_m$, so $P\in K_m[Z]$. The response $w=P(\theta_m)=\theta_m^2-1/16$ determines a unique polynomial in $K_m[Z]$: if $Q(\theta_m)=w$, then $(Q-P)(\theta_m)=0$, and the transcendence of $\theta_m$ over $K_m$ implies $Q=P$.

No prior upper bound on $\deg P$ is assumed. We apply a fixed inverse of an injective evaluation on all of $K_m[Z]$. This inverse is a mathematical rule allowed by the model, but we do not assert it can be computed. Membership of the coefficients in the field is essential. For an arbitrary real polynomial, a suitable index exists only after fixing that polynomial, as established by Corollary~\ref{neg:cor:coefficients}.
\end{example}

\subsubsection{Probes near the initial point and root selection}
\paragraph{Location of the second query}
For $m\ge1$ and $x\in\R^n$, define
\begin{equation}\label{neg:eq:probe}
 h_m=2^{-m},\qquad
 a_j(x,m)=h_m\left\lfloor\frac{x_j}{h_m}\right\rfloor,
 \qquad y(x,m)=a(x,m)+h_m\theta_m.
\end{equation}
Every coordinate of $a(x,m)$ is rational. By the definition of the floor function,
\(a_j(x,m)\le x_j<a_j(x,m)+h_m.\)
From $0<\theta_{m,j}<1$ it follows that $a_j(x,m)<y_j(x,m)<a_j(x,m)+h_m$. Consequently,
\begin{equation}\label{neg:eq:probe-distance}
 \norm{y(x,m)-x}<h_m.
\end{equation}

\begin{lemma}[Recovery from a translated probe]\label{neg:lem:shifted-decode}
Let $F$ be polynomial and suppose its coefficients belong to $K_m$. From $x$, $m$, and a single vector response $W=F(y(x,m))$, all coefficients of $F$ can be recovered exactly by fixed mappings independent of $F$.
\end{lemma}
\begin{proof}
Write $a=a(x,m)$ and $h=h_m$. For each component define, only to verify the construction,
\(Q_i(Z)=F_i(a+hZ).\)
Its coefficients belong to $K_m$, since those of $F_i$ belong to that field and $a,h$ are rational. Moreover,
\(Q_i(\theta_m)=F_i(y(x,m))=W_i.\)
By Lemma~\ref{neg:lem:decoder}, $\mathcal D_m(W_i)=Q_i$. The rules of the method obtain these polynomials by applying $\mathcal D_m$ to the responses; they do not receive $Q_i$ as additional data. Once obtained, they form
\(\widehat F_i(Z)=Q_i\!\left(\frac{Z-a}{h}\right).\)
Since $h>0$, we have $\widehat F_i=F_i$. This proves exact recovery.
\end{proof}

\paragraph{A fixed selection rule that does not receive the desired root}
For a polynomial field $P$, put $Z(P)=\{z\in\R^n:P(z)=0\}$. Define a map $\mathcal S(P,x)$ as follows. If $Z(P)=\varnothing$, put $\mathcal S(P,x)=x$. Otherwise, consider the zeros minimizing $\norm{z-x}$, and select the lexicographically smallest one.

\begin{lemma}[Well-definedness and local selection]\label{neg:lem:root-selector}
The map $\mathcal S$ is well defined. If $\alpha$ is a simple root of $P$, there exists $r>0$ such that
\(\mathcal S(P,x)=\alpha\qquad\text{for }\norm{x-\alpha}<r.\)
\end{lemma}
\begin{proof}
Since $P$ is continuous, $Z(P)=P^{-1}(\{0\})$ is closed. If it is nonempty, take a sequence of zeros whose distances from $x$ approach the infimum. The sequence can be chosen bounded, and a subsequence converges to an element of $Z(P)$ attaining that infimum. The set of minimizers is closed and bounded, hence compact. To select its lexicographically smallest element, minimize the first coordinate, then the second, and so on through the last, retaining at each step the compact subset where the minimum is attained. At the end we obtain a unique vector.

Since $DP(\alpha)$ is invertible, the inverse function theorem gives a neighborhood in which $\alpha$ is the only root. Choose $R>0$ such that there is no other zero in the ball of radius $R$ centered at $\alpha$. If $\norm{x-\alpha}<R/3$ and $\beta\ne\alpha$ is another zero, then
\(\norm{\beta-x}\ge\norm{\beta-\alpha}-\norm{x-\alpha} >2R/3>\norm{x-\alpha}.\)
Thus $\alpha$ is the unique zero nearest to $x$. It suffices to take $r=R/3$.
\end{proof}

$\mathcal S$ uses only $(\widehat F,x)$. The distinguished root of the problem is not supplied to the method. If the field has several simple roots, the same selector works near each of them.

\subsubsection{A locally exact method with budget \texorpdfstring{$(2,1)$}{(2,1)}}
\paragraph{Definition of the rules}
Fix once and for all the preceding subfields, probes, and mappings. From $x\in\R^n$, the method proceeds as follows.

First query $U=F(x)$. If $U=0$, return $x$. If $U\ne0$, query $A=DF(x)$. If $A$ is singular, return $x$. Otherwise define
\begin{equation}\label{neg:eq:scale}
 d=A^{-1}U,\qquad \rho=\norm d>0,\qquad
 m=\max\{1,\lceil\rho^{-1}\rceil\}.
\end{equation}
Construct $a,h_m,y=y(x,m)$ by \eqref{neg:eq:probe} and query $W=F(y)$.

Now apply $\mathcal D_m$ to each component of $W$. If any returns $\bad$, terminate with output $x$. Otherwise obtain polynomials $Q_1,\ldots,Q_n$, form
\begin{equation*}
 \widehat F_i(Z)=Q_i\!\left(\frac{Z-a}{h_m}\right),
 \qquad 1\le i\le n,
\end{equation*}
and return
\begin{equation}\label{neg:eq:method}
 \M_F(x)=\mathcal S(\widehat F,x).
\end{equation}
The return branches ensure that all decisions are defined for every field of class $C^1$ throughout $\R^n$. On numerical histories incompatible with earlier decisions, we may define the output to be $x$. No new query is made within $\mathcal D_m$ or $\mathcal S$.

\begin{theorem}[Local exactness with three queries]\label{neg:thm:exact}
For each $n\ge1$, the preceding rules define a deterministic, causal, memoryless method, total on global fields $C^1$, which uses at most two evaluations of $F$ and one of $DF$. Its Jacobian evaluation is causally paired. For every problem $(\R^n,F,\alpha)\in\mathscr P_n$, there exists $\delta_{F,\alpha}>0$ such that
\begin{equation}\label{neg:eq:local-exact}
 \M_F(x)=\alpha\qquad\text{if }\norm{x-\alpha}<\delta_{F,\alpha}.
\end{equation}
If there is a second function query, $\norm{y-x}\to0$ whenever $x\to\alpha$.
\end{theorem}
\begin{proof}
The only requests to the field are $F(x)$, $DF(x)$, and $F(y)$. The second is made after knowing $F(x)$, and the third is determined by already received data. All remaining operations are fixed mappings of those data. Thus the rules are deterministic and causal. Each cycle begins with $x$ and the same fixed algebraic objects, without using responses from previous cycles. The method is memoryless. Its sole Jacobian evaluation has an earlier function evaluation at the same point, and the budgets are respected globally.

Now let $F$ be polynomial and let $\alpha$ be a simple root. Fix the integer $m_F$ of Corollary~\ref{neg:cor:coefficients}. Near $\alpha$, the matrix $DF(x)$ is invertible and its inverse depends continuously on $x$. Moreover, $F(x)\to0$ as $x\to\alpha$, so
\(DF(x)^{-1}F(x)\longrightarrow0.\)
The root is isolated. Shrink the neighborhood so that $F(x)\ne0$ whenever $x\ne\alpha$. In that neighborhood, $d\ne0$ and $\rho>0$. It follows from \eqref{neg:eq:scale} that
\(m(x)\longrightarrow+\infty\qquad(x\to\alpha,\ x\ne\alpha).\)
Therefore there exists a punctured neighborhood in which $m(x)\ge m_F$. For each point of that neighborhood, all coefficients of $F$ belong to $K_{m(x)}$. Lemma~\ref{neg:lem:shifted-decode} shows that no decoding fails and the recovered field satisfies $\widehat F=F$.

By Lemma~\ref{neg:lem:root-selector}, we can shrink the neighborhood once more so that $\mathcal S(F,x)=\alpha$. Thus \eqref{neg:eq:method} gives exactly \eqref{neg:eq:local-exact}. At $x=\alpha$, the first response is zero and the method returns $x=\alpha$, as required.

Finally, \eqref{neg:eq:probe-distance} gives
\(\norm{y-x}<2^{-m(x)}\longrightarrow0.\)
In particular, $y\to\alpha$. We have not used a probe that remains away from the root or a global condition on the method's trajectories.
\end{proof}

\begin{corollary}[No finite barrier on the polynomial class]\label{neg:cor:refutation}
For each $n\ge1$, there is no finite upper bound, valid for all methods with arbitrary rules and budget $(N_F,N_J)=(2,1)$, on the universal local order required only on $\mathscr P_n$.
\end{corollary}
\begin{proof}
The method of Theorem~\ref{neg:thm:exact} uses budgets $(2,1)$ and is locally exact near every simple root of every polynomial field of dimension $n$. Therefore, given any $p>1$ and any fixed problem in this class, there is a neighborhood in which
\(\norm{\M_F(x)-\alpha}=0\le\norm{x-\alpha}^{p}.\)
The same method thus satisfies a universal estimate of every finite exponent $p$, with radii that may depend on the problem. There can be no finite bound for all methods in this class with that budget.
\end{proof}

\begin{remark}[The method uses no hidden information about the field]
The objects $K_m,\theta_m,\mathcal D_m,\mathcal S$ are fixed independently of $F$. The integer $m_F$ is used only to prove the existence of a neighborhood of exactness; the method neither knows nor receives it. The index it uses is computed from $F(x)$ and $DF(x)$ through \eqref{neg:eq:scale}. The coefficients are recovered by a fixed map of the second function response. They are not supplied as external data.

Allowing arbitrary rules means admitting that map just as any other function of the responses. Excluding it would require an additional restriction on the rules. The theorem does not assert that the operations can be carried out on a computer or that a real number can be stored with finite precision.
\end{remark}

\begin{remark}[Two function queries also suffice]
The Jacobian evaluation is not essential for this phenomenon. In \eqref{neg:eq:scale}, we may replace $\rho=\norm{DF(x)^{-1}F(x)}$ by $\rho=\norm{F(x)}$, omit the Jacobian query, and keep the rest of the construction. Since $F(x)\to0$ and the root is isolated, $m(x)\to\infty$. The same proof gives local exactness with two function values. In Theorem~\ref{neg:thm:exact}, we retained the Jacobian to obtain an example with initial queries $F(x)$ and $DF(x)$ and budget $(2,1)$.
\end{remark}

\subsubsection{Why this does not contradict indistinguishability}
Two distinct polynomials may produce the same record in an execution. The construction does not claim that both can be recovered correctly from that record. It asserts that, for each fixed polynomial, all executions sufficiently near a simple root use indices for which its coefficients already belong to the decoding field. The radius of that neighborhood may change with the polynomial.

The following example makes this difference explicit.

\begin{example}[Lack of uniformity in a family of degree three]\label{neg:ex:nonuniform}
Fix $m\ge2$, put $t_m=1/m$, $h_m=2^{-m}$, and $\eta_m=h_m\theta_{m,1}$. Consider the fields of any dimension $n\ge1$
\(L_m(z)=(z_1+t_m,z_2,\ldots,z_n)^T\)
and
\[
 G_m(z)=\bigl(z_1+t_m+z_1^2(z_1-\eta_m),z_2,\ldots,z_n\bigr)^T.
\]
From $x=0$, the two fields give
\(F(0)=(t_m,0,\ldots,0)^T,\qquad DF(0)=I.\)
Thus the method chooses exactly index $m$ and the same probe $y=h_m\theta_m$. Since $y_1=\eta_m$, their function responses at $y$ also coincide.

All coefficients of $L_m$ are rational. Decoding recovers $L_m$ and the output is its root $(-t_m,0,\ldots,0)^T$. By equality of the records, this is also the output on $G_m$.

Let
\(g_m(s)=s+t_m+s^2(s-\eta_m).\)
Its derivative satisfies $g_m'(s)=1+3s^2-2\eta_ms=1+3(s-\eta_m/3)^2-\eta_m^2/3>0$.
since $0<\eta_m<2^{-m}\le1/4$. The polynomial has a unique real root $b_m$. Since $g_m(-t_m)=-t_m^3-\eta_mt_m^2<0$ and $g_m(0)=t_m>0$, we have $-t_m<b_m<0$. The mean value theorem gives
\[
 b_m+t_m=\frac{t_m^3+\eta_mt_m^2}{g_m'(\xi_m)},
 \qquad -t_m<\xi_m<b_m.
\]
Now, $\eta_m/t_m\le m2^{-m}\to0$ and $g_m'(\xi_m)\to1$. Consequently,
\(b_m+t_m\sim t_m^3,\qquad |b_m|\sim t_m.\)
If $\beta_m=(b_m,0,\ldots,0)^T$, we obtain
\[
 \frac{\norm{\M_{G_m}(0)-\beta_m}}{\norm{0-\beta_m}^{p}}
 \sim t_m^{3-p}\longrightarrow+\infty\qquad(p>3).
\]
The coefficients of $G_m$ converge to those of $(z_1+z_1^3,z_2,\ldots,z_n)^T$, and $\beta_m\to0$. All fields have degree at most three and simple roots. Thus even on this compact family of coefficients and roots there is no uniform estimate of exponent $p>3$ on a common radius. Nevertheless, each fixed field $G_m$ has its own neighborhood of exactness, by Theorem~\ref{neg:thm:exact}.

The coefficient $-\eta_m$ uses precisely the block $\theta_m$ excluded from $K_m$. The selected execution occurs at a scale that is not yet suitable for this field. As we approach its root more closely while holding $G_m$ fixed, the index of the construction exceeds $m$ and that coefficient enters the corresponding subfield. The two statements do not contradict one another.
\end{example}

\subsubsection{Scope and reformulations that remain open}
The result depends on the polynomial class, arbitrary rules, and validity radii that may vary with the problem; it works in every finite dimension. Finiteness of the number of coefficients is essential: although the coefficients of an analytic series may individually belong to different subfields $K_m$, there need not be a common index for all of them. The map $\mathcal D_m$ decodes polynomials, and we have not proved an analogous recovery for general analytic functions.

Thus the counterexample does not affect the barriers established on classes containing our entire test fields. Among the open questions, we distinguish imposing regularity on rules in the polynomial class and studying uniform guarantees over families of coefficients. Example~\ref{neg:ex:nonuniform} shows that our procedure does not preserve a uniform estimate of order higher than three even on a compact family of bounded degree. These questions differ from the per-problem local guarantee used in this article.

\subsection{Order of one cycle and order along trajectories}\label{sec:trajectory-order}
We compare the universal local estimate \eqref{eq:order} with the order along iterative trajectories, without assuming that regularity of the field is inherited by the output map.

\subsubsection{Local estimates and error equations}
Fix a problem $(\Omega,F,\alpha)$ on which the method is defined near the root. Write
\(E(e)=\M_F(\alpha+e)-\alpha\qquad(0<\norm e_\infty<\rho)\)
and, when iteration from the root is needed, define $E(0)=0$. This extension fixes the root without making an additional query or altering any execution with input different from it.

\begin{proposition}[Consequences of a one-cycle estimate]\label{prop:cycle-consequences}
Let $p>1$. The following assertions hold.
\begin{enumerate}[label=\textup{(\roman*)},leftmargin=2.3em]
\item Si
\[
 \limsup_{e\to0,\ e\ne0}\frac{\norm{E(e)}_\infty}{\norm e_\infty^p}<+\infty,
\]
then there exist $C,\delta>0$ such that $\norm{E(e)}_\infty\le C\norm e_\infty^p$ for $0<\norm e_\infty<\delta$.
\item Let $m\ge2$ be an integer. If $P_m:\R^n\to\R^n$ is a homogeneous polynomial of degree $m$ and
\begin{equation}\label{eq:local-error-equation}
 E(e)=P_m(e)+o(\norm e_\infty^m)\qquad(e\to0),
\end{equation}
then a local estimate with exponent $m$ holds. If $P_m\not\equiv0$, no local estimate with an exponent greater than $m$ holds.
\item If $\norm{E(e)}_\infty\le C\norm e_\infty^p$ near zero, every trajectory starting sufficiently near zero remains in that neighborhood and converges to zero, or reaches zero after finitely many steps. For each nonterminating trajectory,
\[
 \limsup_{k\to\infty}
 \frac{\norm{e_{k+1}}_\infty}{\norm{e_k}_\infty^p}\le C,
 \qquad
 \frac{\norm{e_{k+1}}_\infty}{\norm{e_k}_\infty}\longrightarrow0.
\]
\end{enumerate}
\end{proposition}
\begin{proof}
For the first assertion, let $L$ be the limsup, which is finite by hypothesis. By the definition of the limit superior there exists $\delta>0$ such that the quotient is less than $L+1$ in the punctured neighborhood of radius $\delta$. Taking $C=L+1$ gives the estimate.

For the second, continuity of $P_m$ on the compact sphere $\norm v_\infty=1$ gives
\(A=\max_{\norm v_\infty=1}\norm{P_m(v)}_\infty<+\infty.\)
By homogeneity, $\norm{P_m(e)}_\infty\le A\norm e_\infty^m$. The remainder in \eqref{eq:local-error-equation} has norm at most $\norm e_\infty^m$ in a sufficiently small neighborhood. Thus $\norm{E(e)}_\infty\le(A+1)\norm e_\infty^m$ there. If $P_m\not\equiv0$, choose $v$ with $\norm v_\infty=1$ and $P_m(v)\ne0$. For $t>0$,
\[
 E(tv)=t^mP_m(v)+o(t^m),\qquad
 \frac{\norm{E(tv)}_\infty}{t^m}\longrightarrow\norm{P_m(v)}_\infty>0.
\]
For every $r>m$, the quotient $\norm{E(tv)}_\infty/t^r$ tends to $+\infty$. This rules out a local estimate of exponent $r$.

Finally, choose $0<\delta_1<\delta$ so that $C\delta_1^{p-1}\le1/2$. If $\norm{e_0}_\infty<\delta_1$, then, as long as $e_k\ne0$,
\(\norm{e_{k+1}}_\infty \le C\norm{e_k}_\infty^{p-1}\norm{e_k}_\infty \le\tfrac12\norm{e_k}_\infty.\)
Induction shows that all iterates remain in the neighborhood and $\norm{e_k}_\infty\le2^{-k}\norm{e_0}_\infty$. If some iterate is zero, all later iterates are zero under the chosen extension. Otherwise, the quotients are defined and the initial estimate gives the stated limit superior. Furthermore,
\(\frac{\norm{e_{k+1}}_\infty}{\norm{e_k}_\infty} \le C\norm{e_k}_\infty^{p-1}\longrightarrow0.\)
This proves all three assertions.
\end{proof}

\begin{remark}[One-cycle limits and historical formulations]\label{rem:historical-order}
In \cite[Section 2, p.~644]{KT1974}, order is formulated by a limit as the input $x$ approaches the root, not merely by a limit along a trajectory. Every finite limit of the normalized error on a neighborhood gives the estimate of the first part of the proposition. Similarly, a vector error equation of the form \eqref{eq:local-error-equation}, valid as $e\to0$ without restricting directions, gives the corresponding estimate. Thus, when a formulation of order satisfies these conditions and the method belongs to the defined query model, the barriers apply directly; no equivalence with trajectory order needs to be established.

This remark does not, by itself, identify all historical classes of methods with our model. Causality, absence of memory, accounting of information, and universality over the function class must be preserved. Nor does it replace a neighborhood limit with an equation verified only along a particular sequence.
\end{remark}

\subsubsection{Trajectories and limits of the reverse implication}
\begin{definition}[Trajectory and exact $Q$ order]\label{def:Q-order}
Let $E$ be a map defined near zero with $E(0)=0$. A trajectory is a sequence satisfying $e_{k+1}=E(e_k)$ whose iterates lie in the domain. It is nonterminating if $e_k\ne0$ for every $k$. A nonterminating trajectory converging to zero has \emph{exact $Q$ order $r>1$} if there exists $\ell\in(0,+\infty)$ such that
\(\frac{\norm{e_{k+1}}_\infty}{\norm{e_k}_\infty^r}\longrightarrow\ell.\)
We say it satisfies a \emph{trajectory estimate with exponent $r$} if there exist $C>0$ and $k_0$ such that $\norm{e_{k+1}}_\infty\le C\norm{e_k}_\infty^r$ for $k\ge k_0$. Exact $Q$ order implies this latter property, but the converse does not require a finite positive limit.
\end{definition}

The constant in a trajectory estimate may depend on the initial point. It must not be confused with a constant valid for every input in a neighborhood. Nor do we assign infinite exact $Q$ order to a trajectory reaching zero: thereafter the quotients in the definition have zero denominators.

\begin{example}[$C^1$ regularity and different orders]\label{ex:scalar-Q}
Consider the scalar map in \cite[Example 10.7]{CastroKT},
\[
 \psi(e)=\begin{cases}-e^2,&e\ge0,\\-e^4,&e<0.\end{cases}
\]
The function is of class $C^1$: it is continuous at zero, its derivative there is zero, and the derivatives $-2e$ and $-4e^3$ on both sides tend to zero. If $0<|e_0|<1$, the first error is negative and nonzero. Thereafter,
\(e_{k+1}=-e_k^4,\qquad |e_{k+1}|=|e_k|^4.\)
All such trajectories converge to zero and have exact $Q$ order four, with constant one. However, for $e>0$,
\(\frac{|\psi(e)|}{|e|^4}=e^{-2}\longrightarrow+\infty.\)
The neighborhood estimate of order four fails. The exact local order of this map is two: $|\psi(e)|\le|e|^2$ for $|e|<1$, and the equality on the positive half-line rules out every higher exponent.
\end{example}

\begin{example}[Vector analyticity is not enough]\label{ex:vector-Q}
Consider the polynomial map
\(E:\R^2\longrightarrow\R^2,\qquad E(u,v)=(v^2,v^3).\)
For $\norm{(u,v)}_\infty<1$, we have $\norm{E(u,v)}_\infty=|v|^2\le\norm{(u,v)}_\infty^2$. On $(u,v)=(0,t)$, with $0<t<1$, the error norm is exactly $t^2$. Hence the exact local order of $E$ is two.

On the other hand, if $0<|v_0|<1$, the trajectory from $(u_0,v_0)$ satisfies, for $k\ge1$,
\[
 v_k=v_0^{3^k},\qquad u_k=v_0^{2\cdot3^{k-1}},\qquad
 \norm{(u_k,v_k)}_\infty=|v_0|^{2\cdot3^{k-1}}.
\]
These identities hold for $k=1$ and are preserved by substitution in $E$. Therefore,
\(\norm{(u_{k+1},v_{k+1})}_\infty =\norm{(u_k,v_k)}_\infty^3\qquad(k\ge1).\)
Every sufficiently near nonterminating trajectory has exact $Q$ order three, with constant one. If $v_0=0$, the trajectory reaches zero in one step.

Thus, even for an analytic output map, the order of all nonterminating trajectories does not automatically determine the local order of the vector map. The lowest-degree homogeneous term, $(v^2,0)$, vanishes in nonzero directions; the iterates approach those directions. Both examples in this subsection are abstract error maps, not query methods on an unknown function. They do not refute any barrier: they refute the logical implication that would allow us to switch notions of order without additional hypotheses.
\end{example}

\subsubsection{A sufficient vector criterion for trajectory order}
The preceding examples show that a leading term can vanish along nonzero directions and allow trajectories of higher order. We study a sufficient condition excluding this phenomenon. It is not among the hypotheses of the universal local-order barriers.
Let $E$ be of class $C^{b+1}$ near zero, with $E(0)=0$. For $1\le j\le b$, write
\(P_j(e)=\frac1{j!}D^jE(0)[e,\ldots,e].\)
$D^jE(0)$ is a symmetric $j$-linear map. We say that the first nonzero homogeneous term up to degree $b$ is \emph{nondegenerate in directions} if, for the smallest $m\le b$ with $P_m\not\equiv0$,
\(P_m(v)\ne0\qquad\text{for every }v\ne0.\)
If $P_j\equiv0$ ($1\le j\le b$), the condition is vacuous.

\begin{lemma}[A trajectory and a nondegenerate leading term]\label{lem:vector-Q-bridge}
Let $b\in\N$, and suppose $E$ satisfies the preceding conditions. If there is a nonterminating trajectory $e_k\to0$ with a trajectory estimate of exponent $r>b$, then all $P_j$, $1\le j\le b$, are zero. Moreover,
\[
 \norm{E(e)}_\infty=O(\norm e_\infty^{b+1})\qquad(e\to0).
\]
\end{lemma}
\begin{proof}
If there is a first nonzero $P_m$, the hypothesis and compactness of the unit sphere give
\(c=\min_{\norm v_\infty=1}\norm{P_m(v)}_\infty>0.\)
By Taylor's theorem and homogeneity,
\[
 E(e)=P_m(e)+o(\norm e_\infty^m),\qquad
 \norm{P_m(e)}_\infty\ge c\norm e_\infty^m.
\]
Thus $\norm{E(e)}_\infty\ge(c/2)\norm e_\infty^m$ for sufficiently small $e\ne0$. Along the trajectory,
\[
 \frac{\norm{e_{k+1}}_\infty}{\norm{e_k}_\infty^r}
 \ge\frac c2\norm{e_k}_\infty^{m-r}\longrightarrow+\infty,
\]
contradicting the trajectory estimate. Therefore all homogeneous polynomials up to degree $b$ vanish.

For each small $e$, apply Taylor's theorem with integral remainder to the one-variable function $s\mapsto E(se)$ at $[0,1]$. We obtain
\(E(e)=\frac1{b!}\int_0^1(1-s)^b D^{b+1}E(se)[e,\ldots,e]\,\dd s.\)
On a sufficiently small closed ball, continuity gives a bound $M$ for the norm of $D^{b+1}E$. Hence
\(\norm{E(e)}_\infty\le\frac{M}{(b+1)!}\norm e_\infty^{b+1}.\)
This proves the assertion.
\end{proof}

\begin{corollary}[$Q$-order barrier under nondegeneracy]\label{cor:Q-vector}
Let $\M$ be a method with local budget $(N_F,N_J)$ on $\Hclass(I_0)$, with $d=N_F+N_J\ge1$. Consider $b=2^{d-1}$. Alternatively, if $N_F,N_J\ge1$ and the method is locally paired, we may take $b=B(N_F,N_J)$. Suppose that, for every problem of $\Hclass(I_0)$, the vector error map admits an extension $C^{b+1}$ at zero and its first nonzero homogeneous term up to degree $b$, when one exists, is nondegenerate in directions.

Suppose further that there exists $p>1$ such that, for every problem, either the error map vanishes in a neighborhood or there is a convergent nonterminating trajectory toward the root with exact $Q$ order $p_F\ge p$. Then $p\le b$.
\end{corollary}
\begin{proof}
If $p>b$, each problem not solved exactly in one cycle satisfies the hypotheses of Lemma~\ref{lem:vector-Q-bridge}, with $r=p_F$. We obtain a local estimate of exponent $b+1$. The same estimate is immediate for problems whose error map vanishes in a neighborhood. Thus the method has universal local order at least $b+1$ on $\Hclass(I_0)$, contradicting Corollary~\ref{cor:CT} or Theorem~\ref{thm:paired-general}, according to the choice of $b$.
\end{proof}

\begin{remark}[Scope of the trajectory connection]
Nondegeneracy is imposed on $E$, not on $F$; Example~\ref{ex:vector-Q} shows it does not follow from analyticity. In one dimension it holds automatically when there is a first nonzero homogeneous term. We do not claim that this condition is necessary or that every historically considered method satisfies it. Without it, the barriers proved here remain barriers for the universal local order of one cycle.
\end{remark}

\clearpage
\section{Conclusions}\label{sec:conclusions}
We have proved a barrier for the universal local order of deterministic, causal, memoryless methods that solve nonlinear systems through function and Jacobian queries. Under the specified function classes and well-definedness conditions, we obtain
\[
 p\le 2^{N_F+N_J-1},\qquad N_F,N_J\in\N_0,\quad N_F+N_J\ge1.
\]
We require neither pairing, regularity of the rules, nor a prescribed initial order of queries. The Cordero--Torregrosa inequality holds without $1\le N_J\le N_F$, including when $N_J>N_F$. This extension concerns the range of budgets for which the upper bound holds; it does not assert attainability of the bound for each pair of budgets.

The general proof is self-contained. In the scalar case, we first construct one fixed real entire function, a simple root, and one sequence of inputs along which the error quotients diverge simultaneously for every exponent above the barrier. The passage from scale-dependent perturbations to a fixed witness preserves exactly the responses to all queries made. Thus we need neither continuity nor computability of the rules to compare outputs. Simulation on diagonal fields transfers this obstruction to systems without increasing the information budget.

Our third result is a refinement of the barrier under pairing of function and Jacobian nodes. For $N_F,N_J\ge1$, with $s=\min\{N_F,N_J\}$, we prove
\[
 p\le 2^{N_F-s+1}3^{s-1}.
\]
The bound remains valid when pairing is required only at the end of the cycle, without temporal precedence. Causal normalization queries and stores the required function values in advance but delivers to the original rules only the information requested at each stage. Hence it reproduces their decisions and output without increasing local budgets.

The construction is related to the methodology we developed earlier for the Kung--Traub problem \cite{CastroKT}. We prove the case of queries of orders zero and one within this article and separately explain how the earlier scalar result also yields the general barrier. The citation shows methodological continuity but is not a premise of the proof presented here.

Finally, our examples show why the assumptions on function classes cannot be omitted: on all polynomial fields, a noncomputable algebraic rule can achieve local exactness with few queries. Likewise, the universal order of one cycle does not by itself determine the exact $Q$ order of all trajectories. Open problems include general attainability of the bounds, restrictions intermediate between arbitrary and regular rules, and extensions of these arguments to other types of information. The resulting bounds are upper bounds subject to the explicit hypotheses of our theorems.

\clearpage
\phantomsection


\addcontentsline{toc}{section}{References}

\begin{thebibliography}{99}

\bibitem{CastroKT}
R. A. Castro Marín,
\emph{Order barriers under query budgets:
the Kung--Traub conjecture as a corollary},
preprint, version 1.0, Zenodo,
November 2, 2026.
\url{https://doi.org/10.5281/zenodo.23873507}.

\bibitem{KT1974}
H. T. Kung and J. F. Traub,
\emph{Optimal Order of One-Point
and Multipoint Iteration},
Journal of the ACM
\textbf{21}(4) (1974), 643--651.
\url{https://doi.org/10.1145/321850.321860}.

\bibitem{ACT2011}
V. Arroyo, A. Cordero and J. R. Torregrosa,
\emph{Approximation of artificial satellites'
preliminary orbits: The efficiency challenge},
Mathematical and Computer Modelling
\textbf{54}(7--8) (2011), 1802--1807.
\url{https://doi.org/10.1016/j.mcm.2010.11.063}.

\bibitem{ACT2010}
V. Arroyo, A. Cordero and J. R. Torregrosa,
\emph{Optimal iterative methods to calculate
artificial satellites preliminary orbits},
in L. Jódar (ed.),
\emph{Modelling for Addictive Behaviour,
Medicine and Engineering 2010},
Multidisciplinary Mathematics Institute,
Polytechnic University of Valencia,
2010, pp. 13--16.
\url{https://imm.webs.upv.es/jornadas/proceedings/Modelling2010.pdf}.


\bibitem{Woz1975}
H. Wo\'{z}niakowski,
\emph{Generalized Information and Maximal Order of Iteration for Operator Equations},
SIAM Journal on Numerical Analysis \textbf{12}(1) (1975), 121--135.
\url{https://doi.org/10.1137/0712012}.

\bibitem{Woz1976}
H. Wo\'{z}niakowski,
\emph{Maximal Order of Multipoint Iterations Using $n$ Evaluations},
in \emph{Analytic Computational Complexity}, Academic Press, 1976,
pp.~75--107.
\url{https://doi.org/10.1016/B978-0-12-697560-4.50010-1}.

\bibitem{Wasilkowski1981}
G. W. Wasilkowski,
\emph{$n$-Evaluation Conjecture for Multipoint Iterations
for the Solution of Scalar Nonlinear Equations},
Journal of the ACM (1981).
\url{https://doi.org/10.1145/322234.322240}.

\bibitem{TraubWoz1980}
J. F. Traub and H. Wo\'{z}niakowski,
\emph{A General Theory of Optimal Algorithms},
Academic Press, 1980.

\end{thebibliography}
\end{document}